\documentclass[a4paper,11pt]{article}
\usepackage[top=25mm,left=22mm]{geometry}
\pdfoutput=1 
\usepackage{amsmath,amssymb,amsthm,amsfonts}
\usepackage{lineno,url}\usepackage{float}
\floatstyle{plaintop}\restylefloat{table}
\usepackage[tableposition=top]{caption}
\usepackage{graphicx}
\usepackage{paralist, longtable}
\newtheorem{theorem}{Theorem}[section]
\newtheorem{remark}[theorem]{Remark}
\def\bexe{\begin{exercise}}\def\eexe{\eex\end{exercise}}
\def\bsol{\begin{solution}}\def\esol{\eex\end{solution}}
\def\bexa{\begin{example}}\def\eexa{\end{example}}
\def\brem{\begin{remark}}\def\erem{\end{remark}}
\def\bthm{\begin{theorem}}\def\ethm{\end{theorem}}
\def\blem{\begin{lemma}}\def\elem{\end{lemma}}
\def\bcor{\begin{corollary}}\def\ecor{\end{corollary}}
\def\bdefi{\begin{definition}}\def\edefi{\end{definition}}
\newcommand{\IDEA}{\textbf{Idea of the Proof.} }

\def\bmip{\begin{minipage}{\textwidth}}\def\emip{\end{minipage}}
\def\huga#1{\begin{gather} #1 \end{gather}}
\def\hugast#1{\begin{gather*} #1 \end{gather*}}
\def\hual#1{\begin{align} #1 \end{align}}

\newcommand{\R}{{\mathbb R}}
\newcommand{\C}{{\mathbb C}}\newcommand{\N}{{\mathbb N}}
\newcommand{\Q}{{\mathbb Q}}\newcommand{\Z}{{\mathbb Z}}

\def\CA{{\cal A}}\def\CD{{\cal D}}  
\def\CG{{\cal G}}\def\CH{{\cal H}}
\def\CO{{\cal O}}
\def\CM{{\cal M}}

\def\setm{\setminus}

\def\uti{\tilde{u}}
\def\yt{\tilde{y}}

\def\ga{\gamma}\def\om{\omega}
\def\noi{\noindent}\def\ds{\displaystyle}
\def\vt{\vartheta}\def\pa{{\partial}}\def\lam{\lambda}
\newcommand{\bi}{\begin{itemize}}\newcommand{\ei}{\end{itemize}}
\newcommand{\ben}{\begin{enumerate}}\newcommand{\een}{\end{enumerate}}
\newcommand{\bce}{\begin{center}}\newcommand{\ece}{\end{center}}
\newcommand{\bci}{\begin{compactitem}}\newcommand{\eci}{\end{compactitem}}
\newcommand{\bcen}{\begin{compactenum}}\newcommand{\ecen}{\end{compactenum}}
\newcommand{\bcena}{\begin{compactenum}[(a)]}
\newcommand{\reff}[1]{(\ref{#1})}
\newcommand{\ov}[1]{{\overline {#1}}}

\newcommand{\spr}[1]{\left\langle #1 \right\rangle}
\newcommand{\hs}[1]{{\hspace{#1}}}\newcommand{\vs}[1]{{\vspace{#1}}}

\def\eps{\varepsilon}
\def\ra{\rightarrow}

\newcommand{\barr}{\begin{array}}\newcommand{\earr}{\end{array}}
\newcommand{\bpm}{\begin{pmatrix}}\newcommand{\epm}{\end{pmatrix}}
\newcommand{\bsm}{\left(\begin{smallmatrix}}
\newcommand{\esm}{\end{smallmatrix}\right)}
\newcommand{\ba}{\begin{array}}\newcommand{\ea}{\end{array}}
\def\dd{\, {\rm d}}\def\ri{{\rm i}}

\def\er{{\rm e}}
\def\re{{\rm Re}}\def\im{{\rm Im}}
\def\om{\omega}\def\Om{\Omega}
\def\ddt{\frac{\rm d}{{\rm d}t}}
\def\del{\delta}

\def\eex{\hfill\mbox{$\rfloor$}}

\def\phiti{\tilde{\phi}}

\def\al{\alpha}

\def\bd{\begin{displaymath}} \def\ed{\end{displaymath}}
\def\ba{\begin{array}} \def\ea{\end{array}}
  
\def\eps{\varepsilon}

\newcommand{\argmax}{\operatornamewithlimits{argmax}}

\def\rds{{\rm ds}}
\def\ig{\includegraphics}\def\pdep{{\tt pde2path}}\def\pdepo{{\tt pde2path}}
\def\sign{{\rm sign}}\def\Re{{\rm Re}}

\def\oop{{\tt OOPDE}}
\def\mlab{{\tt Matlab}}\def\ptool{{\tt pdetoolbox}}

\def\bcs{BC}

\def\ind{{\rm ind}}\def\emu{{\rm err}_{\ga_1}}
\renewcommand{\arraystretch}{1}\renewcommand{\baselinestretch}{1.0}
\def\medskip{}\def\bigskip{}

\def\Ati{{\tilde A}}\def\Bti{{\tilde B}}
\def\ati{{\tilde a}}\def\bti{{\tilde b}}\def\hoxi{\xi_{\text{H}}}
\def\neig{n_\text{eig}}\def\fltol{{\rm tol}_{{\rm fl}}}
\def\heda{{\bf HD1}}\def\hedb{{\bf HD2}}
\def\fla{{\bf FA1}}\def\flb{{\bf FA2}}

\usepackage{hyperref}
\usepackage{blindtext}

\begin{document}
\text{}\vspace{5mm}\begin{center}\Large
Hopf bifurcation and time periodic orbits with \pdep\ -- algorithms and applications \\[4mm]
\normalsize Hannes Uecker \\[2mm]
\small 
Institut f\"ur Mathematik, Universit\"at Oldenburg, D26111 Oldenburg, 
hannes.uecker@uni-oldenburg.de\\[2mm]
\normalsize
\today
\end{center}
\noi
\begin{abstract} 
We describe the algorithms used in the \mlab\ continuation and 
bifurcation package \pdep\ for Hopf bifurcation and 
continuation of branches of periodic orbits in systems of PDEs in 1, 2, and 3 spatial dimensions, including the computation of 
Floquet multipliers. 
We first test the methods on three reaction diffusion examples, namely 
a complex Ginzburg--Landau equation as a toy problem, a reaction diffusion 
system on a disk with rotational waves including 
stable (anti) spirals bifurcating out of the trivial solution, 
and a Brusselator system with interaction 
of Turing and Turing--Hopf bifurcations. Then we consider a system 
from distributed optimal control, which is ill-posed as an 
initial value problem and thus needs a particularly stable 
method for computing Floquet multipliers, for which we use 
a periodic Schur decomposition. The implementation details how to use 
\pdep\ on these problems are given in an accompanying tutorial, which, 
together with all other downloads (function libraries, demos and 
further documentation) can be found at \url{www.staff.uni-oldenburg.de/hannes.uecker/pde2path}. 
\end{abstract}

\tableofcontents 

\section{Introduction}\label{i-sec}
The package \pdepo\ \cite{p2pure, p2p2, p2phome} 
has originally been developed as a continuation/bifurcation 
package for stationary problems of the form 
\huga{\label{gform}
G(u,\lam):=-\nabla\cdot(c\otimes\nabla u)+a u-b\otimes\nabla u-f=0. 
}
Here $u=u(x)\in\R^N$, $x\in\Omega$ with $\Om\subset\R^d$ some bounded domain, 
$d=1,2,3$, 
$\lam\in\R^p$ is a parameter (vector), 
and $c\in\R^{N\times N\times 2\times 2}$,
$b\in\R^{N\times N\times 2}$, 
 $a\in\R^{N\times N}$ and $f\in\R^N$
 can  depend on $x,u,\nabla u$, and 
parameters.%
\footnote{We have 
$[\nabla\cdot(c\otimes\nabla u)]_i:=
\sum_{j=1}^N [\pa_xc_{ij11}\pa_x+\pa_x c_{ij12}\pa_y+\pa_yc_{ij21}\pa_x
+\pa_yc_{ij22}\pa_y]u_j$ ($i^{{\rm th}}$ component), and similarly 
$[au]_i=\sum_{j=1}^N a_{ij}u_j$, 
$[b\otimes\nabla u]_i:=\sum_{j=1}^N [b_{ij1}\pa_x+b_{ij2}\pa_y]u_j$, 
and $f=(f_1,\ldots,f_N)$ as a column vector. }
The  boundary conditions (BC) are of 
``generalized Neumann'' form, i.e., 
\begin{align}\label{gnbc}
{\bf n}\cdot (c \otimes\nabla u) + q u = g,
\end{align}
where ${\bf n}$ is the outer normal and again $q\in \R^{N\times N}$
and $g\in \R^N$ may depend on $x$, $u$, $\nabla u$ and
parameters. These BC include zero flux \bcs, 
and a ``stiff spring'' approximation of Dirichlet \bcs\ via 
large prefactors in $q$ and $g$, and  periodic BC are also supported
over suitable domains. 
Moreover, there are interfaces to couple \reff{gform} with 
additional equations, such as mass conservation, or phase 
conditions for considering co-moving frames, and to set up extended 
systems, for instance for fold point and branch point continuation. 

\pdep\ has been applied to various research problems, e.g., 
patterns in 2D reaction diffusion systems 
\cite{uwsnak14, k15, k15b, sder15,  w16, ZUFM17},  
some problems in fluid dynamics and nonlinear optics 
\cite{ZHR15,DU14,ET17} and in optimal control \cite{U16, GU17}. 
Here we report on features and algorithms in \pdepo\ to treat Hopf (or Poincar\'e--Andronov--Hopf) bifurcations and the continuation  
of time--periodic orbits for systems of the form 
\huga{\label{tform} 
\pa_t u=-G(u,\lam), \quad u=u(x,t),\ x\in\Om\subset\R^d,\ d=1,2,3,\  
t\in\R\ \ 
\text{ ($d+1$ dimensional problem),}
}
with $G$ from \reff{gform} and BC from \reff{gnbc}. 
Adding the time dimension makes computations more expensive, 
such that here we focus on 1D and 2D, and only give one 3D 
example to illustrate that all user interfaces are essentially dimension 
independent. 

For general introductions to and reviews of 
(numerical) continuation and bifurcation we recommend 
\cite{gov2000, kuz04, doedel07, seydel}, and \cite{mei2000}, 
which has a focus on reaction--diffusion systems.
The treatment of large scale problems, typically 
from the spatial discretization of PDEs, including the continuation of time periodic orbits, has for instance 
been discussed in \cite{LRSC98, TBark00, LR00}, and has recently been 
reviewed in \cite{Detal14}. There, 
the focus has been on matrix--free methods where 
the periodic orbits are computed by a shooting method, which can 
conveniently be implemented if a time--stepper for the given 
problem is available. In many cases, 
shooting methods can also be used 
to investigate the bifurcations from periodic orbits, 
and to trace bifurcation curves in parameter space, 
by computing the Floquet multipliers of the periodic orbits. 
In this direction, see in particular 
\cite{ BTh00, SanGM2013, WIJ13, NSan2015} for impressive results in 
fluid problems.  

Here we proceed by a collocation (in time) method for the continuation 
of periodic orbits. With respect to computation time and in particular 
memory requirements such methods 
are often more demanding than (matrix free) shooting methods. However, 
one reason for the efficiency of shooting methods in the works cited above is 
that there the problems considered are strongly dissipative, with 
only few eigenvalues of the linearized evolution near the imaginary axis. 
We also treat such problems, and show that up to moderately large scale 
they can efficiently be treated by collocation methods as well. However, 
another class of problems we deal with are 
canonical systems obtained from distributed 
optimal control problems with infinite time horizons. 
Such problems are ill-posed as initial value problems, 
which seems quite problematic for genuine shooting methods. 

We also compute the Floquet multipliers for periodic orbits. 
For this, a direct approach is to explicitly construct 
the monodromy matrix 
from the Jacobian used in the collocation solver 
for the periodic orbit. We find that this works well for dissipative 
problems, but 
completely fails for the ill--posed optimal control problems, and thus we also 
provide a method based on a periodic Schur decomposition, which can 
handle this situation. Currently, our Floquet computations 
can be used to assess the stability of periodic orbits, and for 
{\em detection} of possible bifurcations from them. However, we do not (yet) 
provide tools for {\em localization} of, or branch switching at, 
such bifurcation points,  which is work 
in progress.

To illustrate the performance of our {\tt hopf} library  we consider four 
example problems, 
with the \mlab\ files included as demo directories in the package download at \cite{p2phome}, where also a \pdep\ user-guide with installation instruction, 
a tutorial on Hopf bifurcations with implementation details, and various other tutorials on how to run \pdep\ are available.  
The first example is a cubic--quintic  complex Ginzburg--Landau (cGL) equation,  
which we consider over 1D, 2D, and 3D cuboids with homogeneous 
Neumann and Dirichlet BC, such that we can explicitly  
calculate all Hopf bifurcation points (HBP) from the trivial 
branch. 
For periodic BC we also have the Hopf branches 
explicitly, which altogether 
makes the cGL equation a nice toy problem to validate 
and benchmark our routines. Next we consider 
a reaction diffusion system from \cite{GKS00}  
on a circular domain and with somewhat non-standard 
Robin BC, which lead to rotating waves, and in particular to 
the bifurcation of (anti) spiral waves out of the trivial solution. 
Our third example is a 
Brusselator system from \cite{yd02}, which shows interesting 
interactions between Turing branches and Turing--Hopf branches.  
As a non--dissipative example we treat the canonical system 
for a simple control problem of ``optimal pollution''. This is still of 
the form \reff{tform}, but is ill--posed as an initial value problem, 
since it includes ``backward diffusion''. Nevertheless, we continue 
steady states, and obtain Hopf bifurcations 
and branches of periodic orbits.

Many of the numerical results on periodic orbits in PDE in the literature, 
again see \cite{Detal14} for a review, are based on custom made 
codes, 
which sometimes do not seem easy to access and 
modify for non--expert users. 
Although in some of our research applications we consider 
problems with on the order of $10^5$ unknowns in space, 
\pdep\ is not primarily intended for very large scale problems. 
Instead, the goal of \pdep\ is to provide 
a general and easy to use (and modify and extend) toolbox to 
investigate bifurcations in PDEs of the (rather large) class 
given by \reff{tform}. 
With the {\tt hopf} library we provide some basic functionality for Hopf 
bifurcations and continuation of periodic orbits for such PDEs  over 1D, 2D, and 
3D domains, where at least the 1D cases and simple 2D cases are sufficiently 
fast to use \pdep\ as a quick (i.e., interactive) tool 
for studying interesting problems. 
The user interfaces reuse the standard \pdep\ setup, 
and no new user functions are necessary. 
Due to higher computational costs 
in 2+1D, in 3D, or even 3+1D, compared to the 2D case from \cite{p2pure}, 
in the applications given here and the associated tutorials 
we work with quite coarse meshes, but give a number 
of comments on how to adaptively generate and work with finer meshes. 

In \S\ref{nsec} we review some basics of the Hopf bifurcation, 
of periodic orbit continuation and multiplier computations, and 
explain their numerical treatment in \pdep.  
In \S\ref{exsec} we present the examples, and \S\ref{dsec} contains 
a brief summary and outlook. The \pdep\ setup, 
data structures and help system are reviewed in \cite{qsrcb}, and 
implementation details for the Hopf demos are given in \cite{hotutb}.  
For comments, questions, and bugs, please mail to 
{\tt hannes.uecker@uni-oldenburg.de}. 

\vs{3mm}
\noindent 
{\bf Acknowledgment.}  Many thanks to Francesca Mazzia 
for providing TOM \cite{MT04}, which was 
essential help for setting up the {\tt hopf} library; to Uwe Pr\"ufert 
for providing \oop; 
to Tomas Dohnal, Jens Rademacher and Daniel Wetzel for some 
testing of the Hopf examples; to Daniel Kressner for {\tt pqzschur}; 
to Arnd Scheel for helpful comments on the system in \S\ref{rotsec}; 
and to Dieter Grass for the cooperation 
on distributed optimal control problems, which was one of my main 
motivations to implement 
the {\tt hopf} library. 

\section{Hopf bifurcation and periodic orbit continuation in \pdep}
\label{nsec}
Our description of the algorithms is 
based on the spatial FEM discretization of \reff{tform}, which, with a slight 
abuse of notation, we write as 
\huga{\label{tform2} 
M\dot u(t)=-G(u(t),\lam), 
}
where $M\in\R^{n_u\times n_u}$ is the mass matrix, $n_u=Nn_p$ is the 
number of unknowns (degrees of freedom DoF) with $n_p$ is the number 
of mesh-points,  and, for each $t$, 
$$
u(t)=(u_{1,1},\ldots,u_{1,n_p}, u_{2,1},\ldots,u_{N,1},\ldots u_{N,n_p})(t)
\in \R^{n_u},$$ and similarly $G:\R^{n_u}\times\R^{p}\ra \R^{n_u}$. 
We use the 
generic name $\lam$ for the parameter vector, {\em and} the 
{\em active} continuation parameter, again 
see \cite{p2p2} for details. When in the following we discuss eigenvalues $\mu$ 
and eigenvectors $\phi$ of the linearization 
\huga{\label{lin1}
M\dot v=-\pa_u G(u,\lam)v
}
of \reff{tform2} around some (stationary) solution of \reff{tform2}, 
or simply eigenvalues of $\pa_u G=\pa_u G(u,\lam)$, 
we always mean the generalized eigenvalue problem 
\huga{\label{evp1}
 \mu M\phi=\pa_u G\phi.
}
Thus eigenvalues of $\pa_u G$ with {\em negative} real parts 
give dynamical {\em instability} of $u$. 

\brem\label{qrem}{\rm  For, e.g., the continuation of 
traveling waves in translationally invariant problems, 
the PDE \reff{tform} is typically transformed 
to a moving frame $\xi=x-\ga(t)$, with BC that respect the 
translational invariance, and where $\dot \ga$ is an unknown 
wave speed, which yields an additional term $\dot\ga\pa_x u$ on the rhs of 
\reff{tform}. The reliable continuation of traveling then also 
requires a phase condition, i.e., an additional equation, 
for instance of the form $q(u)=\spr{\pa_x \uti, u}\stackrel!=0$, where 
$\uti$ is a reference wave (e.g.~$\uti=u_{{\rm old}}$, where 
$u_{{\rm old}}$ is from a previous continuation step), 
and $\spr{u,v}=\int_\Om uv\dd x$.  
Together we obtain a differential--algebraic system instead 
of \reff{tform2}, and similarly for other constraints such 
as mass conservation, see \cite[\S2.4,\S2.5]{p2p2} for examples, 
and for instance \cite{BT07, symtut} for equations with 
continuous symmetries and the associated ``freezing method''. 
Hopf bifurcations can occur in such systems, see e.g.~the Hopf 
bifurcations from traveling ($\dot \ga\ne 0$) or standing ($\dot \ga=0$) 
fronts and pulses in \cite{hm94-pattern, GAP06, BT07, GF2013}, 
but are somewhat more difficult to treat numerically than the case without 
constraints. 
Thus, here we restrict to %
problems of the form \reff{tform} without constraints, 
and hence to \reff{tform2} on the spatially discretized level, and refer 
to \cite{symtut,hotut} for examples of Hopf bifurcations with 
constraints in \pdep. For instance, in \cite[\S4]{symtut} we 
consider Hopf bifurcations to modulated traveling waves in a 
model for autocatalysis, and the Hopf bifurcation of standing 
breathers in a FitzHugh--Nagumo system, and in \cite[\S5]{hotut} 
the Hopf bifurcation of modulated standing and traveling waves 
in the Kuramoto-Sivashinky equation with periodic boundary conditions. 
}\eex \erem

\subsection{Branch and Hopf point detection and localization} 
Hopf bifurcation means the bifurcation of a branch of time periodic 
orbits from a branch $\lam\mapsto u(\cdot,\lam)$ of 
stationary solutions of \reff{tform}, or numerically \reff{tform2}. 
This generically  occurs if at some $\lam=\lam_H$ 
a pair of simple complex conjugate 
eigenvalues $\mu_j(\lam)=\ov{\mu}_{j+1}(\lam)$ of %
$G_u=\pa_u G(u,\lam)$ crosses the imaginary axis with 
nonzero imaginary part and nonzero speed, i.e., 
\huga{\label{hocon}
\mu_j(\lam_H)=\ov{\mu}_j(\lam_H)=\ri \om\ne 0, 
\quad \text{ and }\re\mu_j'(\lam_H)\ne 0.
}
Thus, the first issue is to define a suitable test function $\psi_H$ to 
numerically detect \reff{hocon}. Additionally, we also want to detect real 
eigenvalues crossing the imaginary axis, i.e., 
\huga{\label{bcon}
\mu_j(\lam_{\text{BP}})=0, 
\quad \text{ and }\re\mu_j'(\lam_{\text{BP}})\ne 0.
}
A fast and simple method for \reff{bcon} is to monitor sign changes of 
the product 
\huga{\label{bdp1} 
\psi(\lam)=\prod_{i=1,\ldots,n_u}\mu_i(\lam)=\det(G_u)}
of all eigenvalues, which even for large 
$n_u$ can be done quickly via the $LU$ factorization of $G_u$, 
respectively the extended matrix in arclength continuation, 
see \cite[\S2.1]{p2pure}.  This so far has been the standard setting 
in \pdepo, but the 
drawback of \reff{bdp1} is that the sign of $\psi$ only changes 
if an odd number of real eigenvalues crosses $0$. 

Unfortunately, there is no general method for \reff{hocon} 
which can be used for large $n_u$. For small systems, one option would be 
\huga{\label{psiH}
\psi_H(\lam)=\prod_{i}(\mu_i(\lam)+\mu_{i+1}(\lam)), 
}
where we assume the eigenvalues to be sorted by their real parts. 
However, this, unlike \reff{bdp1} requires the actual computation 
of all eigenvalues, which is not feasible for large $n_u$. 
Another option are dyadic products, \cite[Chapter 10]{kuz04}, 
which again is not feasible for large $n_u$. 

If, on the other hand, \reff{tform} is a dissipative problem, 
 then we may try to just compute 
$n_\text{eig}$ eigenvalues of $G_u$ of smallest modulus, which, 
for moderate $n_\text{eig}$ can be done efficiently, and to 
count the number of these eigenvalues which are 
in the left complex half plane, 
and from this detect both \reff{hocon} and \reff{bcon}. 
 The main issue then is to choose $n_\text{eig}$, which 
unfortunately is highly problem dependent, and for a 
given problem may need to be chosen large again.

The method presented in \cite{GoS96} 
uses complex analysis, namely the winding number 
$W(g(i\om),0,\infty)$ of $g(z)=c^T(G_u-zM)^{-1}b$, which is 
the Schur complement of the bordered system 
$\bpm G_u-zM&b\\c^T&0\epm$ with (some choices of) $b,c\in\R^{n_u}$. 
We have 
\huga{\label{gzdef}
g(z)=\frac{N(z)}{\det(G_u-zM)}, \text{ where 
$W(g(i\om),0\infty)=\pi(Z_l-Z_r+P_r-P_l)/2$}, 
} 
where $Z_{l,r}, P_{l,r}$ 
are the zeros and poles of $g(z)$ in the left and right complex 
half planes, respectively, and 
where $N$ is a polynomial in $z$ which depends on $b,c$. 
 Since $\det(G_u-zM)$ does not depend on 
$b,c$, using some clever evaluation  \cite{GoS96} of \reff{gzdef} for some 
choices of $b,c$ one can count the poles of $g$, 
i.e.~the eigenvalues of $G_u$  in the left half plane. 

Here we combine the 
idea of counting small eigenvalues with 
suitable spectral shifts $\ri\om_{1,2,\ldots}$. To estimate 
these shifts, given a current solution $(u,\lam)$ we follow \cite{GoS96} and 
compute 
\huga{\label{reso}
[0,\om_{\max} ]\ni\om\mapsto g(u,\lam,\ri \om; b):=b^T(G_u-\ri\om)^{-1}b, 
}
for one or several suitably chosen $b\in\R^{n_u}$. Generically, $g$ 
will be large for $\ri\om$ near some complex eigenvalue 
$\mu=\mu_r+\ri\mu_i$ with small $\mu_r$, and thus we may consider this $\ri\om$ 
as a {\em guess} 
for a Hopf eigenvalue during the next continuation steps. 
To accurately compute $g$ from \reff{reso} we again use ideas 
from  \cite{GoS96} to refine the 
$\om$ discretization (and actually 
compute the winding number). Then, after 
each continuation step we compute a few eigenvalues near  
$0,\om_1,\ldots$. 
We can reset the shifts $\om_i$ after a number 
of continuation steps by evaluating \reff{reso} again, and 
instead of using \reff{reso} the user can also 
set the $\om_{i}$ by hand.%
\footnote{In principle, instead of using \reff{reso} 
we could also compute the guesses $\om_i$ 
by computing eigenvalues of $G_u(u,\lam)$ at a given $(u,\lam)$; 
however, this may itself either require a priori information on 
the pertinent $\om_i$ (for shifting), or we may again 
need to compute a large number of eigenvalues of $G_u$. Thus 
we find \reff{reso} more simple, efficient and elegant.}

\begin{figure}[ht]
\bce{\small
\begin{tabular}{llll}
(a) $n=0, n_d=0$& (b) $n=1, n_d=1$ & (c) $n=1, n_d=0$& (d) $n=2, n_d=0$.\\
\ig[width=0.23\textwidth]{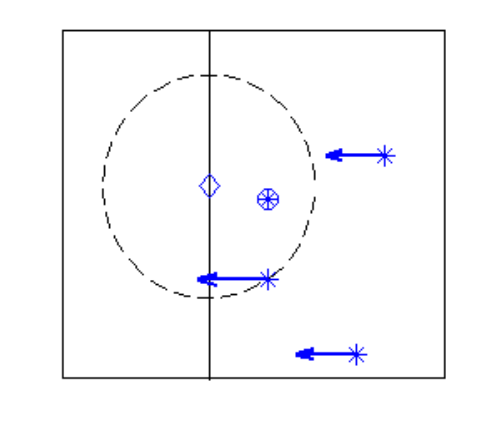}&
\ig[width=0.23\textwidth]{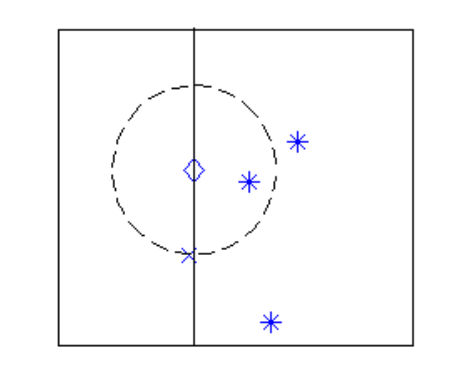}&
\ig[width=0.23\textwidth]{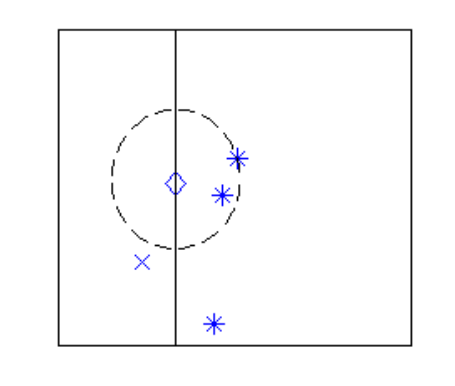}&
\ig[width=0.23\textwidth]{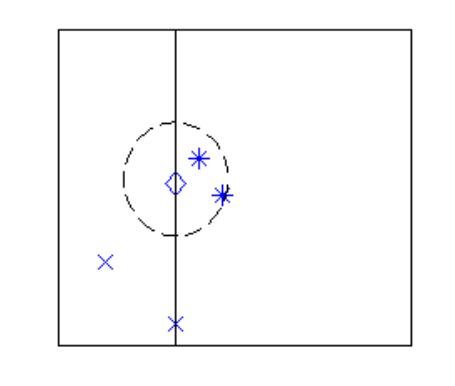}
\end{tabular}}
\ece 

\vs{-5mm}
   \caption{{\small Sketch of the idea, and typical failures, 
of detecting Hopf points 
by counting eigenvalues with negative real parts near some shift $\ri\om_1$, 
marked by $\diamond$. Here, for illustration we use {\tt neig=2}, 
i.e., use the 2 eigenvalues closest to $\ri\om_1$ for 
bifurcation detection, and show 4 eigenvalues 
near $\ri\om_1$, stable ones with $*$ and unstable ones with $\times$. 
$n$ is the total number of negative eigenvalues, 
and $n_d$ the number of detected ones. 
From (a) to (d) we assume that some parameter $\lam$ 
varies, and the shown eigenvalues depend continuously on $\lam$; 
for better illustration we assume that the eigenvalue circled 
in (a) stays fixed. 
The dashed circle has radius $|\mu(\lam)-\ri\om_1|$ with $\mu$ the 
second closest eigenvalue to $\ri\om_1$.  
From (a) to (b) we correctly detect a HBP. From (b) to (c) we incorrectly 
find a HBP, as the only unstable eigenvalue wanders out of the pertinent 
circle. From (c) to (d) we miss a HBP, as the guess $\ri\om_1$ is too 
far off. The failure (b) to (c) can be detected in the localization 
by requiring that at the end the real part of the 
eigenvalue closest to the imaginary axis is sufficiently small. 
The failure from (c) to (d) can be resolved by either (i) computing 
more eigenvalues close to $\ri\om_1$, or (ii) by updating $\ri\om_1$ using 
\reff{reso}.    \label{fbpf}}}
\end{figure}

Of course, the idea is mainly heuristic, and, in this simple 
form, may miss some bifurcation points (BPs, in the sense of 
\reff{bcon}) and 
Hopf bifurcation points (HBPs, in the sense of \reff{hocon}), 
and can and typically 
will detect false BPs and HBPs, see Fig.~\ref{fbpf}, 
which illustrates two ways in which the algorithm can fail.%
\footnote{A third typical kind of failure is that 
during a continuation 
step a number $m$ of eigenvalues  crosses 
the imaginary axis close to $\ri \om_1$, and simultaneously $m$ already 
unstable eigenvalues leave the pertinent circle to the left due to a 
decreasing real part. The only remedy for this is to decrease the step--length 
ds. Clearly, a too large ds can miss bifurcations even if we could 
compute {\em all} eigenvalues, for instance if along the true branch 
eigenvalues cross back and forth. }
However, some of these 
failures can be detected and/or corrected, see Remark \ref{murem}, 
and in practice we found the algorithm to work 
remarkably well in our examples, with a rather 
small numbers of eigenvalues computed near $0$ and $\ri \om_1$, 
and in general to be more robust 
than the theoretically more sound usage of \reff{gzdef}.%
\footnote{However, if additionally to bifurcations one is interested 
in the stability of (stationary) solutions, then the numbers 
of eigenvalues should not be chosen too small; otherwise 
the situation in Fig.~\ref{fbpf}(c,d) may easily occur, i.e., 
undetected eigenvalues with negative real parts.} 
For convenience, in the following we refer to these algorithms as 
\bci 
\item[\heda\ (Hopf Detection 1)] \label{hed} compute the $\neig$ smallest 
eigenvalues of $\pa_u G$ and count those with negative real parts; 
\item[\hedb\ (Hopf Detection 2)] initially compute a number 
of guesses $\om_j$, $j=1,\ldots, g$ for spectral shifts, then compute 
the $n_{\text{eig},j}$ eigenvalues of $\pa_u G$ closest to $\om_j$, and count 
how many have negative real parts to detect crossings of eigenvalues 
near $\om_j$. Update the shifts when appropriate. 
\eci

After detection of a (possible) BP or a (possible) HBP, 
or of several of these along a branch 
between $s_0$ and $s_0+\rds$, it remains to locate the BP or HBP. 
Again, there are various methods to do this, using, e.g., suitably 
extended systems \cite{gov2000}. However, so far we typically 
use a simple bisection, which works well and sufficiently fast in our 
examples.%
\footnote{The only extended systems we deal with in \pdep\ so far are 
those for localization and continuation of stationary branch points, and 
of fold points, see \cite[\S2.1.4]{p2p2}or \cite{lsstutb}.}

\brem\label{murem}{\rm To avoid unnecessary bisections and 
false BPs and HBPs we proceed as follows. After detection of a BP or HBP 
{\em candidate} with shift $\om_j$, we check if the 
eigenvalue $\mu$ closest to $\ri\om_j$ has $|\re\mu|\le \mu_1$, 
with default $\mu_1=0.01$.  If not, 
then we assume that this was a false alarm. Similarly, {\em after completing} 
a bisection we check if the eigenvalue $\mu$ {\em then} closest to $\om_j$ has 
$|\re\mu|<\mu_2$, with default $\mu_2=0.0001$, 
and only then accept the computed point as a BP (if $\om_j=0$) 
or HBP (if $\om_j>0$). In our examples, 
about 50\% of the candidates enter the bisection, 
and of these about 10\% are rejected afterwards, and no 
false BPs or HBPs are saved to disk. 
This seems to be a reasonable compromise between 
speed and not missing BPs and HBPs and avoiding false ones. 
However, the values of $\mu_1,\mu_2$ are of course highly problem dependent. }
\eex
\erem

\subsection{Branch switching}
Branch switching at a BP works as usual by computing an initial guess 
from the normal form of the stationary bifurcation, see \cite[\S2.1]{p2pure}. 
Similarly, to switch to a Hopf branch of time periodic solutions 
we compute an initial guess from an approximation of the normal form 
\huga{\label{honf}
\dot{w}=\mu(\lam)w+c_1(\lam)|w|^2 w, 
} 
of the bifurcation equation on the center manifold associated 
to $(\lam,\mu)=(\lam_H,\ri\om_H)$. Thus we use 
\huga{
\mu(\lam)=\mu_r(\lam)+\ri \mu_i(\lam)=\mu_r'(\lam_H)(\lam-\lam_H)
+\ri(\om_H+\CO(\lam-\lam_H))+\CO((\lam-\lam_H)^2), 
}
and with $w=r\er^{\ri \om_H t}$ replace \reff{honf} by 
\huga{\label{honf2}
0=r\bigg[\mu_r'(\lam_H)(\lam-\lam_H)+c_1(\lam_H)|r|^2\bigg]. 
} 
Following \cite{kuz04}, 
$c_1=c_1(\lam_H)\in\R$ is related to the first Lyapunov coefficient 
$l_1$  by $c_1(\lam_H)=\om_H l_1$, and we use the 
formulas from \cite[p531-536]{kuz04} for the numerical 
computation of $l_1$. Setting $\lam=s\eps^2$ with $s=\pm 1$ 
we then have a nontrivial solution 
\huga{
r=\eps\al, \quad 
\al=\sqrt{-s\mu'(\lam_H)/c_1(\lam_H)}
}
of \reff{honf2} for $s=-\sign(\mu'(\lam_H)/c_1)$, and thus take 
\huga{\label{hoswitr} 
\lam=\lam_H+s\eps^2,
\quad u(t)=u_0+2\eps\al\Re(\er^{-\ri \om_H t}\Psi), }
as an initial guess for a periodic solution of \reff{tform2} 
with period near $2\pi/\om$. 
The approximation \reff{hoswitr} 
of the bifurcating solution in the center eigenspace, also called 
linear predictor, is usually 
accurate enough, and is the standard setting in the \pdep\ 
function {\tt hoswibra}, see \cite{hotutb}. 
The coefficients $s=\pm 1$ and $\al$ in \reff{hoswitr} are computed, 
and $\eps$ is then chosen in such a way 
that the initial step length is {\tt ds} in the norm 
\reff{xinorm} below. %

\subsection{The continuation of branches of periodic orbits}\label{hocontsec}
\subsubsection{General setting}
The continuation of the Hopf branch is, as usual, a predictor--corrector method, 
and for the corrector we offer, essentially, two different methods, namely 
natural and arclength  continuation. 
For both, we reuse the standard \pdep\ settings for assembling $G$ in 
\reff{tform} and Jacobians, such that the user does not have 
to provide new functions. 
In any case, 
first we rescale $t=Tt$ in \reff{tform2} to obtain 
\huga{\label{tform3} 
M\dot u=-TG(u,\lam), \quad u(\cdot,0)=u(\cdot,1), 
}
with unknown period $T$, but with initial guess $T=2\pi/\om$ at bifurcation.

\subsubsection{Arclength parametrization}\label{alsec}
We start with the arclength setting, which is more general and more robust, 
although the continuation in natural parametrization in \pdep\ has other advantages such as 
error control and adaptive mesh refinement for the time discretization, 
see below. We use the 
phase condition 
\huga{\label{pca}
\phi:=\int_0^1\spr{u(t),M\dot u_0(t)}\dd t\stackrel !=0, 
}
where $\spr{\cdot,\cdot}$ is the scalar product in $\R^{n_u}$ and 
$\dot u_0(t)$ is from the previous continuation step, 
and we use the step length condition 
\huga{\label{ala}\psi:=
\hoxi\sum_{j=1}^m\spr{u(t_j){-}u_0(t_j),u_0'(t_j)}+
(1{-}\hoxi)\bigl[w_T(T{-}T_0)T_0'+(1{-}w_T)(\lam{-}\lam_0)\lam_0'\bigr]{-}\rds
\stackrel!=0,} 
where again $T_0,\lam_0$ are from the previous step, $\rds$ is the 
step--length,  $'=\frac{\dd}{\dd s}$ denotes differentiation 
with respect to arclength, $\hoxi$ and $w_T$ denote 
weights for the $u$ and $T$ components of the unknown solution, 
and $t_0=0<t_1<\ldots<t_m=1$ is the temporal discretization. 
Thus, the step length is $\rds$ in the weighted norm 
\huga{\label{xinorm}
\|(u,T,\lam)\|_\xi=\sqrt{\hoxi\left(\sum_{j=1}^m \|u(t_j)\|_2^2\right)+(1-\hoxi)
\bigl[w_TT^2+
(1-w_T)\lam^2\bigr]}. 
}

Even if $\hoxi$ is similar to the (average) mesh--width in $t$, then 
the term $\hoxi\sum_j \|u(t_j)\|_2$ is only a crude approximation of the 
``natural length'' $\int_0^1 \|u(t)\|_2\dd t$. However, the choice 
of the norm is somewhat arbitrary, and we found \reff{xinorm} 
most convenient. 
Typically we choose $w_T=1/2$ such that $T$ and $\lam$ have the same weight 
in the arclength. A possible choice for $\hoxi$ to weight 
 the $m n_u$ components of $u$ is  
\huga{\label{hoxichoi}
\hoxi=\frac 1 {m n_u}. 
}
However, in practice we choose $\hoxi=\frac{10}{m n_u}$, or even 
larger (by another factor 10), since at the Hopf bifurcation 
the branches are ``vertical'' ($\|u-u_0\|=\CO(\sqrt{|\lam-\lam_0|})$, 
cf.~\reff{hoswitr}), and $\hoxi$ tunes the search direction 
in the extended Newton loop  between 
``horizontal'' (large $\hoxi$) and ``vertical'' (small $\hoxi$).  
See \cite[\S2.1]{p2pure} for the analogous 
role of $\xi$ for stationary problems. 

The integral in \reff{pca} is discretized 
as a simple Riemann sum, such that the derivative of $\phi$ 
with respect to $u$ is, with $\tilde u_0(t)=M\dot u_0(t)$, 
\def\uti{\tilde u}
\hual{\label{pauphi}
\pa_u \phi&=(h_1 \uti(t_1)_1,\ldots,h_1\uti(t_1)_{n_u}, 
h_2\uti(t_2)_1,\ldots,h_2\uti(t_2)_{n_u},\ \ldots, 
h_{l-1}\uti(t_{m-1})_{n_u},0,\ldots,0), 
}
$n_u$ zeros at the end, 
where $h_l=t_{l+1}{-}t_l$ is the mesh--size in the time discretization. 
Similarly, denoting the tangent along the branch as 
\huga{\label{taudef}
\tau=(\tau_u,\tau_T,\tau_\lam),\quad \tau_u\in\R^{1\times 
mn_u} \text{ (row vector as in \reff{pauphi}),}\quad \tau_T,\ \tau_\lam\in\R,}
we can rewrite $\psi$ in \reff{ala} as 
\huga{\label{ala2} 
\psi=\hoxi \tau_u(u-u_0)+(1-\hoxi)(w_T \tau_T(T-T_0)+(1-w_T)\tau_\lam(\lam-\lam_0))-\rds. 
}

Setting $U=(u,T,\lam)$, and writing 
\reff{tform3} as $\CG(u,T,\lam)=0$, in each continuation step 
we thus need to solve 
\huga{\label{fsa} 
H(U):=\bpm \CG(U)\\\phi(u)\\\psi(U)\epm\stackrel{!}{=}\bpm 0\\0\\0\epm
\in\R^{mn_u+2}, }
for which we use Newton's method, i.e., 
\huga{\label{news} 
U^{j+1}{=}U^{j}{-}\CA(U^j)^{-1}H(U^j), \quad 
\CA{=}\bpm \pa_u \CG&\pa_T\CG&\pa_\lam\CG\\
\pa_u\phi&0&0\\
\hoxi\tau_u&(1{-}\hoxi)w_T\tau_T&(1{-}\hoxi)(1{-}w_T)\tau_\lam\epm.
}
After convergence of $U^j$ to $U$, i.e., $\|H(U)\|\le $''tolerance'' in some suitable 
norm, the next tangent $\tau_1$ 
with preserved orientation $\spr{\tau_0,\tau_1}>0$ can be calculated as 
usual from 
\huga{\label{newtan}
\CA(U)\tau_1=(0,0,1)^T. 
}

It remains to discretize in time and 
assemble $\CG$ in \reff{tform3} and the Jacobian 
$\pa_u \CG$. For this we use (modifications of) routines from TOM \cite{MT04}, 
which assumes the unknowns in the form 
\huga{\label{ydef0}
u=(u_1,\ldots,u_m)=(u(t_1), u(t_2),\ldots, u(t_m)),
\quad\text{($m$ time slices),}
} 
Then, using the TOM $k=1$ setting, we have, for $j=1,\ldots,m-1$, 
the implicit backwards in time finite differences 
\huga{\label{gjdef}
(\CG(u))_j=-h_{j-1}^{-1}M(u_j-u_{j-1})-\frac 1 2 T(G(u_j)+G(u_{j-1})), 
}
where $u_0:=u_{m-1}$, and additionally the periodicity condition 
\huga{
G_m(u)=u_m-u_1.}
The Jacobian is $\pa_u\CG=A_1$, where  we set, 
as it reappears below for the Floquet multipliers, 
\huga{\label{fl0}
A_\ga=
\bpm M_1&0&0&0&\ldots&-H_1&0\\
-H_2&M_2&0&0&\ldots&0&0\\
0&-H_3&M_3&0&\ldots&0&0\\
\vdots&\ldots&\ddots&\ddots&\ddots&\vdots&\vdots\\
0&\ldots&\ldots&\ddots&\ddots&0&0\\
0&\ldots&\ldots&0&-H_{m-1}&M_{m-1}&0\\
-\ga\,I&0&\ldots&\ldots&\ldots&0&I
\epm,}
where $\ds M_j=-h_{j-1}^{-1}M-\frac 1 2 T G_u(u_j)$, 
$\ds H_j=-h_{j-1}^{-1}M+\frac 1 2 T G_u(u_{j-1})$, and $I$ is 
the ${n_u\times n_u}$ identity matrix. 
The derivatives $\pa_T \CG, \pa_\lam \CG$ in \reff{news} are cheap from 
numerical differentiation, and $\pa_u\phi$ and $\tau$ do not change 
during Newton loops, and are easily taken care of anyway. 

\brem\label{belrem}{\rm 
$\CA\in\R^{(mn_u+2)\times (mn_u+2)}$ in \reff{news}, \reff{newtan} 
consists of $A=A_1=\CG_u\in \R^{mn_u\times mn_u}$, which is large but sparse, 
and borders of widths 2, i.e., symbolically, 
$$
\CA=\bpm A&B\\
C&D\epm, \text{ with large and sparse $A$, with } C^T,B\in\R^{mn_u\times 2}, 
\text{ and }D\in\R^{2\times 2}.
$$
There are various methods to solve bordered systems of the form 
\huga{\label{bsys} 
\CA x=b, \qquad b=\bpm f\\g\epm,}
see, e.g., \cite{gov2000}. Here we use the very simple scheme 
\huga{\label{mbel}
V=A^{-1}B, x_1=A^{-1}f, \tilde{D}=D-CV, y_1=g-Cx_1, 
y_2=\tilde{D}^{-1}y_1, x_2=x_1-Vy_2, x=\bpm x_2\\ y_2\epm.
}
The big advantage of such bordered schemes is that solving systems 
such as $Ax_1=f$ (where we either pre-factor $A$ for repeated 
solves, or use a preconditioned iterative method) 
is usually much cheaper due to the structure of $A$ than solving 
$\CA x=b$ (either by factoring $\CA$, or by an iterative method 
with some preconditioning of $\CA$).%
\footnote{The special structure of $A$ from \reff{fl0} can also 
be exploited to solve $Ax=f$ in such a way that subsequently the 
Floquet multipliers can easily be computed,  see \S\ref{fremsec}, 
and \cite{lust01} for comments on the related algorithms used in AUTO. }

The scheme \reff{mbel} may suffer from some instabilities, 
but often these can be corrected 
by a simple iteration: If $\|r\|$ with $r=\CA x-b$ is too large, 
then we solve $\CA\hat x=r$ (again by \reff{mbel}, which is cheap) 
and update $x=x-\hat{x}$, until $\|r\|\le$ ``tolerance''.  
We in particular sometimes obtain poor solutions of \reff{bsys} 
for $b=(0,0,1)^T$ from \reff{newtan}, but they typically can 
be improved by a few iterations. Altogether we 
found the scheme \reff{mbel} to work well in our problems, 
with a typical speedup of up to 50 compared to the 
direct solution of $\CA x=b$. 
Again, see \cite{gov2000} for 
alternative schemes and detailed discussion. 

For the solutions of $AV=B$ and $Ax_1=f$ in \reff{mbel} we 
give the option to use a preconditioned iterative solver 
from {\tt ilupack} \cite{ilupack}, see also \cite{lsstutb}.%
\footnote{In fact, when using iterative solvers it is often advantageous 
to directly use them for the full system \reff{bsys}, since iterative 
solvers seem rather indifferent to the borders.} 
The number of continuation steps before a new preconditioner is 
needed can be quite large, and often the iterative solvers 
give a  significant speedup. 
}\eex
\erem 

\subsubsection{Natural parametrization}\label{natpasec}
By keeping $\lam$ fixed during correction we cannot pass around folds, 
but on the other hand can take advantage of further useful features of TOM. 
TOM requires 
separated boundary conditions, and thus we use a standard trick and introduce, 
in the notation \reff{ydef0}, 
auxiliary variables 
\def\yt{\tilde u}
$\yt=(\yt_1, \yt_2,\ldots, \yt_m)$ and additional (dummy) ODEs 
$\dot \yt_l=0$. Then setting the boundary conditions 
\huga{\label{artbc}
u_1-\yt_1=0, \quad u_m-\yt_m=0} 
corresponds to periodic boundary conditions for $u$. 
Moreover, we add the auxiliary equation $\dot T=0$, and set up 
the phase condition 
\huga{\label{pcnat}
\phi=\spr{u(0),M\dot u_0(0)}\stackrel!=0. 
} 
as an additional boundary condition. Thus, the complete system to be solved 
is 
\huga{\label{natsys}
\bpm M\dot u\\ \dot \yt\\ \dot T\epm=\bpm -TG(u)\\0\\0\epm,}
together with \reff{artbc} and \reff{pcnat}. 
We may then pass 
an initial guess (from a predictor) at a new $\lam$ to TOM, and let 
TOM solve for $(u,\yt)$ and $T$.  The main advantage is 
that this comes with error control and adaptive  mesh 
refinement for the temporal discretization.%
\footnote{We however also provide an ad-hoc mesh-refinement routine 
for the arclength case, see \cite{hotutb}.}

\subsection{Floquet multipliers}\label{fremsec}
The (in)stability of -- and possible bifurcations from -- a periodic orbit 
$u_H$ are analyzed via the Floquet multipliers $\ga$. These are obtained 
 from 
finding nontrivial solutions $(v,\ga)$ of the variational boundary value 
problem 
\hual{
M\dot v(t)&=-T\pa_u G(u(t))v(t),\label{fl1}\\ 
v(1)&=\ga v(0).
}
Equivalently, the multipliers $\ga$ are the eigenvalues of 
the monodromy matrix $\CM(u_0)=\pa_u \Phi(u_0,T)$, where $\Phi(u_0,t)$ is 
the solution of the initial value problem \reff{tform2} with $u(0)=u_0$ 
from $u_H$. 
Thus, $\CM(u_0)$ depends on $u_0$, but the multipliers $\ga$ do not. 
 By translational invariance, there always is the 
trivial multiplier $\ga_1=1$. 
$\CM(u_0)$ is the linearization of the Poincar\'e map 
$\Pi(\cdot; u_0)$ around $u_0$, which maps a point 
\def\ut{\tilde{u}}
$\ut_0$ in a hyperplane 
$\Sigma$ through $u_0$ and transversal to $u_H$ to its first return 
to $\Sigma$. 
Therefore, a necessary conditions for the bifurcation 
from a branch $\lam\mapsto u_H(\cdot,\lam)$ of periodic orbits 
is that at some $(u_H(\cdot,\lam_0),\lam_0)$, 
additional to the trivial multiplier 
$\ga_1=1$ there is a 
second multiplier $\ga_2$ (or a complex conjugate pair $\ga_{2,3}$) 
with $|\ga_2|=1$, which generically leads to the following bifurcations 
(see, e.g., \cite[Chapter 7]{seydel} or \cite{kuz04} for more details): 
\bci
\item[(i)]\label{biftypes} $\ga_2=1$, yields a fold of the periodic orbit, or a transcritical 
or pitchfork bifurcation of periodic orbits.
\item[(ii)] $\ga_2=-1$, yields a period--doubling bifurcation, i.e., the bifurcation 
of periodic orbits $\ut(\cdot;\lam)$ with approximately double the period, 
$\ut(\tilde T;\lam)=\ut(0;\lam)$, $\tilde{T}(\lam)\approx 2T(\lam)$ for 
$\lam$ near $\lam_0$.  
\item[(iii)] $\ga_{2,3}=\er^{\pm\ri \vt}$ , $\vt\ne 0,\pi$, 
yields a torus (or Naimark--Sacker) bifurcation, i.e., the bifurcation 
of periodic orbits $\ut(\cdot,\lam)$ with two ``periods'' $T(\lam)$ and 
$\tilde T(\lam)$; 
if $T(\lam)/\tilde T(\lam)\not\in \Q$, then $\R\ni t\mapsto \ut(t)$ is dense in 
certain tori.  
\eci 

Here we are first of all interested in the computation of the multipliers. 
Using the same discretization for $v$ as for $u$, it follows that $\ga$ 
and $v=(v_1,\ldots,v_m)$ have to satisfy the matrix eigenvalue 
problem 
\huga{\label{fl2}
A_\ga v=0, 
}
where now $\ga$ in \reff{fl0} is free. 
From this special 
structure it is easy to see, that $\CM(u_{j_0})$ can be obtained 
from certain products involving the $M_j$ and the $H_j$, for instance 
\huga{\label{fl3} 
\CM(u_{m-1})=M_{m-1}^{-1}H_{m-1}\cdots M_1^{-1}H_1. 
} 
Thus, a simple way to compute the $\ga_j$ is to compute the 
product \reff{fl3} and subsequently  (a number of) 
the eigenvalues of $\CM(u_{m-1})$. We call this \fla\ (Floquet Algorithm 1), 
and using 
\huga{\label{emudef}
\emu:=|\ga_1-1|} 
as a measure of accuracy we find that this works 
fast and accurately for our dissipative examples. Typically 
$\emu<10^{-10}$, although at larger amplitudes of $u_H$, 
and if there are large multipliers, this may 
go down to $\emu\sim 10^{-8}$, which is the (default) tolerance we require 
for the computation of $u_H$ itself. Thus, in the software 
 we give a warning if $\emu$ exceeds a certain tolerance 
$\fltol$.   
However, for the optimal control example in \S\ref{ocsec}, 
where we naturally have multipliers $\ga_j$ with 
$|\ga_j|>10^{30}$ and larger\footnote{I.e., $|\ga_{n_u}|\ra\infty$ as 
$n_u\ra\infty$, although the orbits may still be stable in the 
sense of optimal control, see \S\ref{ocsec}}, 
\fla\ completely fails to compute any meaningful multipliers.

More generally, in for instance \cite{fj91,lust01} it is discussed that methods 
based directly on \reff{fl3} 
\bci
\item%
may give considerable numerical errors, in particular if 
there are both, very small and very large multipliers $\ga_j$; 
\item%
discard much useful information, for instance eigenvectors 
of $\CM(u_l)$, $l\ne m-1$, which are useful for branch switching.
\eci
As an alternative, \cite{lust01} suggests to use a periodic Schur decomposition 
\cite{BGD93} to compute the multipliers (and subsequently 
pertinent eigenvectors), and gives examples that in certain cases 
this gives much better accuracy, according to \reff{emudef}. 
See also \cite{kressner01, kress06} for similar ideas and results. 

We thus provide an algorithm \flb\  (Floquet Algorithm 2), 
which, based on {\tt pqzschur} from \cite{kressner01}, 
computes a periodic Schur decomposition of the matrices involved 
in \reff{fl3}, 
from which we immediately obtain the multipliers, see Remark \ref{floqrem}(d).  
For large $n_u$ and $m$, \flb\ gets rather slow, 
and thus we rather use it in two ways. First, 
to validate (by example) \fla, and second to compute 
the multipliers when \fla\ fails, in particular for 
our OC example. 

In summary,  for small to medium sized {\em dissipative} 
problems, i.e., $n_u*m<50000$, say, 
computing (a number of) multipliers with \fla\ 
is a matter of a seconds. 
For the {\em ill-posed} OC 
problems we have to use \flb\ which is slower and for 
medium sized problems can be as slow as the computation of $u_H$. 
In any case, because we do not 
yet consider the localization of the bifurcations (i)--(iii) 
from periodic orbits (this is work 
in progress%
), for efficiency we give the option to switch 
off the simultaneous computation of 
multipliers during continuation of periodic orbits. 

\brem\label{floqrem}{\rm 
(a) To save the stability information on the computed branch we define 
\huga{\label{inddef} 
\ind(u_H)=\text{number 
of multipliers $\ga$ with $|\ga|>1$ (numerically: $|\ga|>1+\fltol$),  
}
}
such that unstable orbits are characterized by $\ind(u_H)>0$.  
Thus, a change in $\ind(u_H)$ 
signals a possible bifurcation, and via 
\def\argmin{{\rm argmin}}
\def\gacand{\ga_{{\rm cand}}}
$$
\gacand:=\argmin\{|\ga_j|: |\ga_j|>1\} 
$$ 
we also get an approximation for the critical multiplier, and 
thus a classification of the possible bifurcation 
in the sense (i)-(iii). 

(b) In \fla\ we compute $n_+$  multipliers 
$\ga_2,\ldots,\ga_{n_+}$ of largest modulus 
(recall that 
we reserve $\ga_1$ for the trivial multiplier), with $|\ga_2|\ge |\ga_3|\ge 
\ldots\ge |\ga_{n_+}|$, and count how many of these have $|\ga_j|>1$, 
which gives $\ind(u_H)$ if we make sure that $|\ga_{n_+}|<1$. 
For dissipative systems, typically $0\le<n_+\ll n_u$, and the large multipliers 
of $\CM$ can be computed efficiently and reliably 
by vector iteration. However, it does happen that some of 
the small multipliers do not converge, in which case we also give a warning, 
and recommend to check the results with \flb.  

(c) The idea of the periodic Schur decomposition is as follows. Given two 
collections $(A_i),(B_i)$, $i=1,\ldots,m$, of matrices 
$A_i,B_i\in \C^{n\times n}$, {\tt pqzschur} computes $Q_i, Z_i, 
\Ati_i, \Bti_i\in \C^{n\times n}$ such that $\Ati_i,\Bti_i$ 
are upper triangular, $Q_i, Z_i$ are orthogonal, and 
\hugast{\barr{ll}
A_1=Q_1\Ati_1 Z_m^H,&B_1=Q_1\Bti_1 Z_1^H\\
A_2=Q_2\Ati_2 Z_1^H,&B_2=Q_2\Bti_1 Z_2^H\\
\ldots,&\ldots\\
A_m=Q_m\Ati_m Z_{m-1}^H,&B_m=Q_m\Bti_m Z_m^H.
\earr
}
Consequently, for the product $\CM=B_m^{-1}A_m\cdots B_1^{-1}A_1$ we 
have 
$$
\CM=Z_m \Bti_m^{-1}\Ati_m\cdots \Bti_1^{-1}\Ati_1 Z_m^H, 
$$ 
and similar for products with other orderings of the factors. 
In particular, the eigenvalues of $\CM$ are given by the products 
$\ds d_i=\prod_{j=1}^m\ati_{ii}^{(j)}/\bti_{ii}^{(j)}$, and, moreover, 
the associated eigenvectors can also be extracted from the $Q_j, Z_j$, 
see \cite{kress06} for further comments.

(d) Alternatively to using Floquet multipliers, we can assess the stability 
of the periodic orbits by using the time--integration routines from 
\pdep, which moreover has the advantage of giving information 
about the evolution of perturbations of unstable solutions; 
see \S\ref{exsec} for examples, where in all cases 
perturbations of unstable periodic orbits lead to convergence to some 
other (stable) periodic orbit. 
}
\eex\erem 

\section{Four examples}\label{exsec}
To illustrate the performance of our algorithms we use 
four examples, included as demos directories in the package download, together 
with the tutorial \cite{hotutb} explaining implementation details. 
See also \cite{p2phome} for a \pdep\ quickstart guide explaining the installation, data structures 
and help system of \pdep, and for other tutorials and further information. 

Thus, here we focus on explaining the results (i.e., the relevant plots), and on 
relating them to 
some mathematical background of the equations.  
In all examples, the meshes are chosen rather coarse, to quickly 
get familiar with the algorithms. 
We did check for all examples that these coarse meshes give reliable 
results by running the same simulations on finer meshes, without 
qualitative changes. 
In some problems we additionally 
switch off the simultaneous computation of Floquet multipliers, 
and instead compute the multipliers a posteriori 
at selected points on branches. 
 Nevertheless, even with the coarse meshes  
some commands, e.g., the continuation 
of Hopf branches in 3+1D, may take several minutes. 
All computational times given in the following 
are from an i5 laptop with Linux Mint 17 and \mlab\ 2013a. Using 
the {\tt ilupack} \cite{ilupack} iterative linear solvers, 
memory requirements are moderate ($<$ 2GB), but using direct solvers 
we need about 11GB for the largest scale problems considered here 
(3D cGL with about 90000 degrees of freedom, see \S\ref{cACsec}).

\subsection{A complex Ginzburg--Landau equation}\label{cACsec}
We consider the cubic-quintic complex Ginzburg--Landau equation 
\huga{\label{cAC0} 
\pa_t u=\Delta u+(r+\ri\nu)u-(c_3+\ri \mu)|u|^2 u-c_5|u|^4u, \quad 
u=u(t,x)\in\C, 
}
with real parameters $r,\nu,c_3,\mu,c_5$. Equations of this type 
are canonical models in physics, and are often derived as amplitude 
equations for more complicated pattern forming systems \cite{AK02,mie99}. 
Using real variables $u_1,u_2$ with $u=u_1+\ri u_2$, 
\reff{cAC0} can be written as a real 
2--component system of the form \reff{tform}, i.e., 
\huga{\label{cAC} 
\pa_t \bpm u_1\\ u_2\epm =\bpm \Delta+r&-\nu\\\nu&\Delta+r\epm
\bpm u_1\\ u_2\epm-(u_1^2+u_2^2)\bpm c_3 u_1-\mu u_2\\ 
\mu u_1+c_3 u_2\epm-c_5(u_1^2+u_2^2)^2\bpm u_1\\ u_2\epm. 
}

We set 
\huga{\label{cglpar}
\text{$c_3=-1, c_5=1, \nu=1, \mu=0.1$,}
}
and use $r$ as the main bifurcation parameter.  Considering 
\reff{cAC} on boxes 
\huga{\label{omdef}
\Om=(-l_1\pi,l_1\pi)\times\cdots\times(-l_d\pi,l_d\pi) 
} 
with homogeneous Dirichlet BC or Neumann BC, or with periodic BC, 
we can explicitly calculate all 
Hopf bifurcation points from the trivial branch $u=0$, and, for periodic BC, 
the bifurcating time periodic branches. For this let 
\huga{\label{cACa}\text{$u(x,t)=a\er^{\ri (\om t-k\cdot x)}$, with wave vector 
$k=(k_1,\ldots,k_d), k_j\in \frac{1}{2l_j}\Z $,} }
and temporal period $2\pi/\om$, which yields 
\hual{\label{cACs}
|a|^2{=}|a(k,r)|^2{=}-\frac{c_3}{2c_5}\pm\sqrt{\frac{c_3^2}{4c_5^2}+r-|k|^2}, 
\quad \om{=}\om(k,r){=}\nu-\mu|a|^2, \quad |k|^2{=}k_1^2+\ldots+k_d^2. 
}
Note that $\om$ and hence the period $T=2\pi/\om$ depend on $|a|$, that 
$u(t,x)$ on each branch is a single harmonic in $x$ and $t$, and 
that the phase of $a$ is free. Using \reff{cglpar} we obtain subcritical 
Hopf bifurcations of solutions \reff{cACa} at
\huga{\label{cACs2}
\text{$r=|k|^2$, with folds at $r=|k|^2-\frac 1 4$.}
} 
Moreover, for these orbits we can also compute the Floquet 
multipliers explicitly. For instance, restricting to $k=0$ in 
\reff{cACa}, and also to the invariant subspace of spatially 
independent perturbations, in polar-coordinates 
\def\ati{\tilde{a}}\def\phiti{\tilde\phi}
$\ut(t)=\ati(t)\er^{\ri\phiti(t)}$ we obtain the 
(here autonomous) linearized ODEs 
\huga{\label{exfl1}
\ddt{\ati}=h(r)\ati,\quad \ddt \phiti=-2\mu a\ati, 
\text{ where } h(r)=r+3a^2-5a^4. 
}
The solution is $\ati(T)=\er^{h(r)T}\ati(0)$, 
$\phiti(T)=\phiti(0)+\frac a {h(r)}(\er^{h(r)T}-1)\ati(0)$, and therefore 
the analytic monodromy matrix (in the $k{=}0$ subspace) is 
$\ds\CM_{k=0}=\bpm \er^{h(r)T}&0\\ \frac a {h(r)}(\er^{h(r)T}-1)&1\epm$ 
with multipliers $\ga_1{=}1$ and $\ga_2{=}\er^{h(r)T}$. 

Thus, \reff{cAC} makes a nice toy problem to validate 
and benchmark our routines, where,  
to avoid translational invariance, cf.~Remark \ref{qrem}, 
we use Neumann and Dirichlet BC. 
For these we still have the formula $r=|k|^2$ for the HBPs, although 
we lose the explicit branches, except the spatially homogeneous branch 
for $k=0$ with Neumann BC.

\begin{figure}[ht]
\bce{\small 
\begin{tabular}{ll}
(a) BD, norm $\|u(\cdot,\cdot;r)\|_*$&(b) Example plots\\
\ig[width=0.19\textwidth, height=38mm]{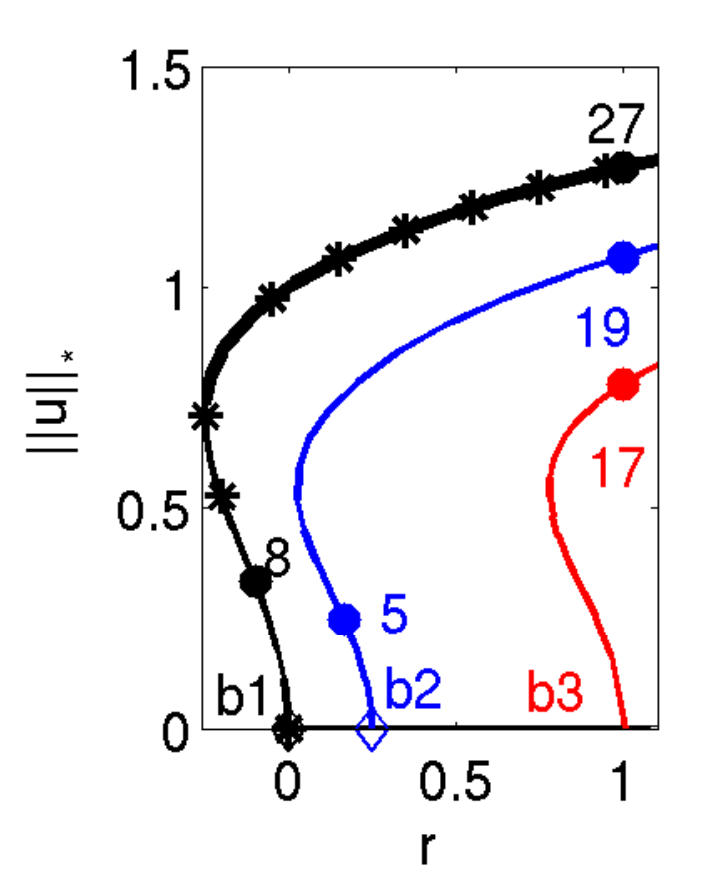}&
\ig[height=38mm]{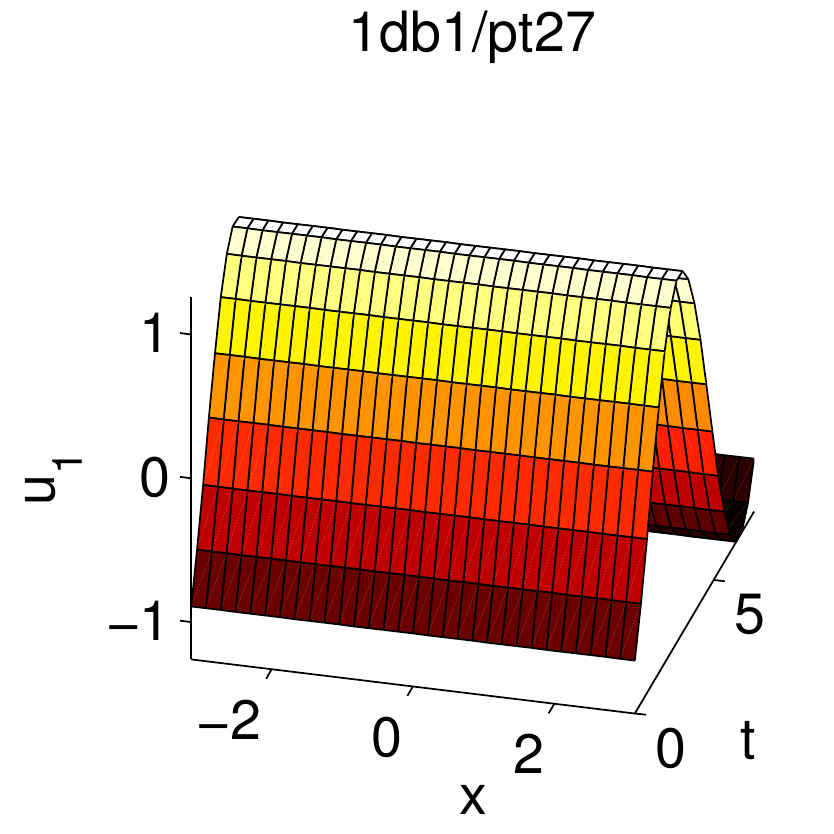}
\ig[width=0.24\textwidth]{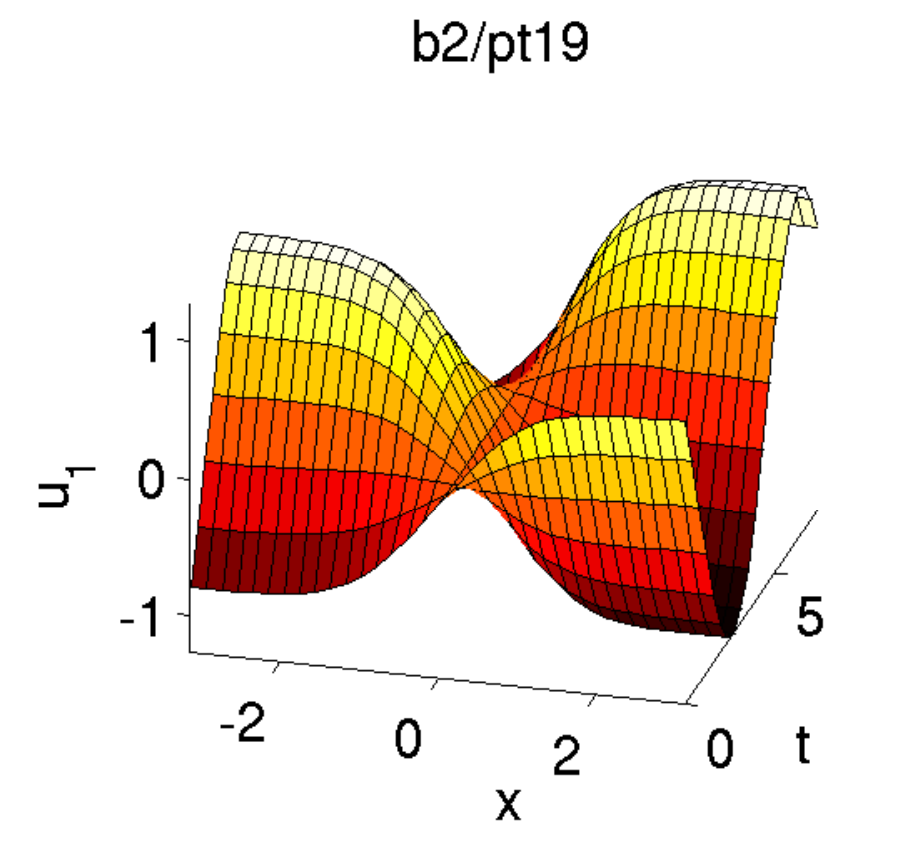}
\ig[width=0.24\textwidth]{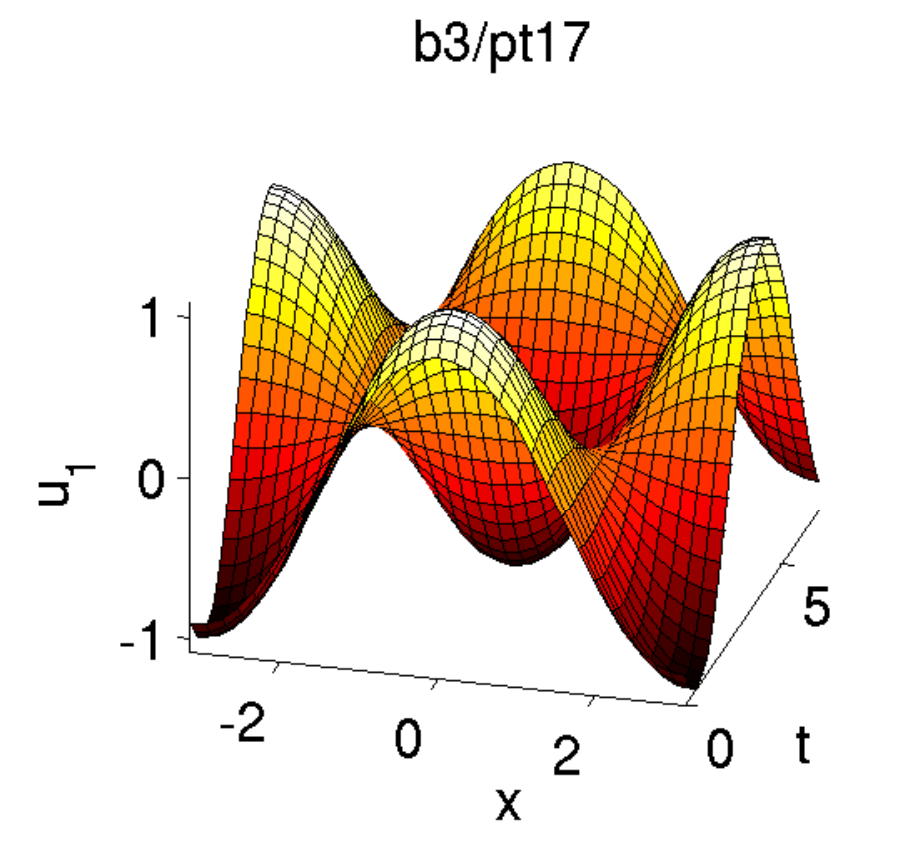} \\
(c) BD; period $T(r)$&
(d) Multipliers at b1/pt8 ($\ind=1$), b1/pt27 ($\ind=0$), 
and b2/pt5 ($\ind=3$)\\[1mm]
\ig[height=28mm]{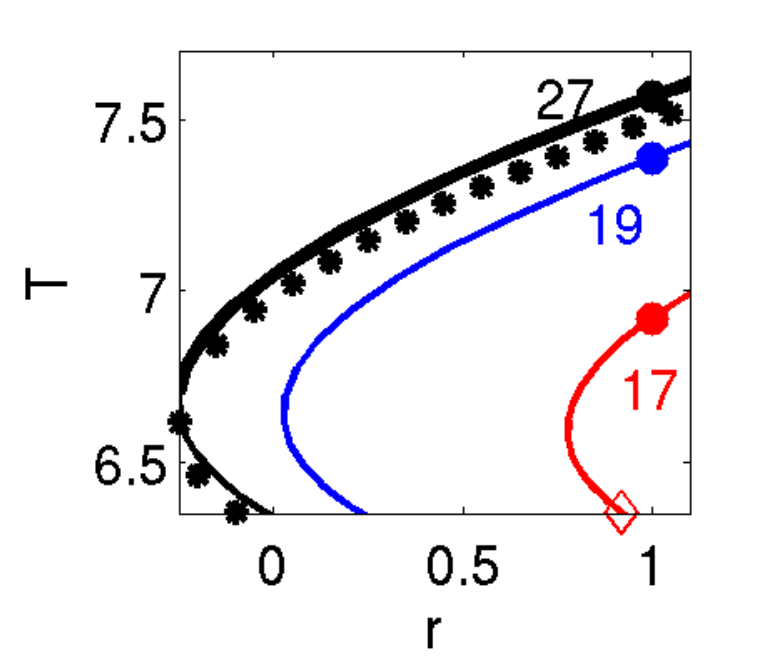}&
\ig[height=25mm]{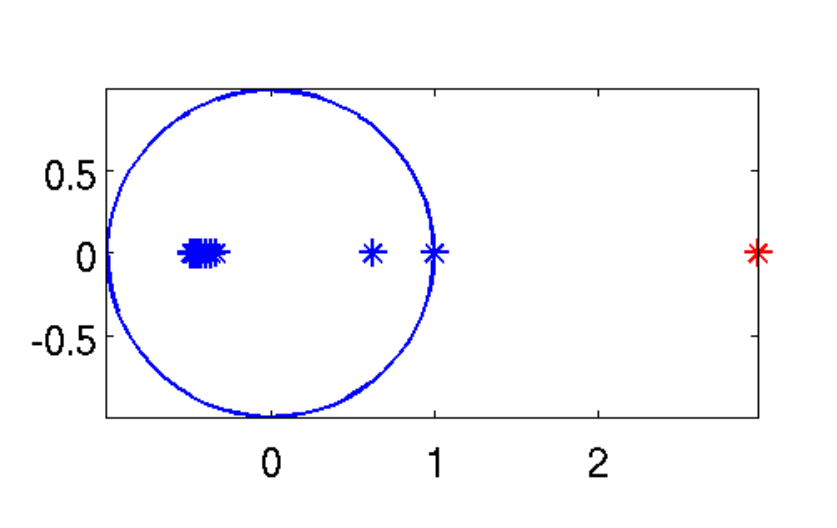}
\raisebox{3mm}{\ig[height=20mm]{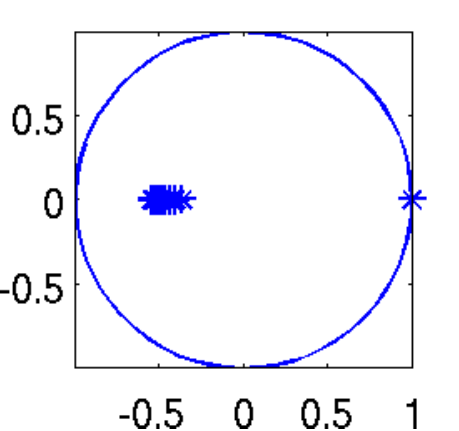}}
\raisebox{3mm}{\ig[height=25mm]{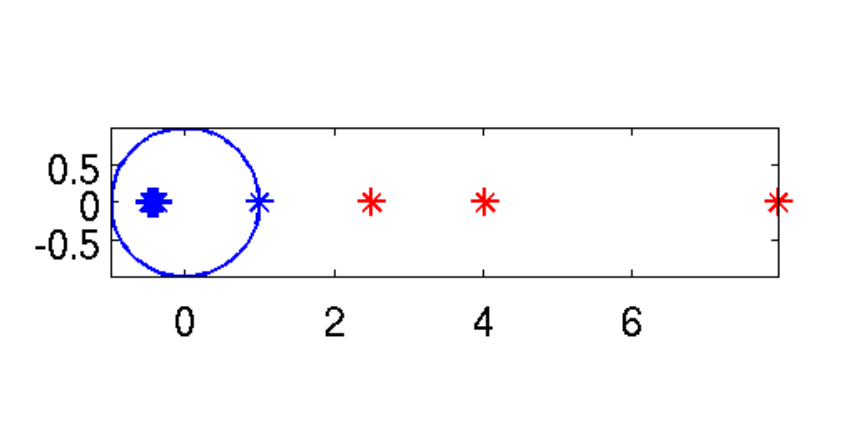}}
\end{tabular}
}
\ece 

\vs{-8mm}
   \caption{{\small Numerical bifurcation diagrams, example plots 
and (leading 20) Floquet multipliers 
for \reff{cAC} on the domain 
$\Om=(-\pi,\pi)$ with Neumann BC, 31 grid--points in $x$.  
Parameters $(\nu,\mu,c_3,c_5)=(1,0.1,-1,1)$, 
hence bifurcations at (restricting to the first three branches) 
$r=0$ ($k=0$, spatially homogeneous branch, black), 
$r=1/4$ ($k=1/2$, blue) and $r=1$ ($k=1$, red), see \reff{cACs2}. 
 The black dots in (a), (b) are from the analytical solution 
\reff{cACs} with $k=0$. The thick part of the black line in (a),(b) 
indicates the only stable periodic solutions. 
  \label{f1}}}
\end{figure}

There are no real eigenvalues of 
$\pa_uG$ on the trivial branch $u=0$ in this example. 
Thus, for the HBP detection 
we can simply use algorithm \heda\ from page \pageref{hed} and 
postpone to \S\ref{brusec} and \S\ref{ocsec} 
the discussion of problems which require \hedb. 
In 1D we use Neumann BC, and $n_x=31$ spatial, and (without mesh-refinement) 
$m=21$ temporal discretization points. 
Just for illustration, we compute the first two bifurcating branches, 
{\tt b1} and {\tt b2}, 
using the arclength setting from the start, while for the third 
branch {\tt b3} we first do 5 steps in natural parametrization, where {\tt TOM} 
refines the starting $t$--mesh of 21 points to 41 points. 
This produces the plots in Fig.~\ref{f1}, where the norm in (a) is 
\huga{\label{defnorm}
\|u\|_*:=\|u\|_{L^2(\Om\times (0,T), \R^N)}/\sqrt{T|\Om|}, 
}
which is our default norm for Hopf branches. 
The simulations run in less than 10 seconds 
per branch, but the rather coarse meshes 
lead to some inaccuracies. For instance, the first three HBPs, 
which analytically are at $r=0, 1/4, 1$, are obtained at 
$r=6*10^{-5}, 0.2503, 1.0033$, and (b) also shows some 
visible errors in the period $T$. 
However, these numerical errors quickly decay if we increase $n_x$ and $m$, 
and runtimes stay small. 

On {\tt b1}, initially there is 
one unstable multiplier $\ga_2$, i.e., $\ind(u_H)=1$, cf.~\reff{inddef}, 
which passes through 1 to enter the unit 
circle at the fold. Its numerical value is $10^{-5}$ close to the 
analytical result from \reff{exfl1}, and this error decreases upon 
refining the $t$--mesh. On {\tt b2} we start with $\ind(u_H)=3$, and 
have $\ind(u_H)=2$ after the fold. Near $r=0.45$ another multiplier 
moves through 1 into 
the unit circle, such that afterwards we have $\ind(u_H)=1$, with, 
for instance $\ga_2\approx 167$ at $r=1$. Thus, we may expect 
a type (i) bifurcation (cf.~p.~\pageref{biftypes}) near $r=0.45$, 
and similarly we can identify 
a number of possible bifurcation on b3 and other branches. 
The trivial multiplier $\ga_1$ is $10^{-12}$ 
close to $1$ in all these computations, using \fla. 

\begin{figure}[ht]
{\small 
\bce
\begin{tabular}{llll}
(a) BD 2D &(b) solution at b2/pt10 &(c) BD 3D&(d) solution at b2/pt15\\ 
\ig[width=0.19\textwidth]{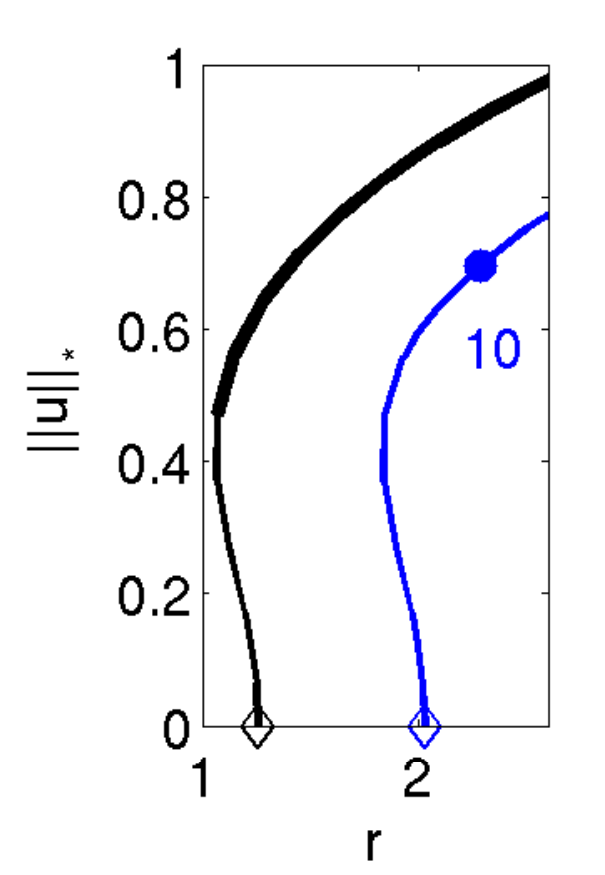}&
\ig[width=0.28\textwidth]{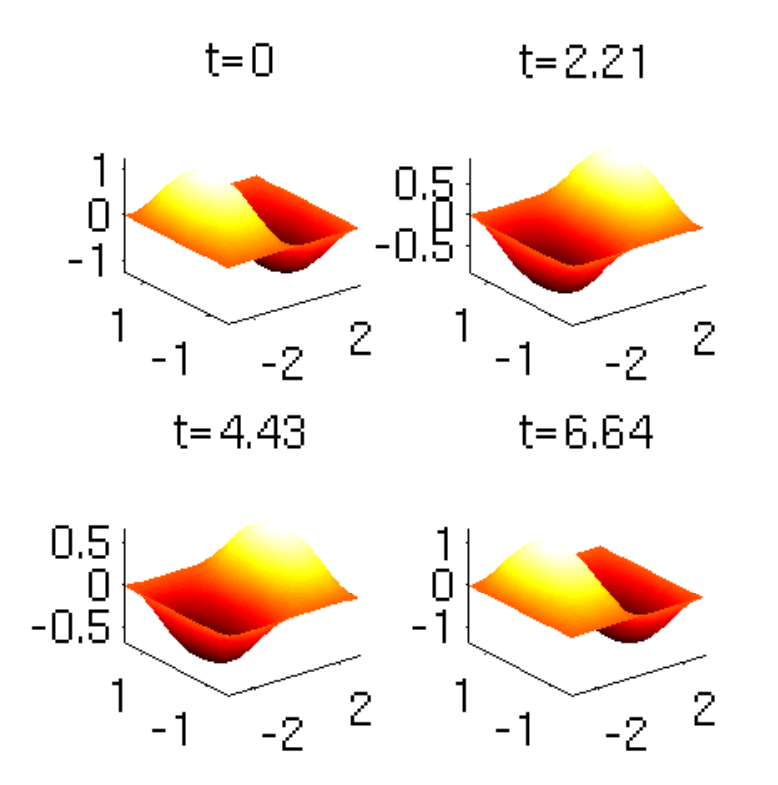}&
\hs{-2mm}\ig[width=0.2\textwidth]{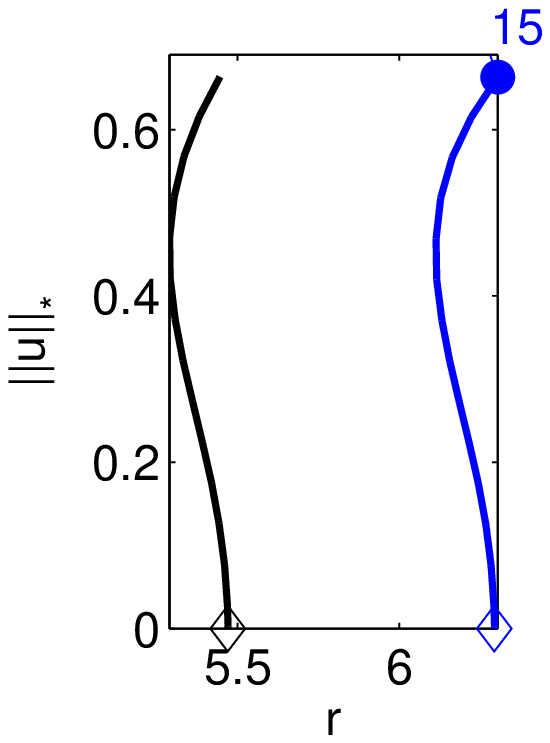}&
\hs{-0mm}\ig[width=0.22\textwidth]{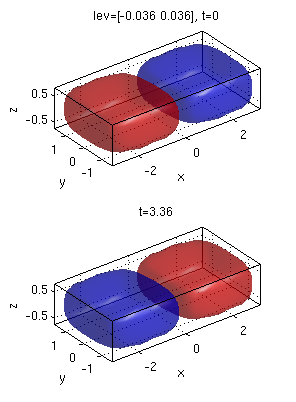}
\end{tabular}
\ece
}

\vs{-5mm}
   \caption{{\small (a) Bifurcation diagrams of the first 2 Hopf branches for 
 \reff{cAC} in 2D. (b) Solution snapshot from b2/pt10, at 
$t=0, \frac 3 {10} T, \frac 6 {10} T, \frac 9 {10} T$. (c),(d) Bifurcation 
diagram and solution snaption in 3D  }
  \label{f2}}
\end{figure}

The basic 1D setup has to be modified only slightly for 2D and 3D. 
In 2D we choose homogeneous Dirichlet BC for $u_1,u_2$. 
Then the first two HBPs are at 
$r_1=5/4$ ($k=(1/2,1)$, and $r_2=2$ ($k=(1,1)$). Figure \ref{f2} shows some 
results obtained  with a coarse 
mesh of $41\times 21$ points, hence $n_u=1722$ spatial unknowns, 
yielding the numerical values $r_1=1.2526$ 
and $r_2=2.01$. With $m=15$ temporal discretization points, the 
computation of each Hopf branch then takes less than a minute. 
Again, the numerical HBPs converge to the exact values when decreasing 
the mesh width, but at the prize of longer computations for the Hopf branches. 
For the Floquet multipliers we obtain a similar picture as in 1D. 
The first branch has $\ind(u_H)=1$ up to the fold, and $\ind(u_H)=0$ 
afterwards, and on b2 $\ind(u_H)$ decreases from 3 to 2 at the fold 
and to 1 near $r=7.2$. 

In 3D, we consider \reff{cAC} 
over $\Om=(-\pi,\pi)\times (-\pi/2,\pi/2)\times(-\pi/4,\pi/4)$. Here we 
use a {\em very} coarse tetrahedral mesh of $n_p=2912$ points, 
thus $5824$ DoF in space. Analytically, the first 2 HBPs are 
$r_1=21/4$ ($k=(1/2,1,2)$) and $r_2=6$ ($k=(1,1,2)$, but with the 
coarse mesh we numerically obtain $r_0=5.47$ and $r_1=6.29$. 
Again, this can be greatly improved by, e.g., halving the spatial 
mesh width, but then the Hopf branches become very expensive.  
Using $m=15$ and {\tt ilupack}, the computation of the branches (with 15
continuation steps each) in Fig.~\ref{f2}(b) takes about 400 seconds%
\footnote{using \reff{mbel} with $LU$ prefactorization of $A$ 
leads to about 900s runtime, and, importantly, much higher 
memory requirements of about 11GB instead of 2GB with {\tt ilupack};}%
, and using \fla\ to 
a posteriori compute the Floquet multipliers about 30 seconds per orbit.  Again, 
on {\tt b1}, $\ind(u_H)=1$ up to fold and $\ind(u_H)=0$ afterwards, 
while on {\tt b2} $\ind(u_H)$ decreases from 3 to 2 at the fold 
and to 1 at the end of the branch, 
and time integration from an IC from {\tt b2} 
yields convergence to a periodic solution from {\tt b1}. 

Additional to the code for the plots in Fig.~\ref{f1} (and Fig.~\ref{f2}), the tutorial \cite{hotutb} explains the basic steps for 
\bci 
\item switching to continuation in another parameter 
\item  using \pdep's time integration routines 
to assess the stability of periodic solutions, and in particular obtain 
the time evolution of perturbations of unstable orbits, 
\eci 
and some additional features such as ad hoc mesh refinement in $t$ for 
the arclength parametrization, and creating movies of Hopf orbits.

\subsection{Rotating patterns on a disk}
\label{rotsec}
While the Hopf bifurcations presented in \S\ref{cACsec} have been to 
(standing) oscillatory 
patterns, %
another interesting class is the 
Hopf bifurcation to 
rotating patterns, in particular to spiral waves. Such spirals are 
ubiquitous in 2D reaction diffusion problems, see, e.g., \cite{pismen06,cg09}. 
Over bounded domains, spiral waves are 
usually found numerically via time integration, see in particular 
EZSPIRAL \cite{Bark91},  with an 
$\CO(1)$ amplitude, i.e., far from bifurcation. On the other hand, 
the bifurcation of spiral waves from a homogeneous solution is usually 
analyzed over all of $\R^2$, e.g., 
\cite{hagan82, KH81, scheel98}, where the spirals are relative 
equilibria, i.e., steady states in a comoving frame. 
Moreover, spiral waves 
often undergo secondary bifurcations such as drift, meandering and 
period doubling, see \cite{Bark95, SSW99, SS07} and the references therein. 
An exception to the rule of finding spirals by time integration 
is \cite{BE07}, where they are found by growing 
them from a thin annulus towards the core using AUTO, i.e., by continuation 
in the inner radius of the annulus. Continuation in other parameters 
is then done as well, but always at $\CO(1)$ amplitude.

Here we study, on the unit disk,  the bifurcation of spiral waves 
from the zero solution 
in a slight modification of a  
real two component reaction diffusion system from \cite{GKS00}, 
somewhat similar 
to the cGL, but with Robin BC. The system reads 
\hual{\begin{split}
\pa_t u&=d_1\Delta u+(0.5+r) u+v-(u^2+v^2)(u-\al v), \\
\pa_t v&=d_2\Delta v+rv-u-(u^2+v^2)(v+\al u), 
\end{split}\label{spir1} \\
&\pa_{{\bf n}} u+10 u=0,\quad \pa_{{\bf n}} v+0.01 v=0, \label{spirbc}
}
where ${\bf n}$ is the outer normal. First (\S\ref{al0sec}) 
we follow \cite{GKS00} and set 
 $\al=0$, $d_1=0.01$, $d_2=0.015$, and take $r$ as 
the main bifurcation parameter. Then (\S\ref{al1sec}) we set $\al=1$, let 
\huga{(d_1,d_2)=\del(0.01,0.015), 
}
and also vary $\del$ which corresponds to changing the domain size 
by $1/\sqrt\del$. 

Due to the BC \reff{spirbc}, the eigenfunctions of the linearization 
around $(u,v)=(0,0)$ are build from Fourier Bessel functions 
\huga{
\phi(\rho,\vt,t)=\Re(\er^{\ri(\om t+m\vt)}J_m(q\rho)),
} 
where $(\rho,\vt)$ are polar-coordinates, and with in general complex 
$q\in\C\setm\R$. Then the modes are growing in $\rho$, 
which is a key idea of \cite{GKS00} to find modes bifurcating 
from $(u,v)=(0,0)$ which resemble spiral waves near their core.

\subsubsection{Bifurcations to rotational modes}\label{al0sec}
\begin{figure}[ht]
{\small 
\bce
\begin{tabular}{ll}
(a) Bifurcation diagram &(b) Spatial mode structure at bifurcation, h1,\ldots,h6 \\
\ig[width=40mm, height=46mm]{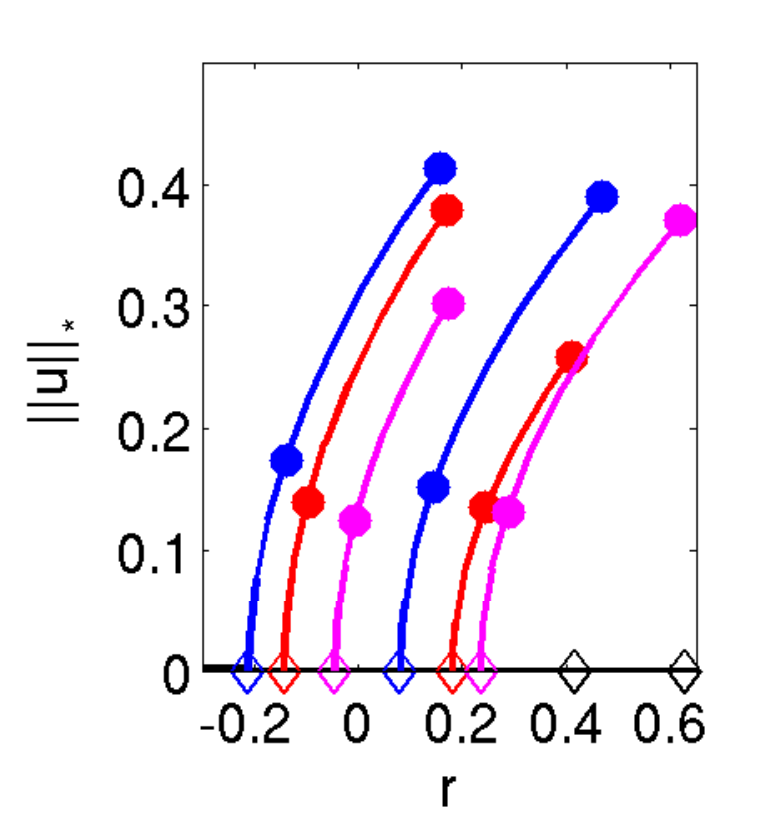}&
\hs{-4mm}\raisebox{25mm}{\begin{tabular}{l}
\ig[width=20mm, height=20mm]{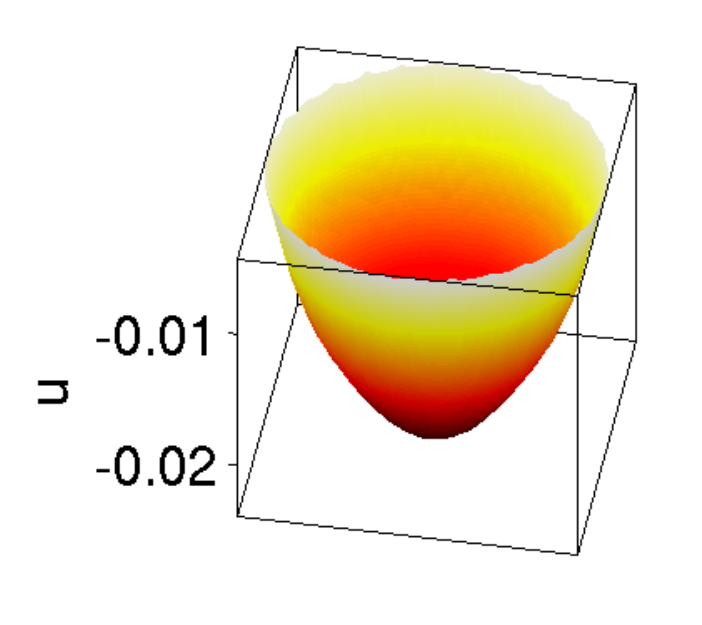}\ig[width=20mm, height=20mm]{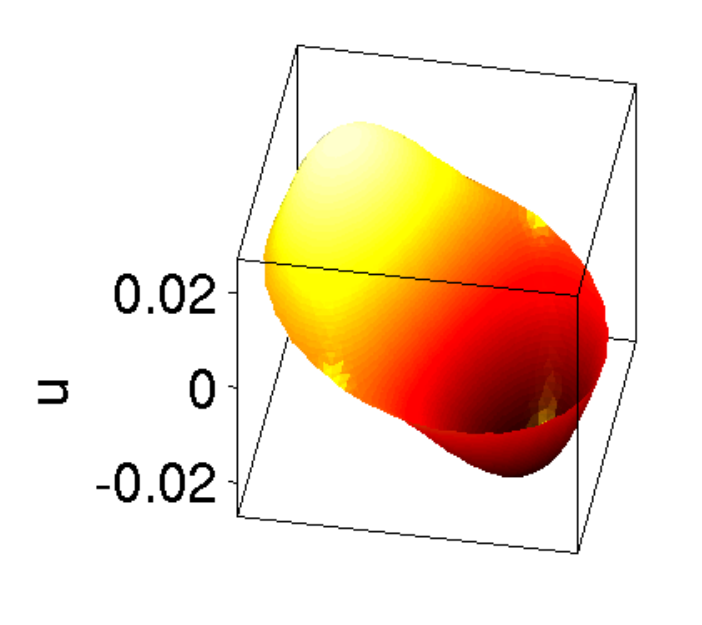}\ig[width=20mm, height=20mm]{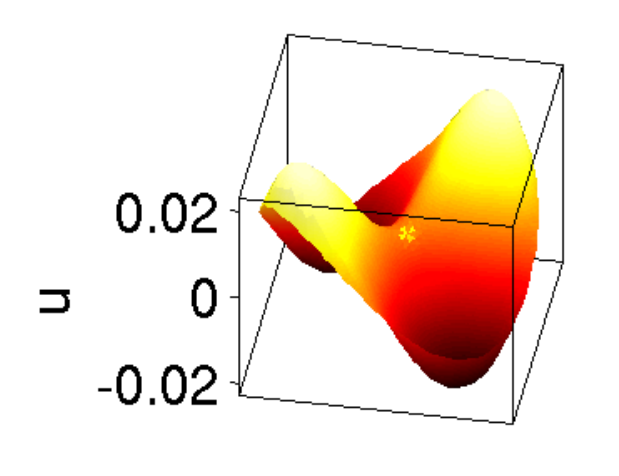}\ig[width=20mm, height=20mm]{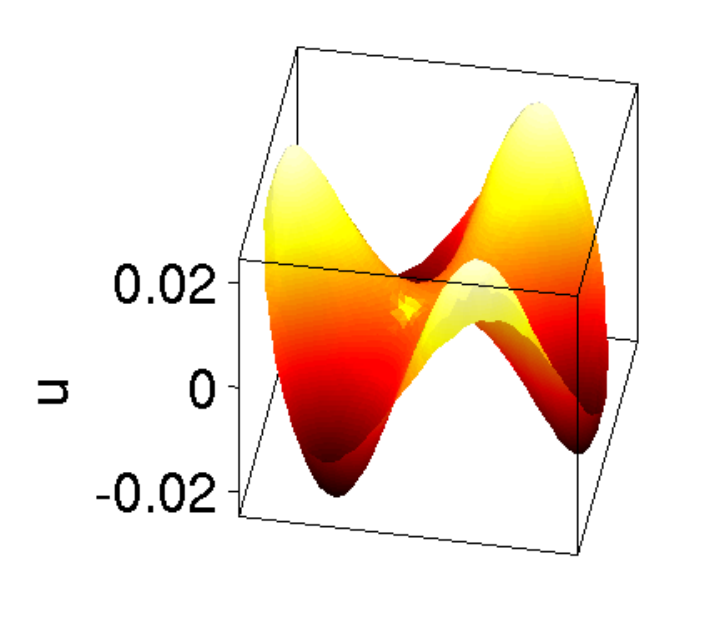}\ig[width=20mm, height=20mm]{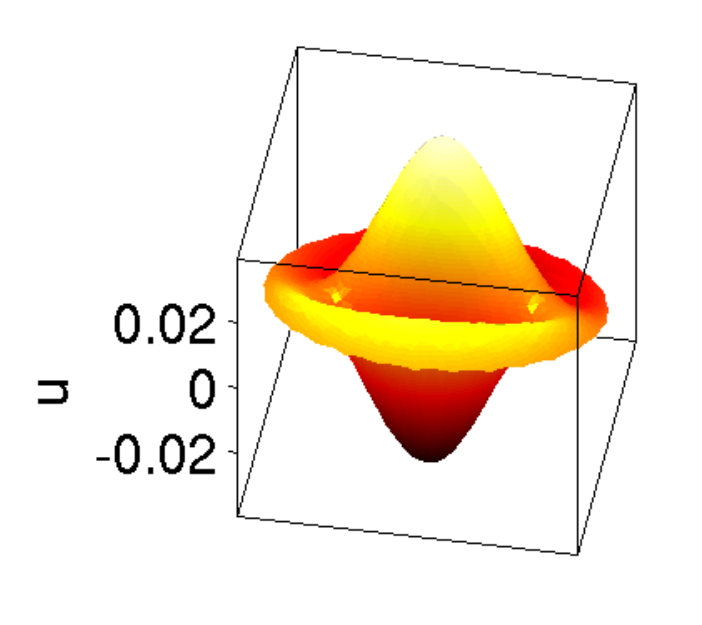}
\ig[width=20mm, height=20mm]{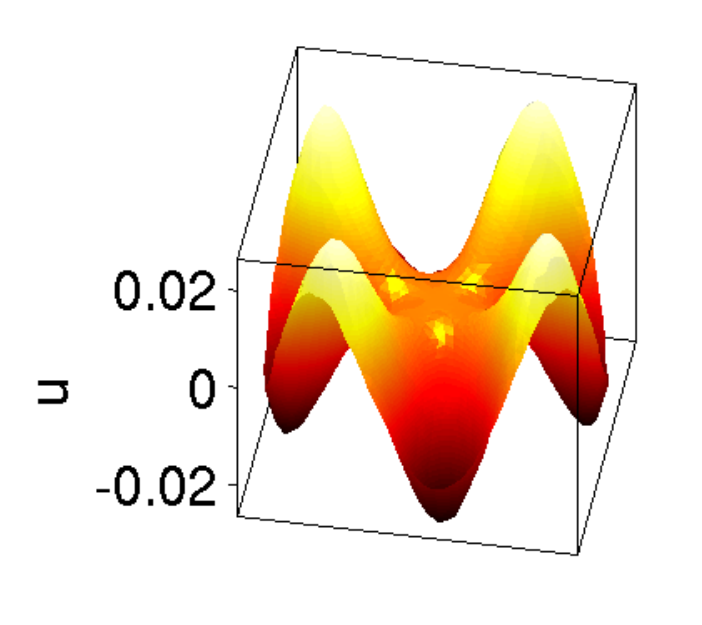}\\
\ig[width=20mm, height=20mm]{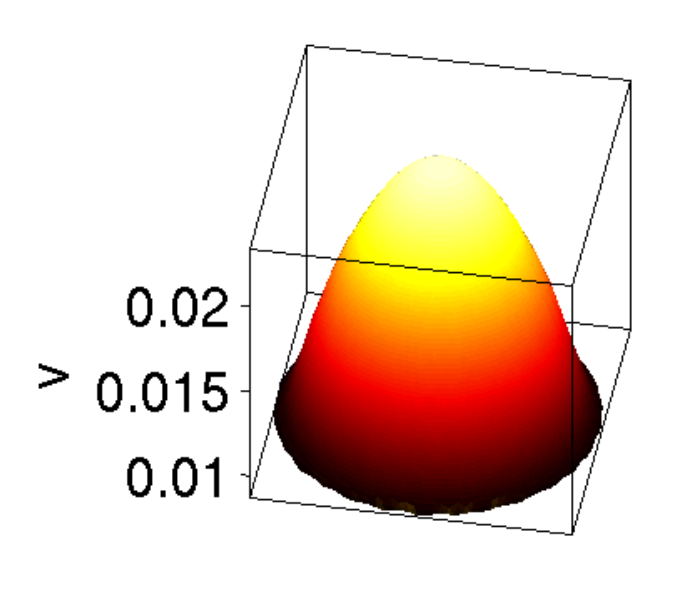}\ig[width=20mm, height=20mm]{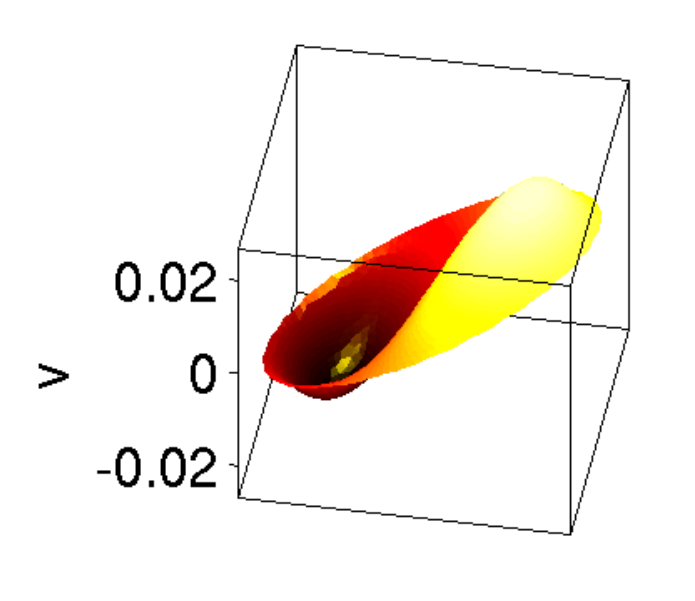}\ig[width=20mm, height=20mm]{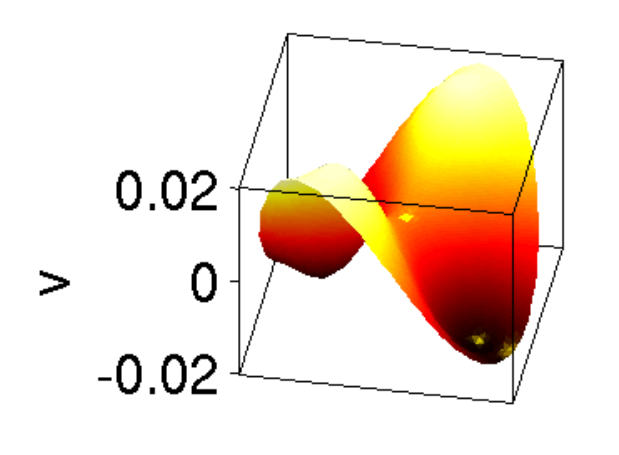}\ig[width=20mm, height=20mm]{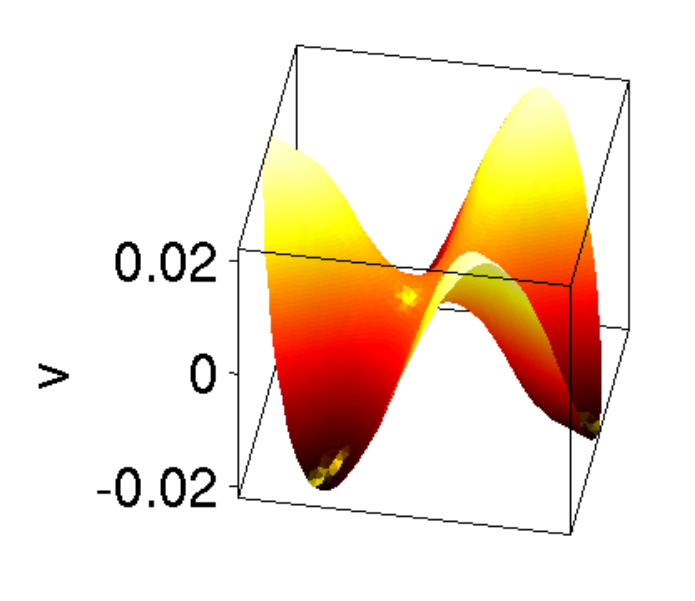}\ig[width=20mm, height=20mm]{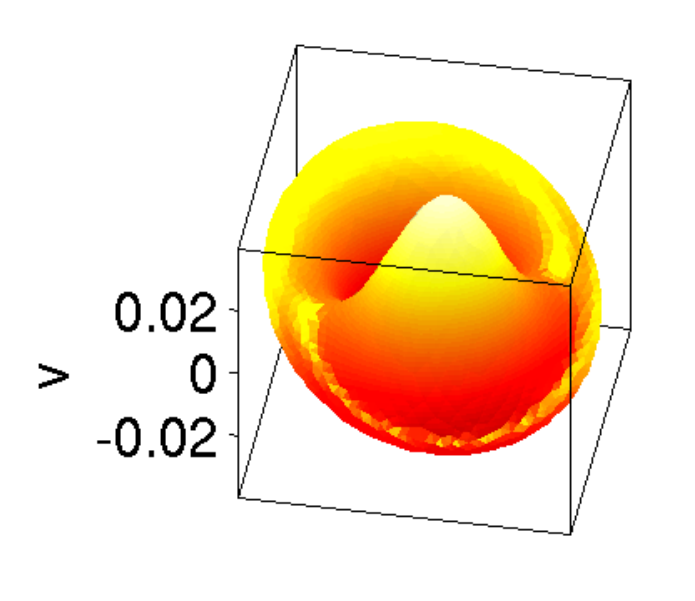}
\ig[width=20mm, height=20mm]{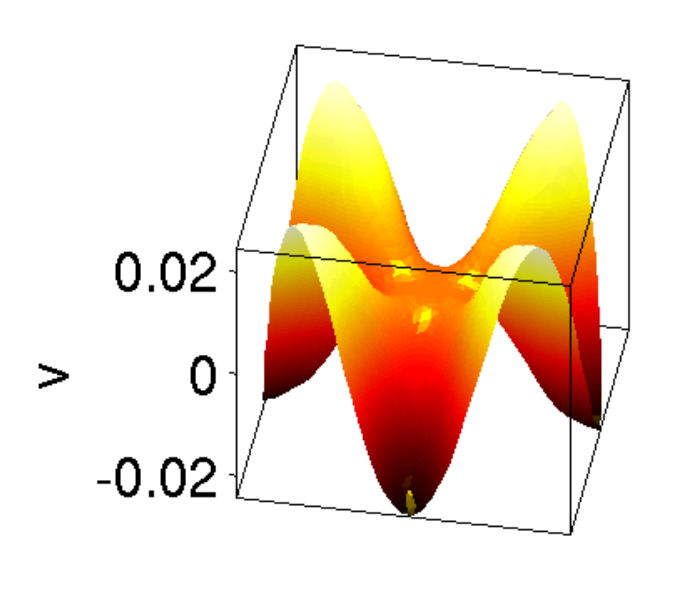}
\end{tabular}}
\end{tabular}
\vs{-0mm}

\begin{tabular}{p{70mm}p{80mm}}
(c) Zero-contours of h2, h3, h4, h5 from (b)&(d) selected snapshots 
from periodic orbits\\ %
\raisebox{25mm}{\begin{tabular}{l}
\ig[width=26mm, height=16mm]{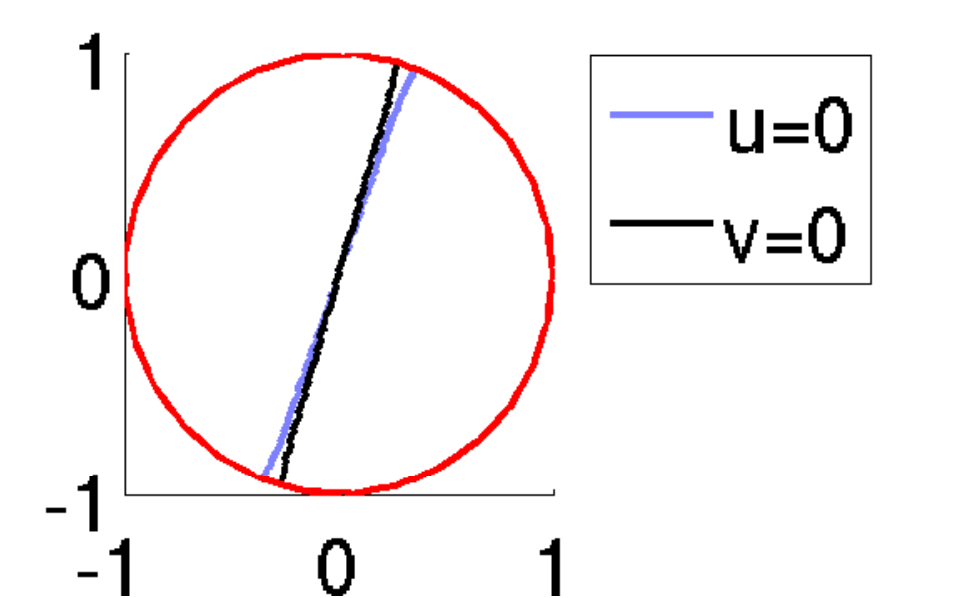}
\ig[width=26mm, height=16mm]{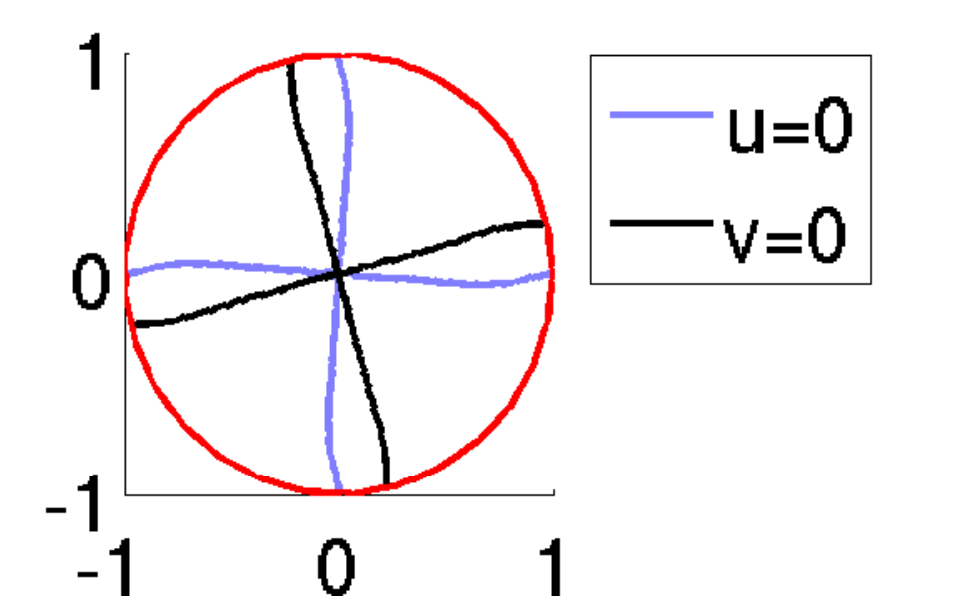}\\
\ig[width=26mm, height=16mm]{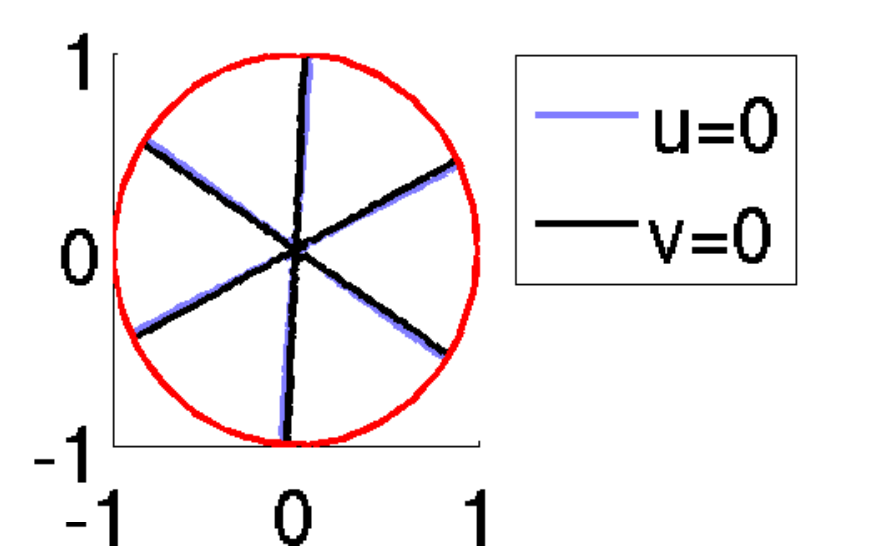}
\ig[width=26mm, height=16mm]{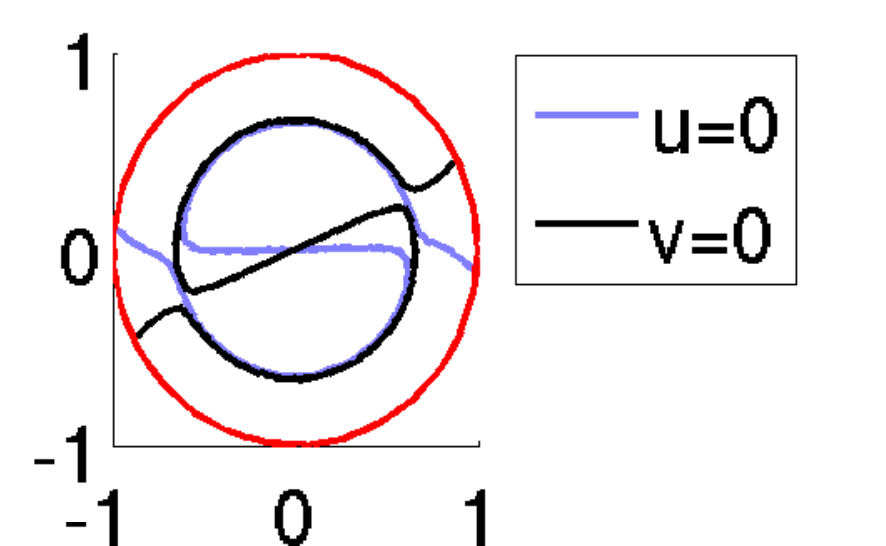}
\end{tabular}}&
\raisebox{25mm}{
\hs{-15mm}\begin{tabular}{ll}
\ig[width=50mm, height=14mm]{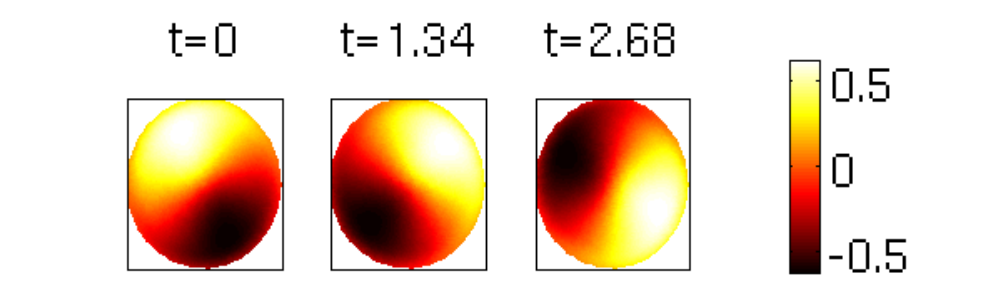}&
\hs{-5mm}\ig[width=50mm, height=14mm]{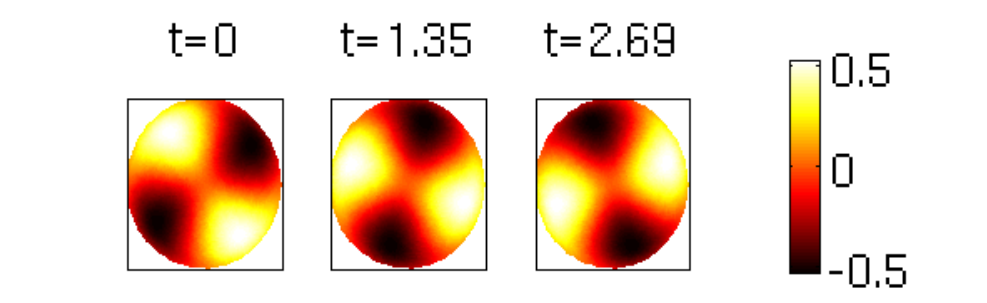}\\
\ig[width=50mm, height=14mm]{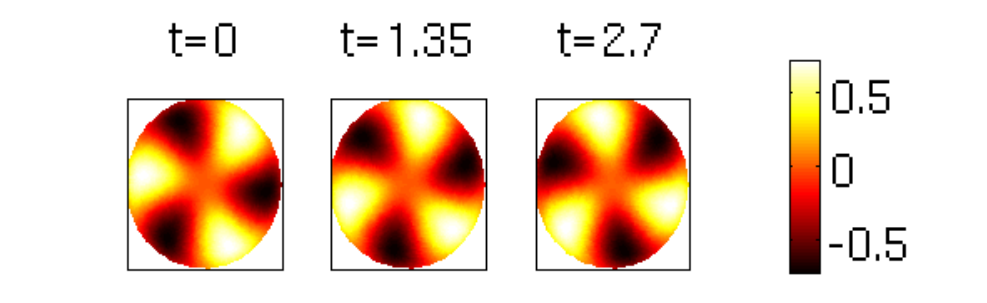}&
\hs{-5mm}\ig[width=50mm, height=14mm]{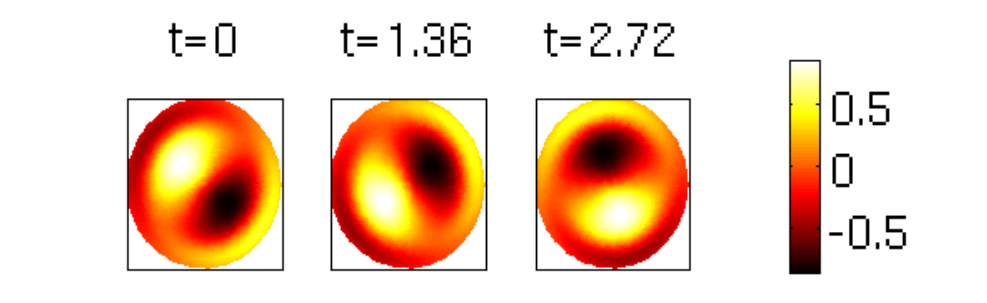}
\end{tabular}}
\end{tabular}
\ece
}

\vs{-10mm}
   \caption{{\small Basic bifurcation diagram (a) for \reff{spir1}, 
\reff{spirbc} with h1-h6 from left to right, 10 continuation 
steps for each. On each branch we mark the points 5 and 10. (b,c): 
information on initial mode structure on 
the first six bifurcating branches. 
 (d)  Example plots last points in h2, h3 (upper row) and 
h4, h5 (lower row). 
Snapshots of $u$ at $t=0, T_j/5$ and $2T_j/5$, with $T_j$ the actual period. 
  \label{spf1}}}
\end{figure}

The trivial homogeneous branch $(u,v)=(0,0)$ is stable up to $r\approx -0.21$, 
and Fig.~\ref{spf1}(a) shows the first 6 bifurcating branches h1,h2,\ldots, 
h6, from left to right, while (b) shows the spatial modes for 
h1-h6 at bifurcation, with  mode numbers $m=0,1,2,3,2,4$. 
We discretized \reff{spir1}, 
\reff{spirbc} with a mesh of 1272 points, hence $n_u=2544$ DoF, and a 
coarse temporal discretization of 11 points, which 
yields about 2 minutes for the computation of each branch, 
with 10 points on each. Example plots 
of solutions on the last points on the branches are given in (e), 
with $T$ near $2\pi$ for all branches. 

The nontrivial solutions from Fig.~\ref{spf1}(a),(d) 
are ``rotations'', except for the spatial $m=0$ mode h1. 
 To discuss this, we return to (c), which shows the nodal lines for the 
components $u$, $v$ at bifurcation of h2 to h5 
(vector $\Psi$ in \reff{hoswitr}).  
The pertinent observation is that h2 to h6 (not shown) do not have nodal lines, 
i.e., $u(x)v(x)\ne 0$ except at $x=0$.%
\footnote{The zero lines for 
h3 are close together, but not equal; for h1 we have $u(x,0)<0$ and 
$v(x,0)>0$ for all $x$.} Thus, the 
branches h2 to h6 cannot consist of oscillatory patterns but must rotate. 
On the other hand, this rotation must involve higher order modes, 
and thus becomes more visible, i.e., {\em almost} 
(but never perfectly) rigid, at larger amplitude.

To assess the numerical accuracy, in Table \ref{sptab1} we compare 
the numerical values for the Hopf points and the temporal wave number $\om$ 
with the values from \cite{GKS00}, who compute $r_4,r_5,r_6$, 
(and three more Hopf points) using semi analytical methods, and some 
numerics based on the \mlab\ \ptool\ with fine meshes. Given our coarse 
mesh we find our results reasonably close, and again our values converge 
to the values from \cite{GKS00} under mesh refinement. 

\begin{table}[ht]
\bce
\begin{tabular}{r|llllll} 
branch&h1&h2&h3&h4&h5&h6\\
$r$&-0.210&-0.141&-0.044&0.079&0.182&0.236\\
$\om$&0.957&0.967&0.965&0.961&0.953&0.957\\
$r^*$&NA&NA&NA&0.080&0.179&0.234\\
$\om^*$&NA&NA&NA&0.961&0.953&0.957\\
\hline
$\ind(u_H)$, pt5&0&2&6&12&16&20\\
$\ind(u_H)$, pt10&0&2&4&8&18&16
\end{tabular}
\ece
\vs{-5mm}
\caption{Comparison of HBPs with \cite{GKS00} (starred values), and Floquet 
indices at points on branches.\label{sptab1}}
\end{table}

The last two rows of Table \ref{sptab1} give the Floquet indices 
of points on the branches, where $\emu$ (cf.~\reff{emudef}) 
is around $10^{-10}$ for each computation. 
All branches except h1 are unstable, and the instability indices 
increase from left to right, and  also vary along the 
unstable branches. 
However, altogether \reff{spir1},\reff{spirbc} with 
$(\al,\del)=(0,1)$ does not appear to be very interesting 
from a dynamical and pattern forming point of view, as 
time--integration yields that  for $r>r_0=-0.21$ solutions 
to generic initial conditions converge to a periodic orbit from h1. 
Thus, we next choose $\al=1$ to switch on a rotation also in the nonlinearity. 

\subsubsection{Spiral waves}\label{al1sec}
\begin{figure}[ht]
{\small 
\bce
\begin{tabular}{lll}
(a) Bifurcation diagram &(b) Profiles at selected points& (c) 
Multipliers\\ %
\ig[width=35mm, height=45mm]{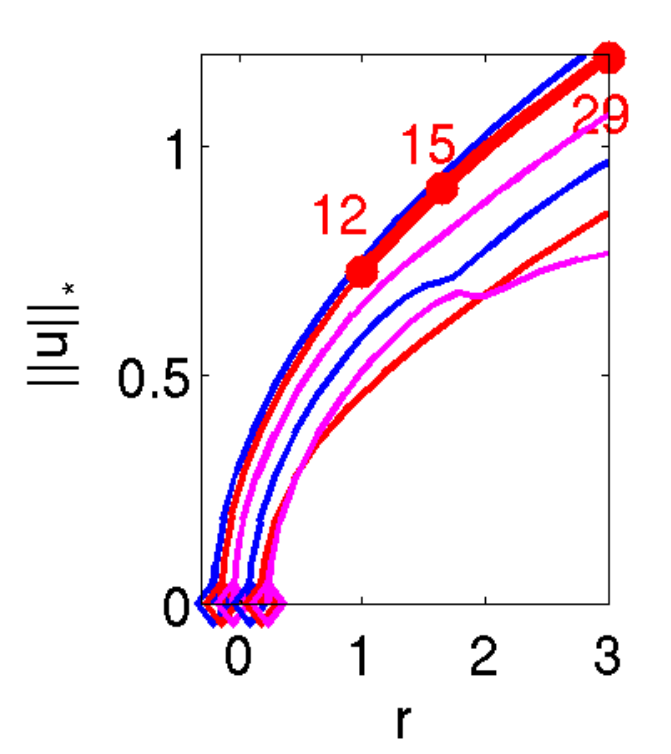}&
\hs{-4mm}\raisebox{23mm}{\begin{tabular}{l}
\ig[width=18mm, height=19mm]{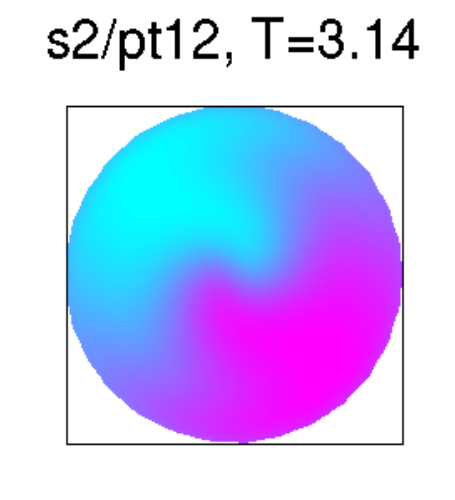}
\ig[width=18mm, height=19mm]{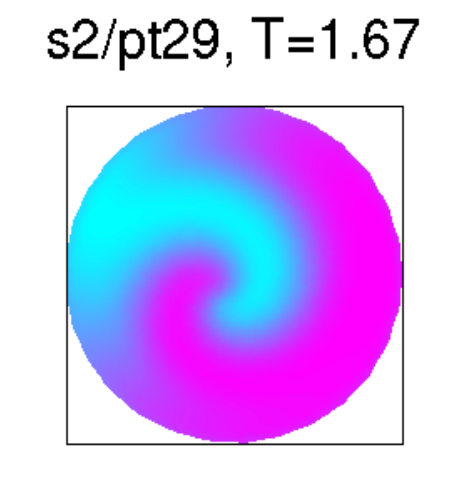}
\ig[width=18mm, height=19mm]{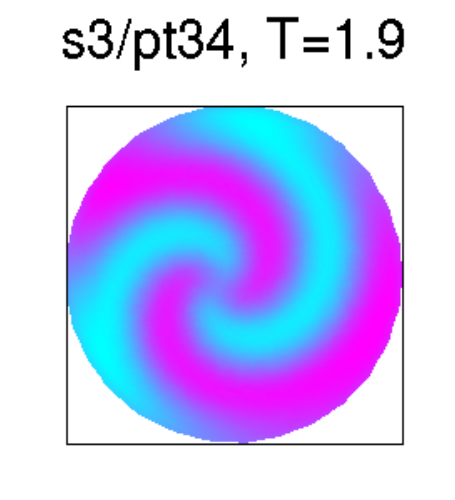}\\
\ig[width=18mm, height=19mm]{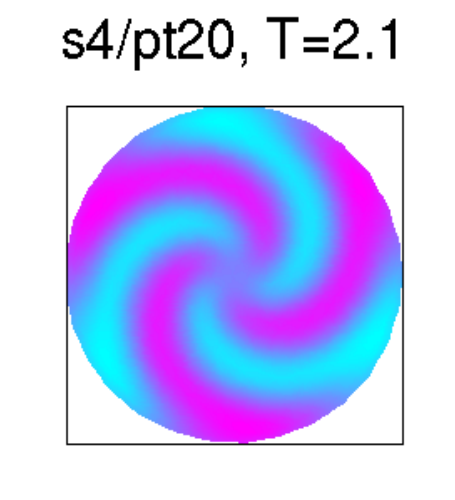}
\ig[width=18mm, height=19mm]{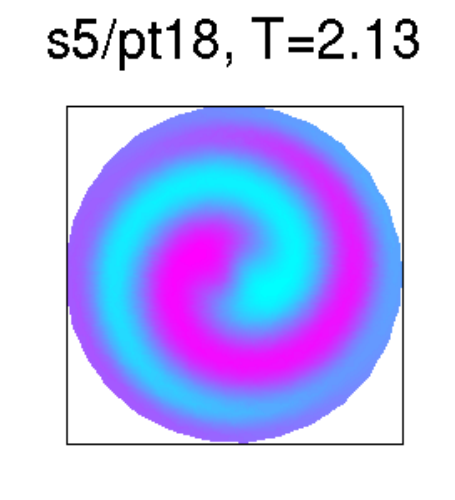}
\ig[width=18mm, height=19mm]{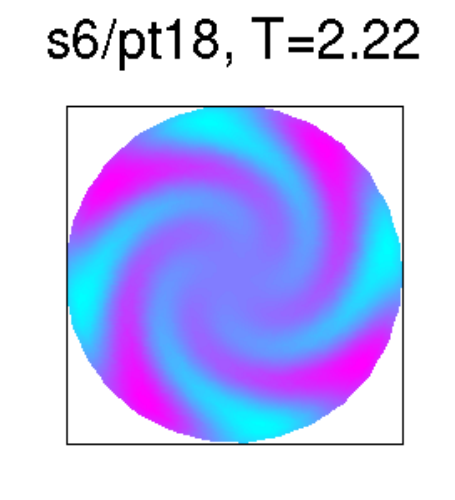}
\end{tabular}}
&\raisebox{25mm}{\begin{tabular}{l}
\ig[width=20mm]{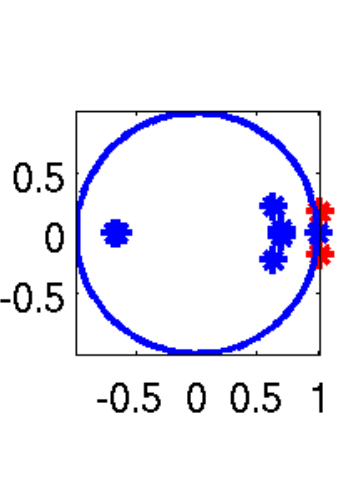}\\[-5mm]
\ig[width=20mm]{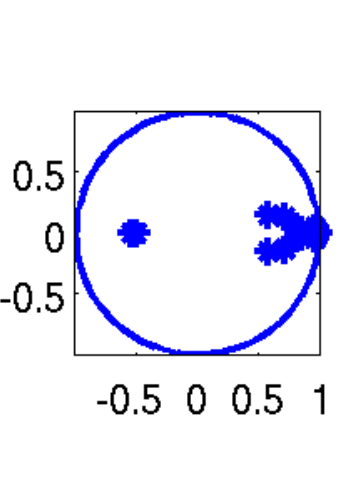}\\[-53mm]
s2/pt12\\[21mm]
s2/pt15\\[20mm]
\end{tabular}}
\end{tabular}
\vs{-0mm}
\ece
}

\vs{-10mm}
   \caption{{\small Bifurcation diagram (a) with branches s1,\ldots,s6 left to 
right, and selected profiles (b) and 
Floquet spectra (c). The (non--rotational) branch s1 is stable for all 
$r$ but plotted as a thin line (first blue line in (a)) for graphical 
reasons. The first two plots in (b) are both from s2, indicating the 
more pronounced spiral nature for larger $r$ (on all branches); remaining 
plots all at $r=3$. $T$ in (b) indicating the period, which decreases in $r$ 
and increases with number $m$ of arms of the spirals. 
  \label{spf10}}}
\end{figure}

For $(\al,\del)=(1,1)$ the linearization around $(u,v)=(0,0)$ and 
thus also the Hopf bifurcation points $r_{h1},\ldots, r_{h6}$ are as 
in \S\ref{al0sec}. However, the nonlinear rotation yields a spiral 
wave structure on the branches s2,\ldots, s6 bifurcating 
at these points, see Fig.~\ref{spf10}(b), 
where we only give snapshots of $u(\cdot,0)$, at $r=1$ and at 
$r=3$ for s2, and $r=3$ for the remaining branches.  
On s2 , s3, s4, and s6 
the solutions rotate almost rigidly in counterclockwise direction 
with the indicated period $T$, while on s5 we have a clockwise rotation. 
Thus, on s2, s3, s4 and s6 we have inwardly moving spirals, also 
called anti-spirals \cite{EpSci2001}. Moreover, again s1 is stable for 
all $r>r_{h1}$, but additionally s2 becomes stable for $r>r_1\approx 1$, 
see Fig.~\ref{spf10}(c), while s5 and the $m$--armed spirals with $m>1$ on 
s3, s4, s6 are unstable, as should be expected \cite{hagan82}; also note 
how the core becomes flatter with an increasing number of arms, again 
cf.~\cite{hagan82} and the references therein. 

\begin{figure}[ht]
{\small 
\bce
\begin{tabular}{lll}
\mbox{(a) Continuation in $\del$, $r=3$,} &  (b) Multipliers
&(c) Continuation in $r$, $\del=0.2$\\
\begin{tabular}{l}
\ig[width=40mm, height=25mm]{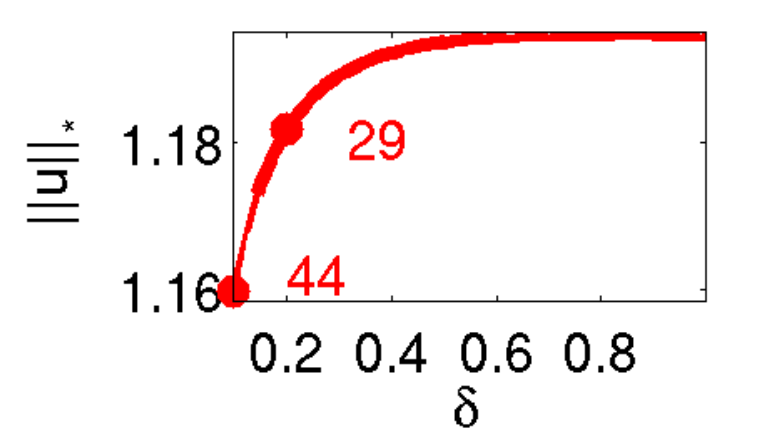}\\
\ig[width=23mm, height=25mm]{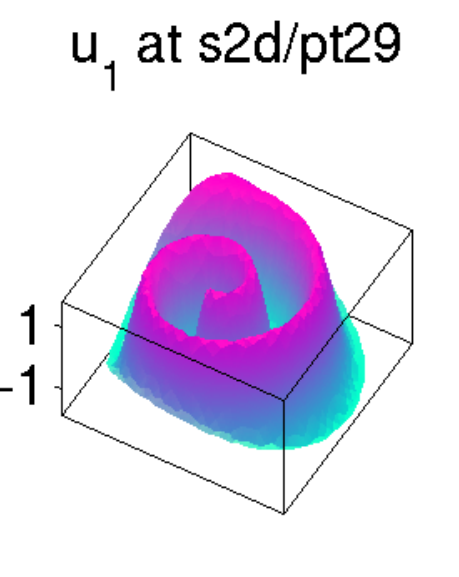}
\ig[width=23mm, height=25mm]{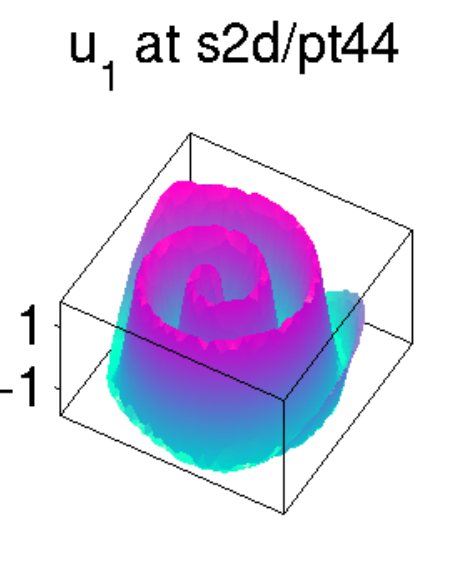}
\end{tabular}&
\hs{-0mm}\raisebox{5mm}{\begin{tabular}{l}
\ig[width=20mm, height=29mm]{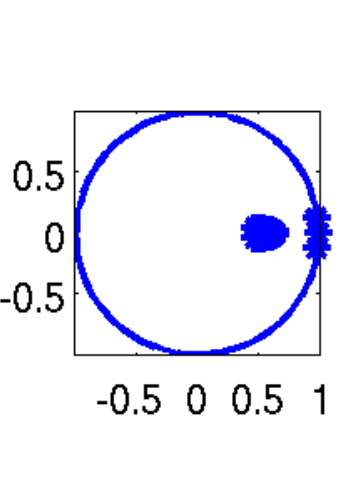}\\[-5mm]
\ig[width=20mm, height=29mm]{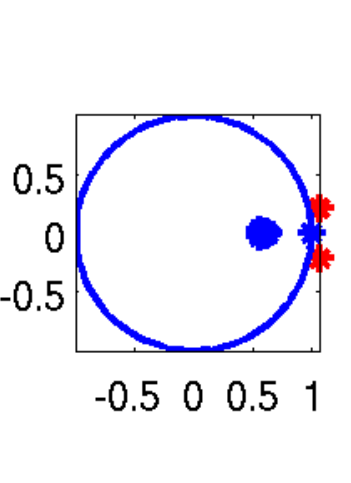}\\[-53mm]
s2d/pt29\\[21mm]
s2d/pt44\\[20mm]
\end{tabular}}&
\begin{tabular}{l}
\ig[width=40mm, height=25mm]{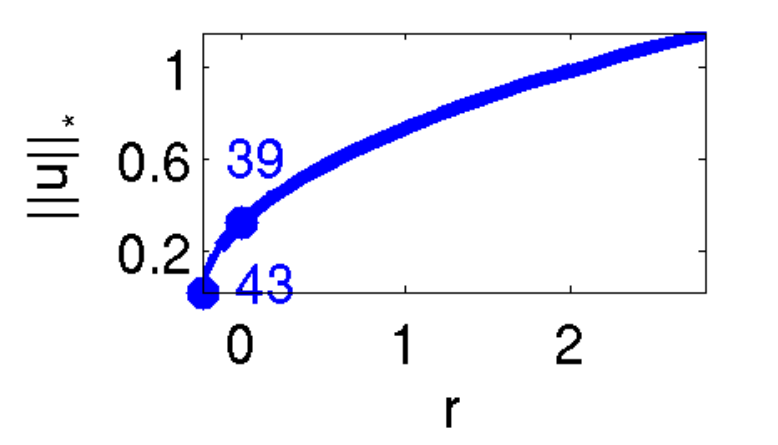}\\
\ig[width=23mm, height=25mm]{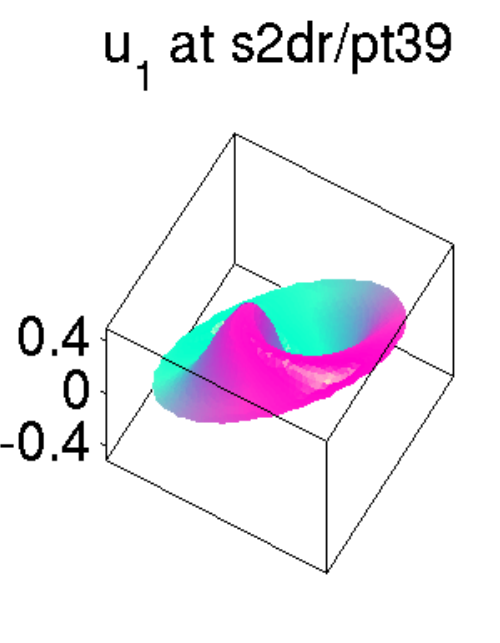}
\ig[width=23mm, height=25mm]{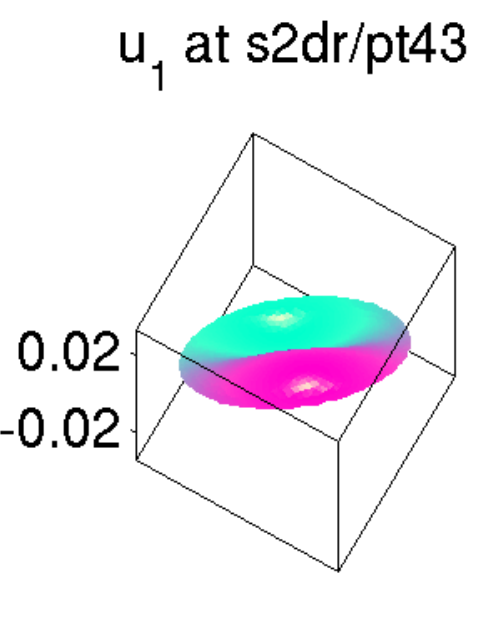}
\end{tabular}
\end{tabular}
\vs{-0mm}
\ece
}

\vs{-4mm}
   \caption{{\small (a) Continuation of the one armed spiral in $\del$ (inverse 
domain-size). Over a larger domain the spiral 
nature (of all spirals) is more visible. (b) Multipliers for points in (a). 
(b) Continuation of pt29 from (a) 
in $r$; over a larger domain the ``one-armed spiral'' is stable 
for lower amplitudes. 
  \label{spf11}}}
\end{figure}

In Fig.~\ref{spf11}(a) we  first continue $(u,v)$ from s2 at $r=3$ in $\del$ 
to $\del=0.1$, i.e., to domain radius $\sqrt{10}$ (branch s2d). 
As expected, with the growing 
domain the spirals become more pronounced (see the example plots in (c)). 
The solutions stay stable down to $\del=\del_1\approx 0.15$, as illustrated 
in (b).  In (c) 
we continue the solution from s2d/pt29 (with $\del=0.2$) 
again in $r$ down to $r=r_{h2}^*\approx 
-0.22$, which is the associated Hopf bifurcation point over 
a circle of radius $\sqrt{5}$, see also the last plot in (c), 
which is very close to bifurcation. Now the 1-armed spiral like solution 
is stable also for rather small amplitude.

The model with $(r,\al,\del)=(3,1,0.1)$ is also be quite rich dynamically. 
Besides solutions converging to s1, the 1-armed spiral 
s2 has a significant domain of attraction, %
but there are also various at least meta-stable 
solutions, which consist of long-lived oscillations (with or 
without rotations). See \cite{hotutb} for comments on how to run such 
dynamical simulations.

\subsection{An extended Brusselator%
}\label{brusec}
As an example with an interesting interplay between 
stationary patterns and Hopf bifurcations, 
with typically many eigenvalues with small real parts, 
 and where therefore detecting 
HBPs  without first setting a guess for a shift $\om_1$ 
is problematic, we consider an ``extended Brusselator'' problem from \cite{yd02}. 
This is a three component reaction diffusion system of the form 
\huga{\label{ebru}
\pa_t u=D_u\Delta u+f(u,v)-cu+dw,\quad  
\pa_t v=D_v\Delta v+g(u,v), \quad \pa_t w=D_w\Delta w+cu-dw,
}
where $f(u,v)=a-(1+b)u+u^2v$, $g(u,v)=bu-u^2v$, with kinetic parameters 
$a,b,c,d$ and diffusion constants $D_u,D_v,D_w$. 
We consider \reff{ebru} on rectangular domains in 1D and 2D, with 
homogeneous Neumann 
BC for all three components. The system has the trivial spatially homogeneous 
steady state 
$$U_s=(u,v,w):=(a,b/a,ac/d),
$$ 
and in suitable parameter regimes it shows co-dimension 2 points between 
Hopf, Turing--Hopf (aka wave), and (stationary) Turing bifurcations 
from $U_s$. We follow \cite{yd02} and fix the parameters 
\huga{\label{brupar}
(c,d,D_u,D_v,D_w)=(1,1,0.01,0.1,1). 
}

\begin{figure}[ht]
\bce{\small 
\begin{tabular}{llll}
(a) &(b)&(c)\\
\ig[width=0.4\textwidth]{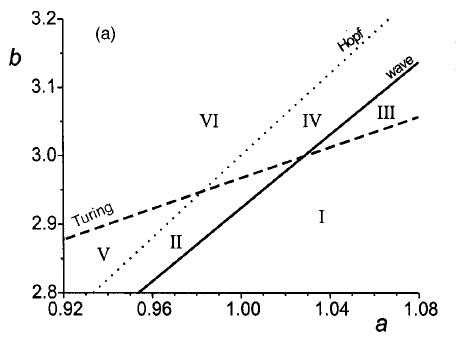}&
\raisebox{25mm}{\begin{tabular}{l}
\ig[width=0.15\textwidth, height=22mm]{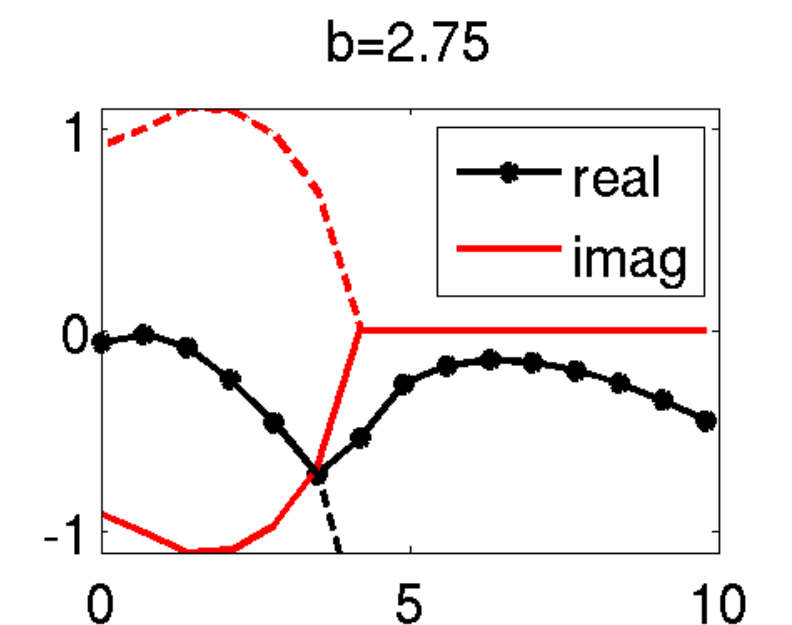}\\
\ig[width=0.15\textwidth, height=26mm]{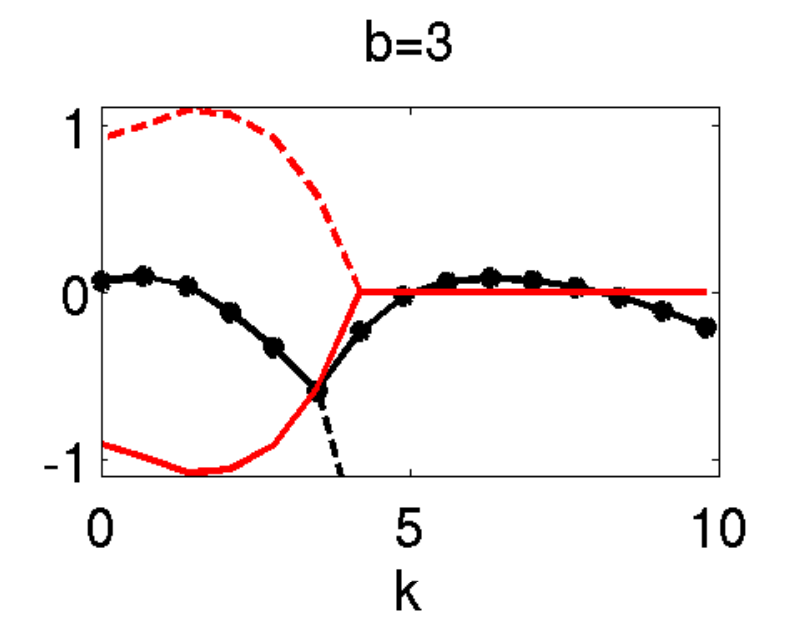}
\end{tabular}}
&\ig[width=0.3\textwidth]{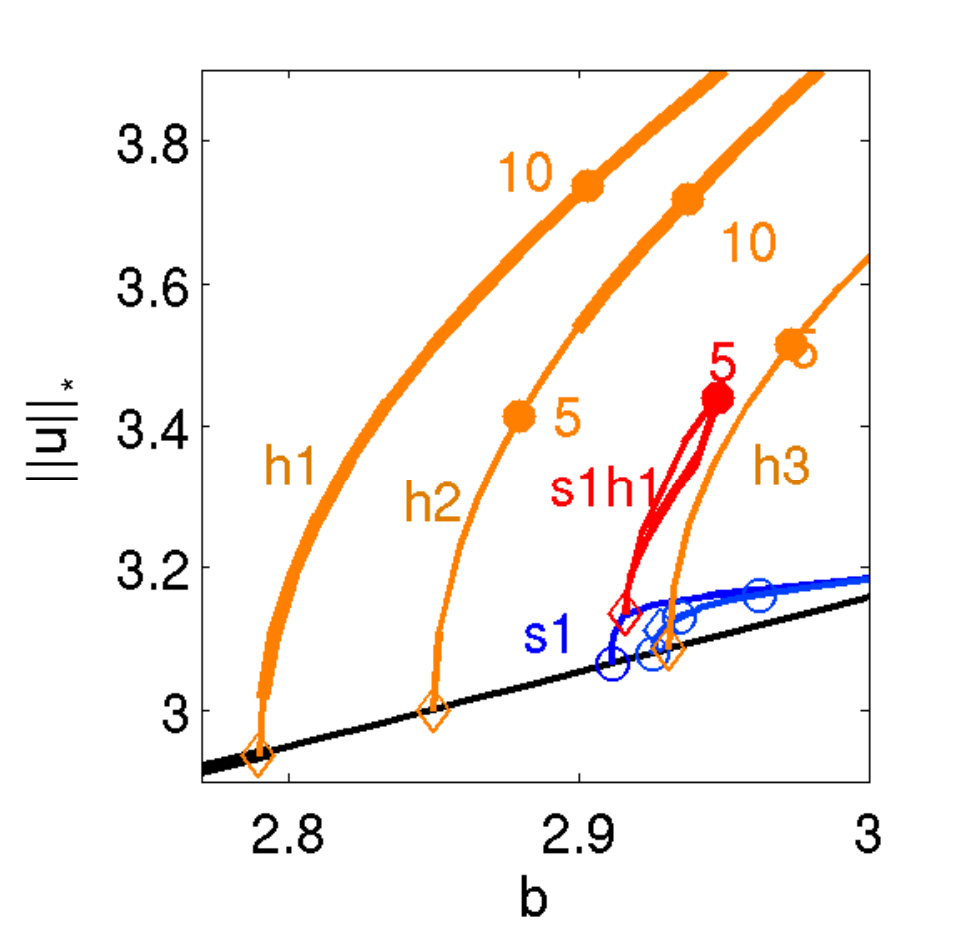}\\
\end{tabular}\\
(d) solution plots\\
\ig[width=0.22\textwidth]{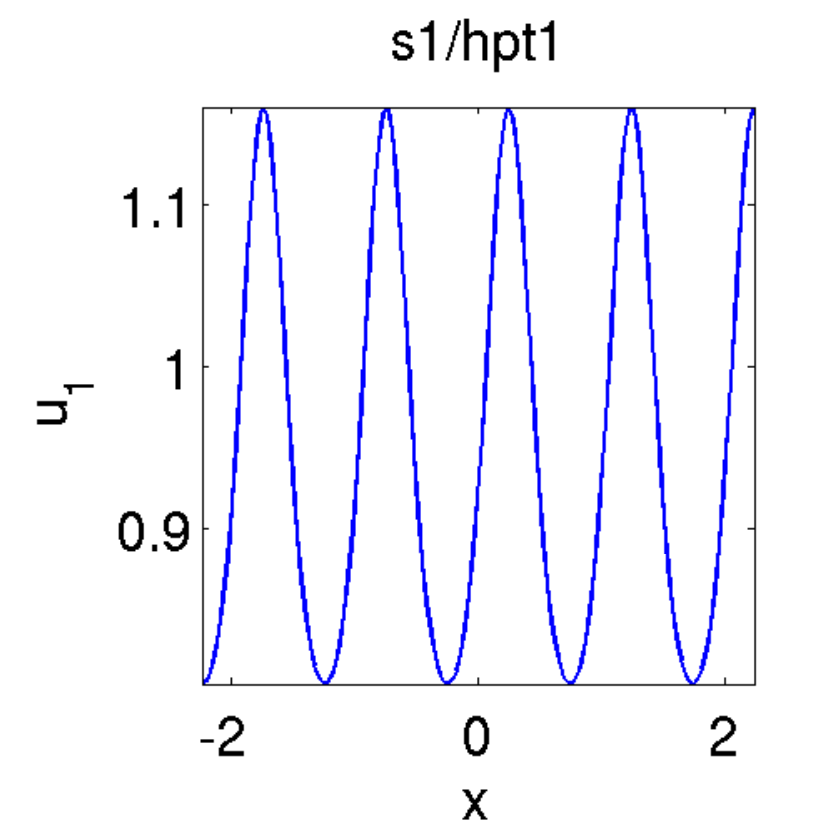}
\ig[width=0.24\textwidth]{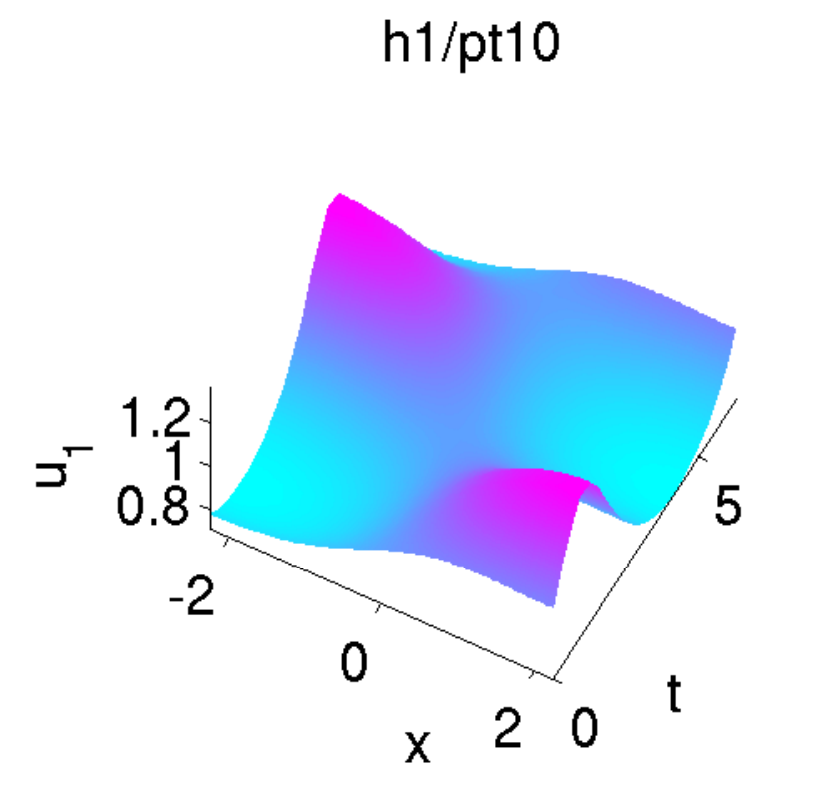}
\ig[width=0.24\textwidth]{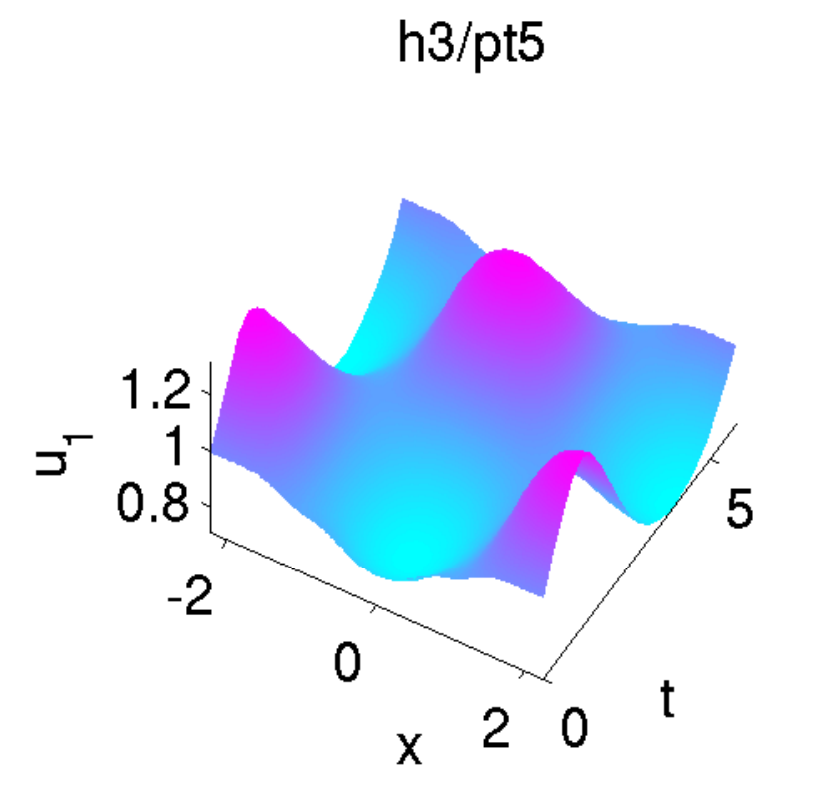}
\ig[width=0.24\textwidth]{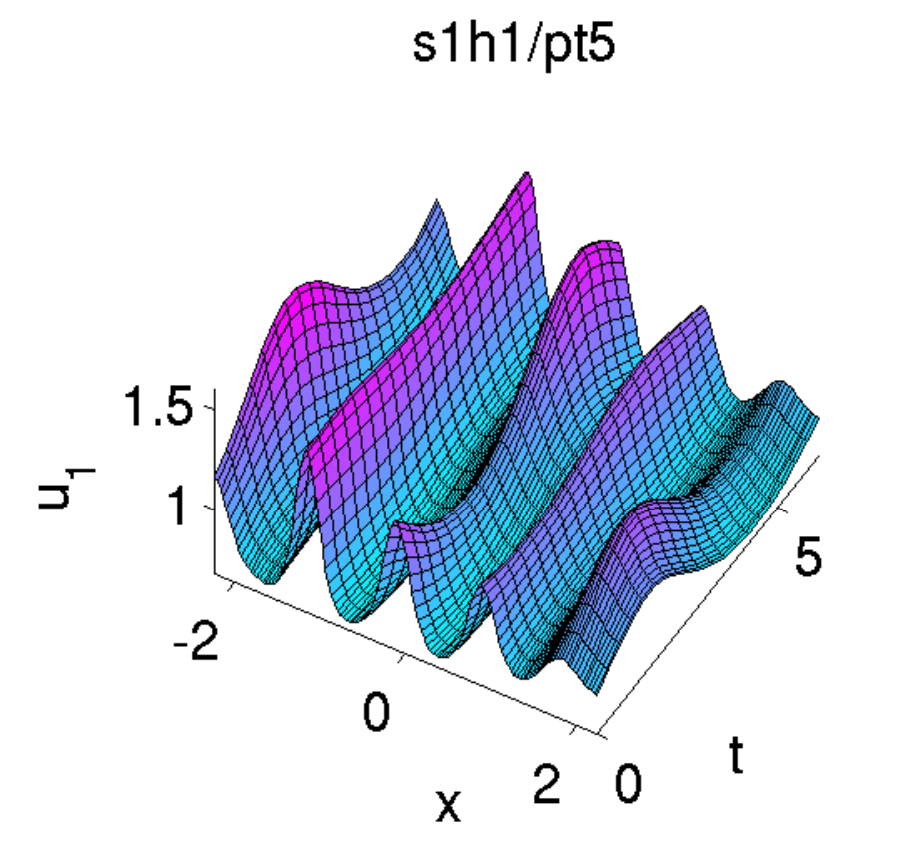}
}
\ece 

\vs{-2mm}
   \caption{{\small (a) Parameter plane with Hopf, Turing--Hopf (wave) and 
Turing instability lines for \reff{ebru}, reprinted with permission 
from \cite{yd02}, copyright 2002, AIP Publishing LLC. 
(b) Spectra for increasing $b$ 
at $a=0.95$. Contrary to the \pdep\ convention that due to 
$\pa_t u=-G(u)$ eigenvalues 
with {\em negative} real parts yield instabilities, here we directly 
plot the spectra of $-\pa_u G$, such that instability occurs for 
eigenvalues with positive real parts. 
The first instability (Turing--Hopf) occurs at $b\approx 2.794$, with 
$k_c\approx 0.7$. The admissible wave-numbers $k$ on a domain $(-l_x,l_x)$ 
with $l_x=0.5\pi/k_c$ are indicated by the dots.  
(c),(d): (partial) bifurcation diagram, and example plots on $\Om=(-l_x,l_x)$. 
  \label{bf0}}}
\end{figure}

Figure \ref{bf0}(a) then shows a characterization of the 
pertinent instabilities of $U_s$ in the $a,b$ plane. $U_s$ is stable 
in region I, and can loose stability by $(a,b)$ crossing the 
Turing line, which yields the bifurcation of stationary 
Turing patterns, or the wave (or Turing--Hopf) line, which 
yields oscillatory Turing patterns. Moreover, there is the 
``Hopf line'' which corresponds to Hopf--bifurcation with spatial 
wave number $k=0$.

In the following we fix $a=0.95$ and take $b$ as the primary 
bifurcation parameter. Figure \ref{bf0}(b) illustrates the different 
instabilities from (a). As we increase $b$ from 2.75, 
we first cross the Turing--Hopf line, with first instability 
at critical spatial wave number 
$k_{{\rm TH}}\approx 0.7$, then the Hopf line, and 
finally the Turing line with critical wave number $k_{\rm T}\approx 6.4$. 
To investigate the bifurcating solutions (and some secondary bifurcations) 
with \pdep, we need to choose a domain $\Om=(-l_x/2, l_x/2)$ (1D), 
where due to the Neumann BC $l_x$ should be chosen as a (half integer) 
multiple of $\pi/k_{{\rm TH}}$. For 
simplicity we take the minimal choice $l_x=0.5\pi/k_{{\rm TH}}$, which restricts the 
allowed wave numbers to multiples of $k_{{\rm TH}}$, as indicated by 
the black dots in Figure \ref{bf0}(b). Looking at the sequence 
of spectral plots for increasing $b$, we may then expect 
first the Turing--Hopf branch h1 with $k=k_{{\rm TH}}$, then a Hopf branch h2 with $k=0$, 
then two Turing branches s1, s2 with $k=6.3$ and $k=7$, then 
a Turing--Hopf branch h3 with $k=2k_{{\rm TH}}$, and so on, and this is what 
we obtain from the numerics, as illustrated in (c) and (d), 
using a coarse mesh with $101$ grid points, hence $3\times 101=303$ DoF in space. 
 
Besides stationary secondary bifurcations we also 
get a rather large number of Hopf points on the Turing 
branches, and just as an example we plot the (Turing)Hopf branch s1h1 
bifurcating from the first Hopf point on s1. The example plots 
in (d) illustrate that solutions on s1h1 look like a superposition of 
solutions on s1 and h1. 
Such solutions were already obtained in \cite{yd02} 
from time integration, such that at least some these solutions also have some 
stability properties, see also \cite{ye03} for similar phenomena. 
By following the model's various bifurcations, this can be studied in 
a more systematic way.

\begin{figure}[ht]
\bce{\small 
\begin{tabular}{llll}
(a) h1/pt10, $\ind=0$&(b) h2/pt5, $\ind=2$&(c) h2/pt10, $\ind=0$&
(d) h3/pt5, $\ind=5$\\[-0mm]
\ig[width=0.2\textwidth]{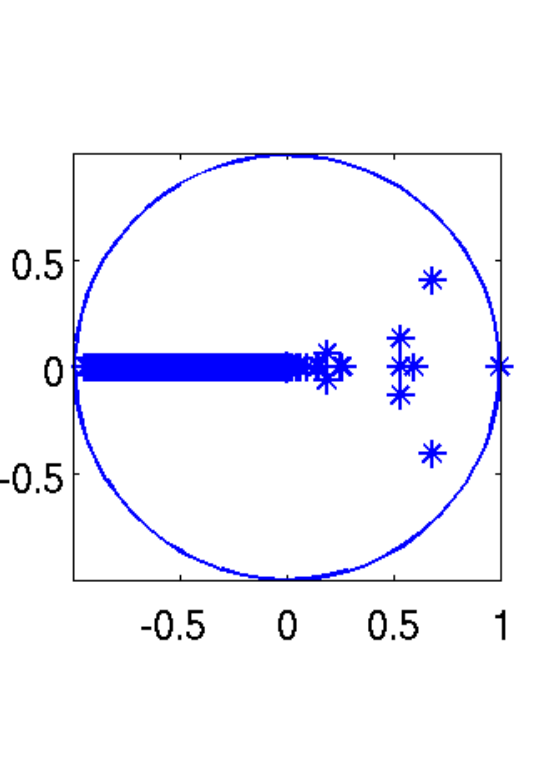}&
\ig[width=0.2\textwidth]{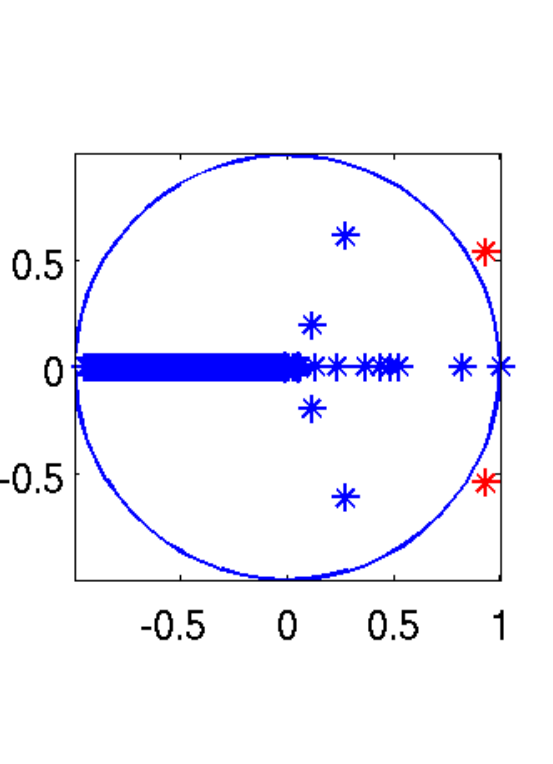}&
\ig[width=0.2\textwidth]{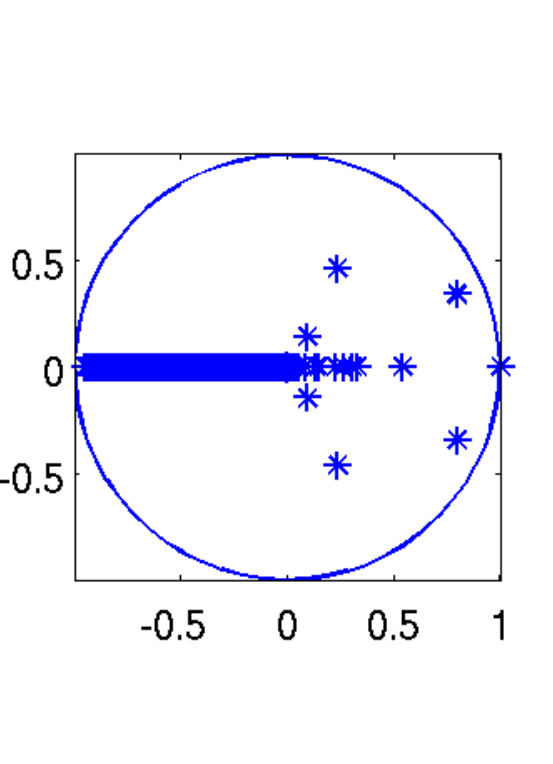}&
\ig[width=0.2\textwidth]{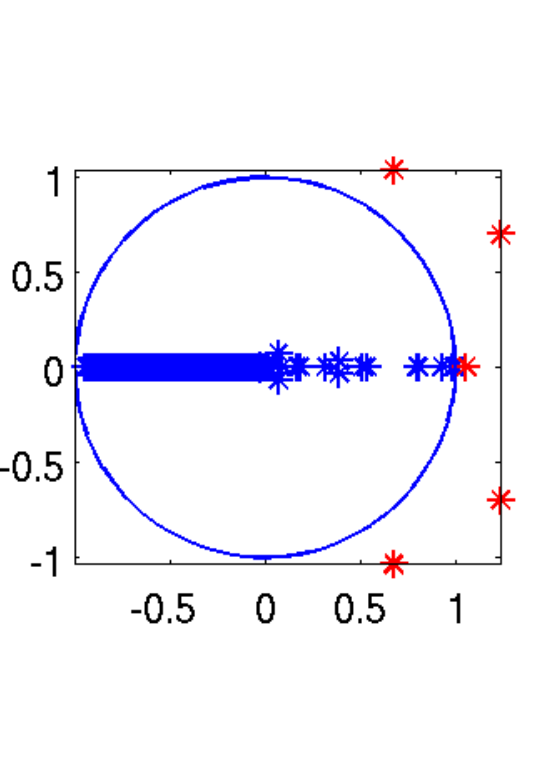}\\
(e) ``error'' time series &(f) initial evolution&(g) transient near h3 &(h) convergence to h1 \\
\ig[width=0.2\textwidth]{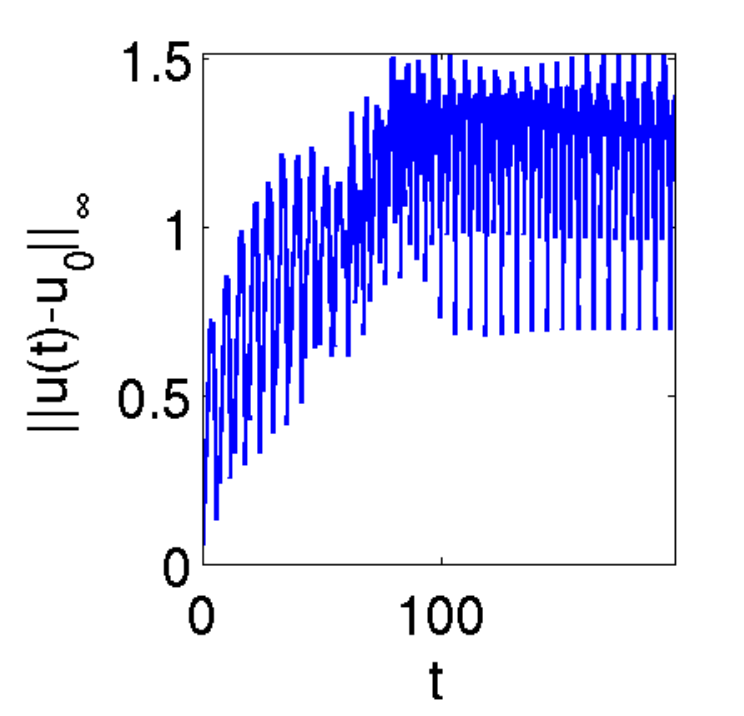}&\ig[width=0.2\textwidth]{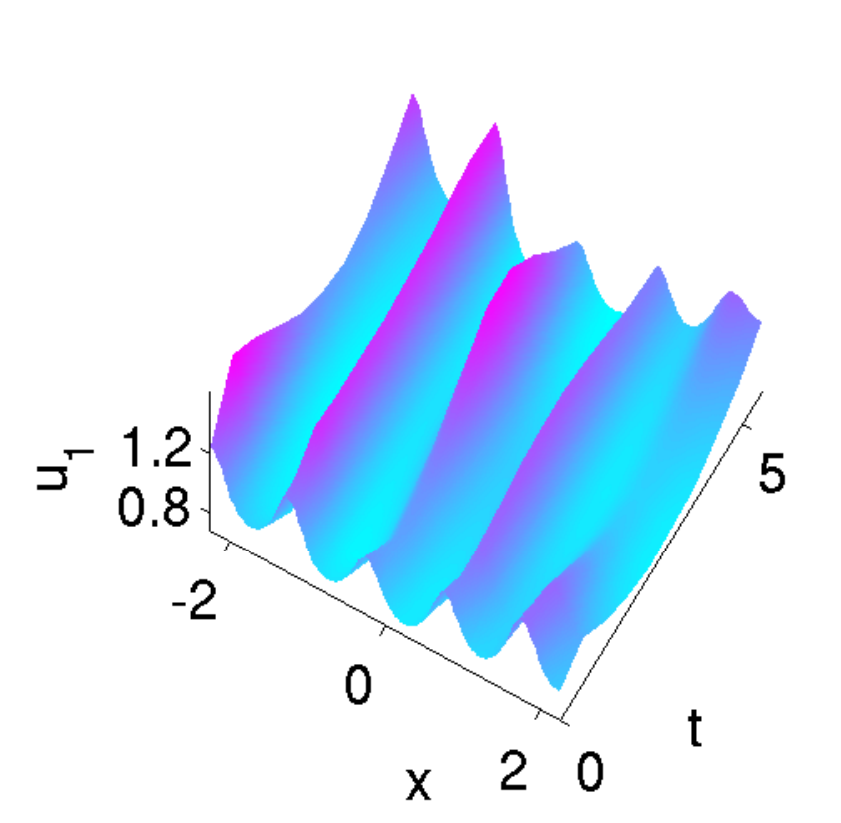}&
\ig[width=0.2\textwidth]{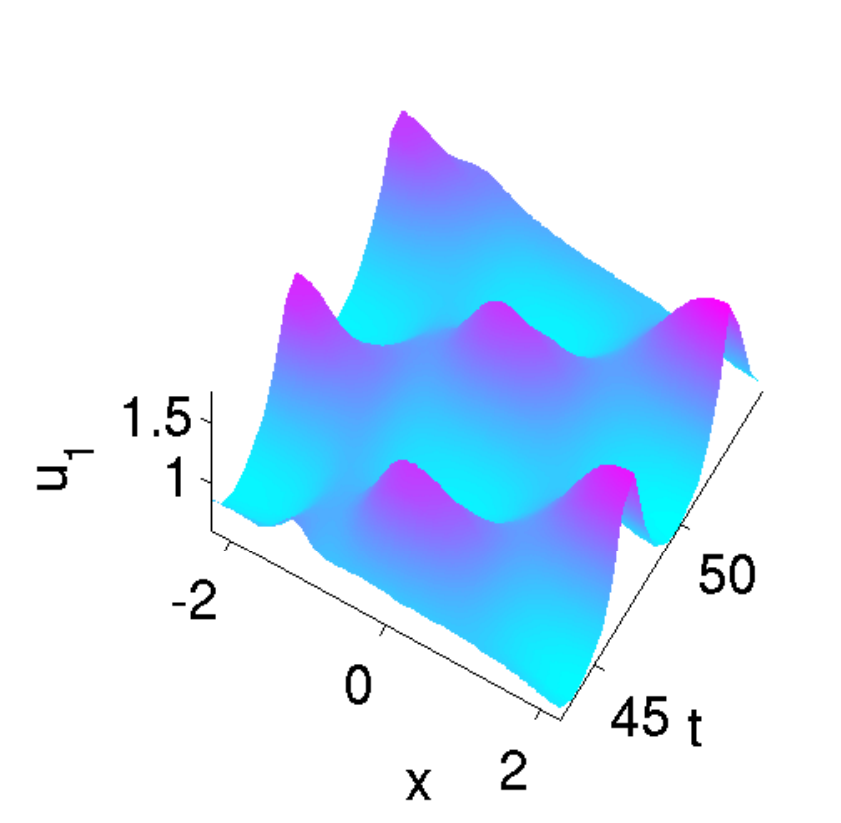}&\ig[width=0.2\textwidth]{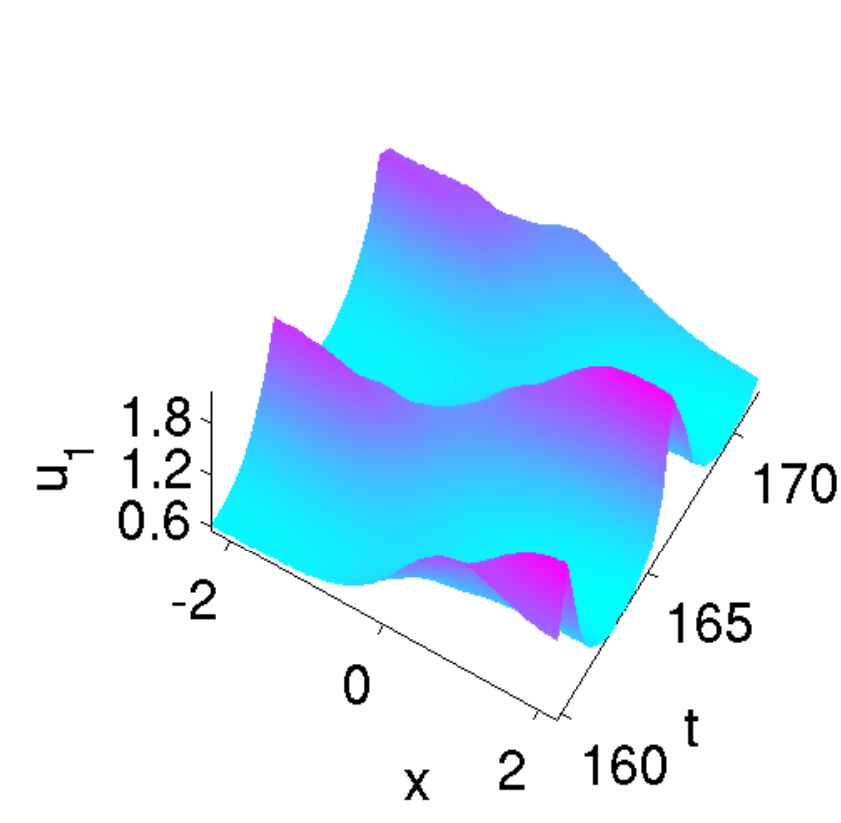}
\end{tabular}}
\ece
\vs{-5mm}
   \caption{{\small (a)-(d) A small sample of Floquet spectra of periodic 
orbits from Fig.~\ref{bf0} (200 largest multipliers computed via {\tt floq}), 
illustrating that a Neimark--Sacker bifurcation 
should be expected near h2/pt5, and similar eigenvalue transitions 
occur on all other Hopf branches except h1. (e)-(h)  Evolution 
of a perturbation of s1h1/pt10. After a rather long transient 
near h3 the solution converges to an orbit on h1. 
  \label{fstab4}}}
\end{figure}

In Fig.~\ref{fstab4}(a)-(d) we give some 
illustration that interesting bifurcations from the Hopf branches should 
occur in \reff{ebru}. It turns out that {\tt h1} is always stable, and that 
(the spatially homogeneous branch) 
h2 is initially unstable with $\ind(u_H)=2$, but close to pt5 on h2 
we find a Neimark--Sacker bifurcation, after which solutions on h2 are 
stable. Similarly, solutions on h3 start with $\ind(u_H)=5$, 
but after a Neimark--Sacker bifurcation, 
and a real  multiplier going through 1 at 
$b\approx 3.35$ we 
find $\ind(u_H)=2$, before $\ind(u_H)$ increases again for larger $b$. 
Also note that there are always many multipliers close to $-1$, but 
we did not find indications for period--doubling bifurcations. 
Finally, in Fig.~\ref{fstab4}(e)--(h) we illustrate the 
evolution of perturbations of s1h1/pt10. After a transient near 
{\tt h3/pt5} (g) the solution converges to a solution from 
the primary Hopf branch {\tt h1} (h), which however itself 
also shows some short wave structure at this relatively large 
distance from bifurcation. 

\begin{figure}[ht]
\bce{\small 
\begin{tabular}{llll}
(a) $\neig=200$&(b) $\neig=300$&(c) $|g|$ from \reff{reso}&(d) \begin{tabular}{l}
$\neig=(3,3)$ with\\ $\om_1=0.9375$\end{tabular}\\
\ig[width=0.19\textwidth, height=35mm]{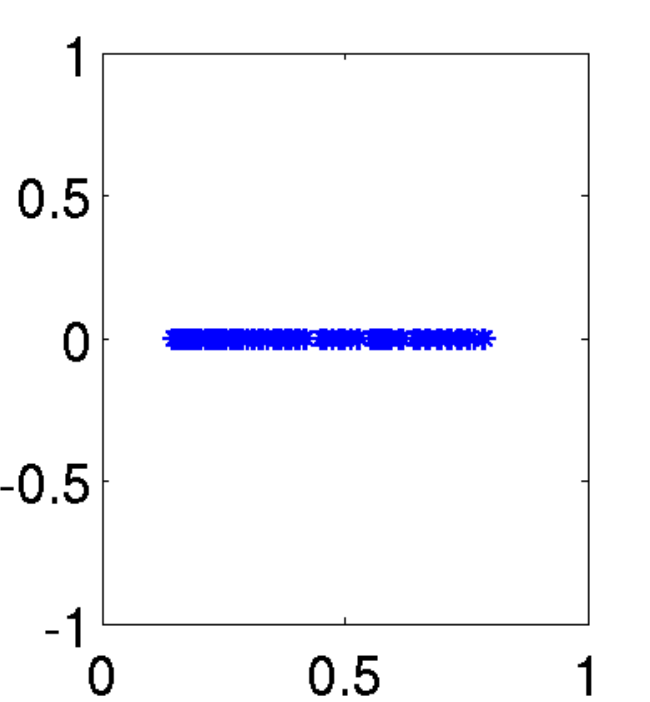}
&\ig[width=0.19\textwidth, height=35mm]{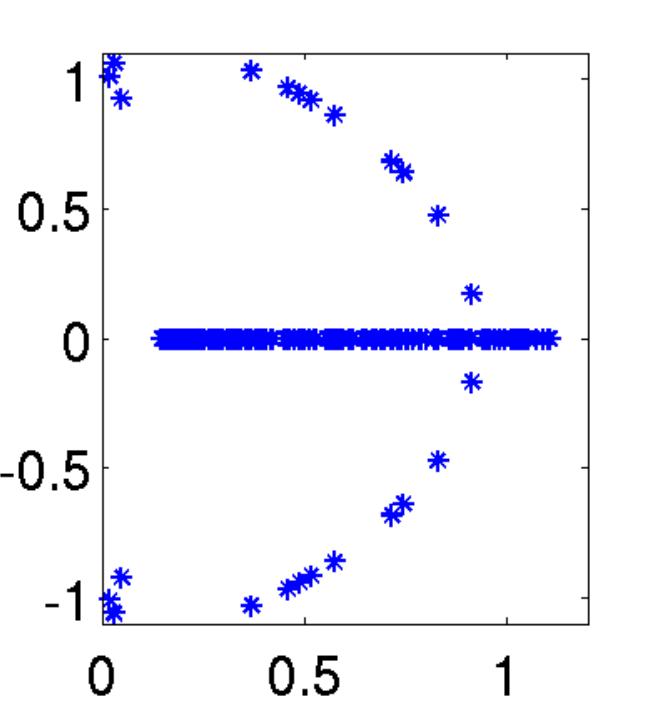}
&\ig[width=0.19\textwidth, height=35mm]{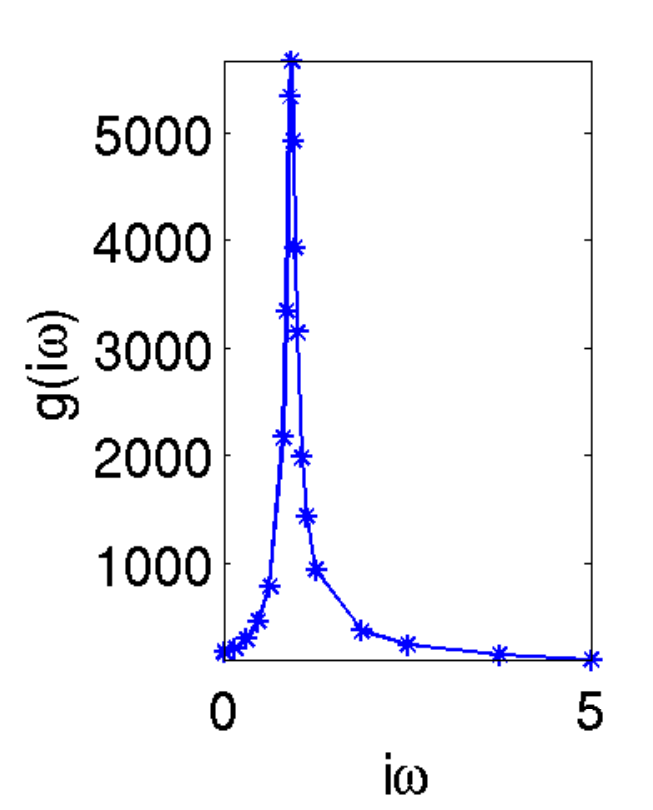}
&\raisebox{5mm}{\ig[width=0.17\textwidth, height=30mm]{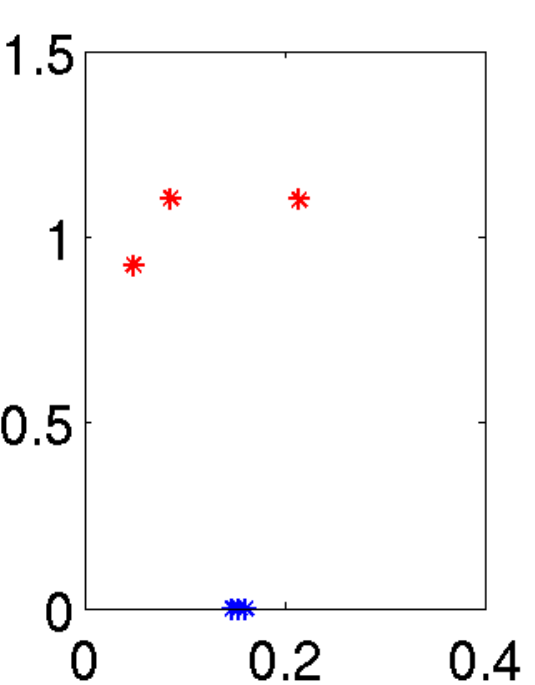}}
\end{tabular}}
\ece 

\vs{-5mm}
   \caption{{\small (a,b) $\neig$ smallest eigenvalues of the linearization 
of \reff{ebru} around 
$U_s$ at $b=2.75$, remaining parameters from \reff{brupar}; 
\heda\  with $\neig=200$  will not detect any Hopf points. 
(c) \reff{reso} yields a guess $\om_1=0.9375$ for the $\om$ value at 
Hopf bifurcation, and then \hedb\ with $\neig=(3, 3)$ is reliable and 
fast: (d) shows the three eigenvalues closest to $0$ in blue, 
and the three eigenvalues closest to $\ri \om_1$ in red. 
  \label{bf1}}}
\end{figure}

In 1D we may still use \heda\  
to detect (and localize) the Hopf bifurcations. 
In 2D this is unfeasible, because even over 
rather small domains we obtain many wave vectors $k=(k_1,k_2)$ with 
modulus $|k|\in(5,8)$, which give leading eigenvalues $\mu_1(k)$ 
with small $\re\mu(k)$ and $\im\mu(k)=0$. 
This is illustrated in Fig.~\ref{bf1}, which shows that for 
$\Om=(-0.5\pi/k_{{\rm TH}}, 0.5\pi/k_{{\rm TH}})^2$ even for $\neig=200$  
(which is quite slow already) we do not 
even see any Hopf eigenvalues, which become ``visible'' at, e.g., $\neig=300$. 
Thus, here we use \hedb\ which runs fast and reliably, even with 
just computing 3 eigenvalues both near 0 and $\om_1$, obtained from \reff{reso}. 

\begin{figure}[ht]
\bce{\small 
\begin{tabular}{p{50mm}p{50mm}p{50mm}}
(a) BD, and $u$ at first HBP on 2ds1 branch&(b) Hopf example plots ($u$)&\hs{-2mm}
(c) Convergence to the primary Hopf branch 2dh1\\
\raisebox{35mm}{\begin{tabular}{l}
\ig[width=0.25\textwidth, height=40mm]{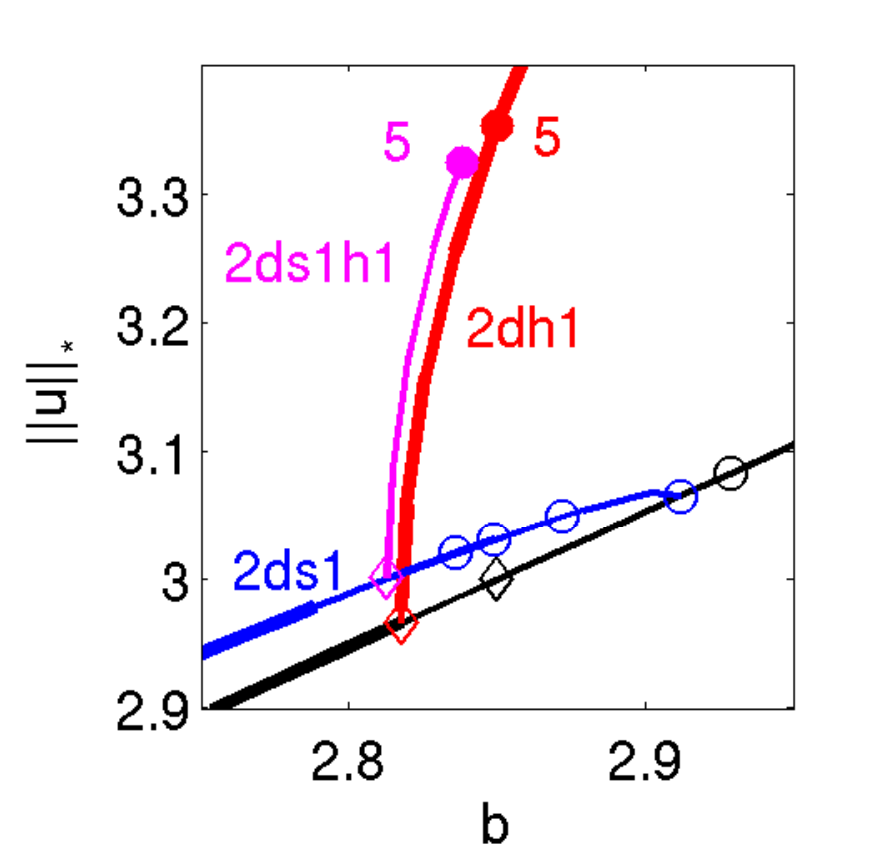}\\
\ig[width=0.24\textwidth,height=16mm]{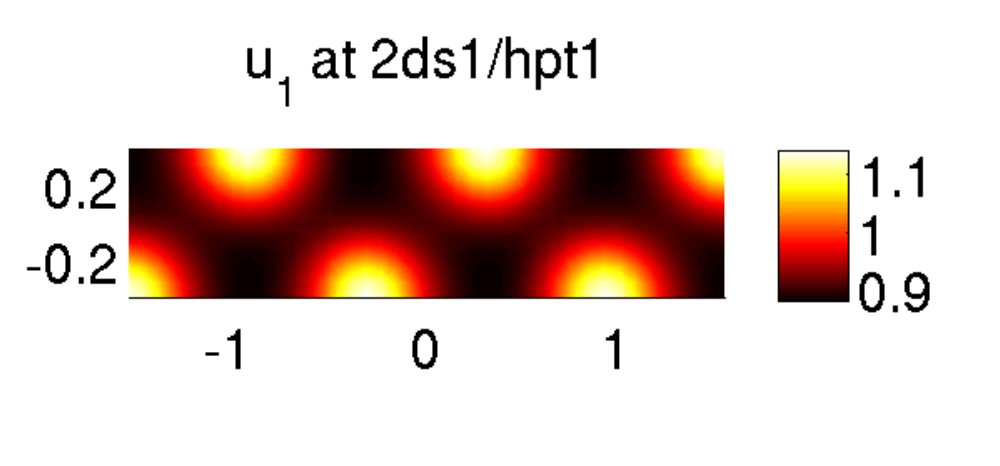}
\end{tabular}}&
\hs{-0mm}\raisebox{35mm}{\begin{tabular}{l}
1) 2dh1/pt5 at $t=0, T/2$\\
\hs{-7mm}\ig[width=0.41\textwidth, height=4mm]{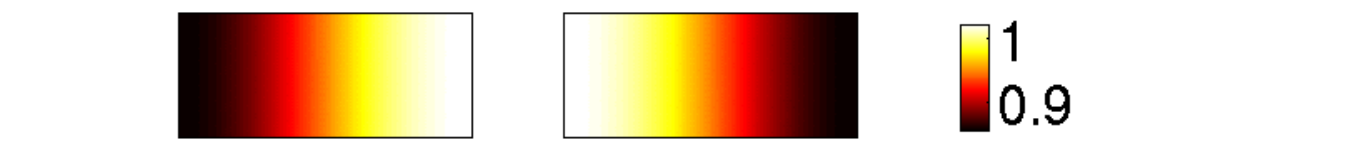}\\[7mm]
2) 2ds1h1/pt5 at $t{=}0,\ldots,3T/4$\\
\hs{-2mm}\ig[width=0.3\textwidth, height=40mm]{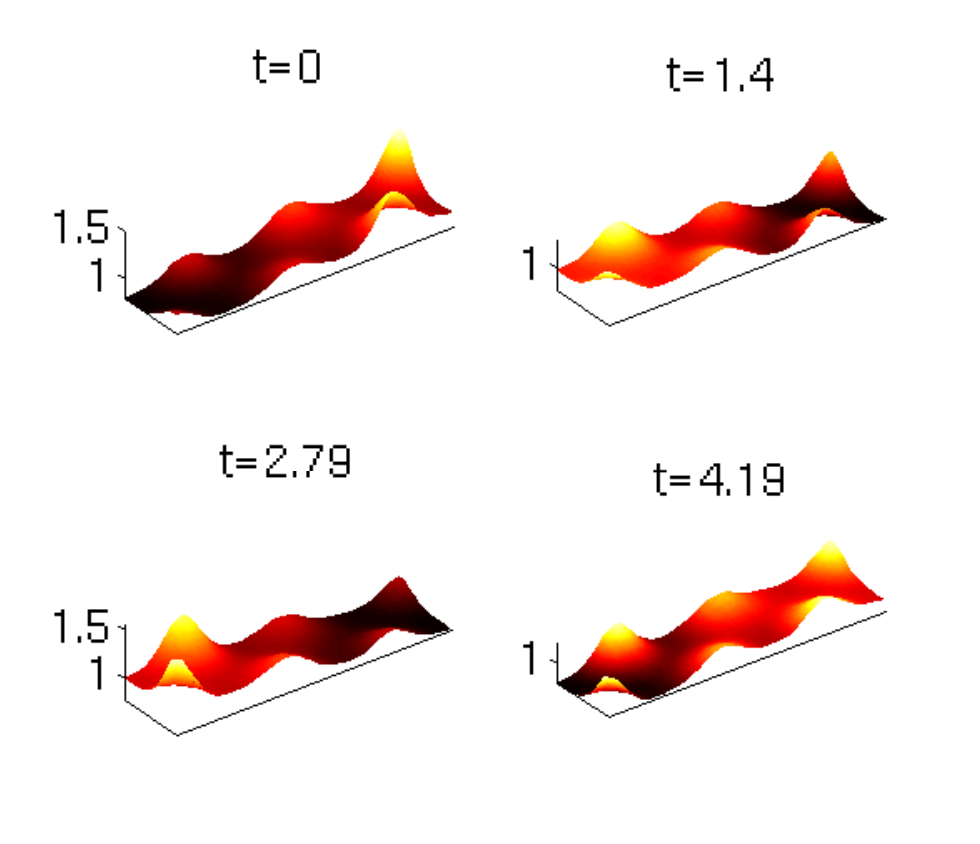}
\end{tabular}}&
\hs{-0mm}\raisebox{10mm}{
\ig[width=0.3\textwidth, height=60mm]{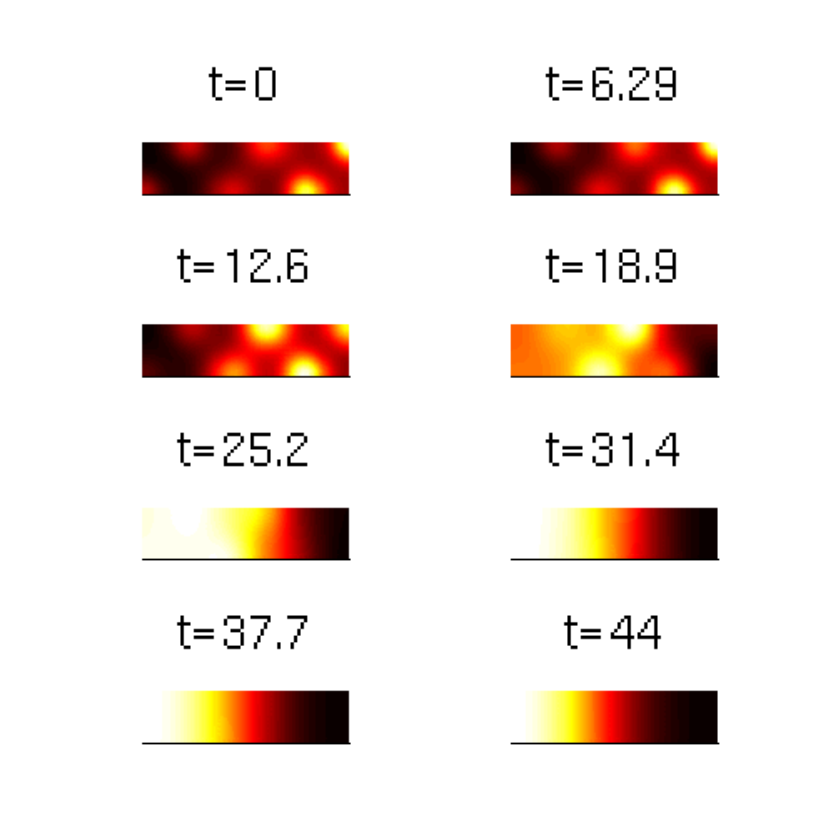}}
\end{tabular}
}
\ece 

\vs{-15mm}
   \caption{{\small (a) Example bifurcations for \reff{ebru} over 
a small 2D domain $\Om=(-\pi/2,\pi/2)\times (-\pi/8,\pi/8)$, and 
example plots of $u$ at  2nd Hopf point on the blue branch. 
(b) Example plots:  solutions on 
primary Hopf branch (1), and 
on the secondary Hopf branch (2) (the amplitude at $t=1.4$ and $t=4.19$ is 
about 0.2). (c) Time integration with 
$u(\cdot,0)$ from 2ds1h1/pt5, snapshots at $0, T, 2T, \ldots, 8T$. 
  \label{bf0b}}}
\end{figure}

Finally, in Fig.~\ref{bf0b} we give an example of just four of the 
many branches which  can be 
obtained for \reff{ebru} in 2D, even over quite small domains. We use 
$\Om=(-l_x,l_x)\times (-l_y,l_y)$, $l_x=\pi/2$, 
$l_y=\pi/8$, where we start wish a mesh of 961 gridpoints, hence 2883 
spatial degrees of freedom. The domain means that admissible wave vectors  are 
$(k_1,k_2)=(n,4m)$, $n,m\in \N_0$. Consequently, no spatial structure 
in $y$ direction occurs in the primary Hopf branches 
(cf.~Fig.~\ref{bf0}b), i.e., the first three are just analogous to 
those in Fig.~\ref{bf0} and occur at $b=2.818$ (with $k=(1,0)$),  
$b=2.859$ (with $k=(0,0)$, i.e., spatially homogeneous, and hence 
$b$ independent of the domain) and $b=3.202$ (with $k=(2,0)$); 
 see (b1) for an example plot on the first Hopf branch. 
The first stationary bifurcation (at $b=2.912$) is now to a spotted 
branch {\tt 2ds1}, and stripe branches analogous to 
{\tt s1} from Fig.~\ref{bf0} bifurcate at larger $b$. 
Interestingly, after some stationary and Hopf bifurcations 
the {\tt 2ds1} branch becomes stable at $b=b_b\approx 2.785$, which illustrates 
that it is often worthwhile to follow unstable branches, 
as they may become stable, or stable branches may bifurcate off.%
\footnote{For continuing this branch we also use a few additional 
features of \pdep\ such as adaptive spatial mesh-refinement and {\tt pmcont}, 
see \cite{hotutb}} 
Figure \ref{bf0b}(b2) shows an example plot 
from the first secondary Hopf branch. This is 
analogous to {\tt s1h1} from Fig.~\ref{bf0}, i.e., the solutions 
look like superpositions of the stationary pattern and solutions on 
the primary Hopf branch {\tt h1}.  Concerning the multipliers we find 
that $\ind(u_H)=0$ on {\tt 2dh1}, and, e.g., $\ind(u_H)=5$ at {\tt 2ds1h1/pt5}, 
where as in 1D (Fig.~\ref{bf1}) there are multipliers suggesting 
Neimark--Sacker bifurcations. 
Figure \ref{bf0b} (c) illustrates the
instability of the spotted Hopf solutions; the spots stay visible  
for about 4 periods, and subsequently 
the solution converges to a periodic orbit from the primary 
Hopf branch, as in Fig.~\ref{fstab4}.

\subsection{A canonical system from optimal control}\label{ocsec}
In \cite{U16,GU17}, \pdep\ has been used to study 
so called canonical steady states and canonical paths 
for infinite time horizon distributed optimal control (OC) problems.  
As an example for such problems with Hopf bifurcations %
 we consider 
\begin{subequations} \label{oc1}
\hual{
    &V(v_0(\cdot))\stackrel!=\max_{k(\cdot,\cdot)}J(v_0(\cdot),k(\cdot,\cdot)), \qquad 
    J(v_0(\cdot),k(\cdot,\cdot)):=\int_0^\infty\er^{-\rho t}
J_{ca}(v(t),k(t)) \dd t,\\
\intertext{where $\ds J_{ca}(v(\cdot,t),k(\cdot,t))=\frac 1{|\Om|}
\int_\Om J_c(v(x,t),k(x,t))\dd x$ 
is the spatially averaged current value function, with }
&J_c(v,k)=pv_1-\beta v_2-C(k) \text{ the local current value,} \quad C(k)=k+\frac 1 {2\ga} k^2, 
 \\
\intertext{
$\rho>0$ is the discount rate (long-term 
investment rate), 
and where the state evolution is}
&\pa_t v_1=-k+d_1\Delta v_1, \quad \pa_t v_2=v_1-\al(v_2)+d_2\Delta v_2, 
\label{woc} 
}
\end{subequations}
with Neumann BC $\pa_{\bf n} v=0$ on $\pa\Om$. 
Here, $v_1=v_1(t,x)$ are the emissions of some firms, 
$v_2=v_2(t,x)$ is the pollution stock, 
 and the control 
$k=k(t,x)$ models the firms' abatement policies. 
In $J_c$, 
$pv_1$ and $\beta v_2$ are the firms' value 
of emissions and costs of pollution, 
$C(k)$ are the costs for abatement, and 
$\al(v_2)=v_2(1-v_2)$ is the recovery function of the  environment. 
The discounted time integral in (\ref{oc1}a) is typical for economic
problems, where ``profits now'' weight more than mid or far 
future profits. 
Finally, the $\max$ in (\ref{oc1}a) runs over all {\em admissible} controls 
$k$; this essentially means that $k\in L^\infty((0,\infty)\times\Om,\R)$, 
and we do not consider active control or state constraints. 

The associated ODE OC problem (no $x$--dependence of $v,k$) 
was set up and analyzed in \cite{TH96,wirl00}; in suitable parameter regimes 
it shows Hopf bifurcations of periodic orbits for the associated 
so called canonical (ODE) system. See also, e.g., \cite{DF91, wirl96, 
grassetal2008} 
for general results about the occurrence of Hopf bifurcations 
and optimal periodic solutions in ODE OC problems. 

Setting $g_1(v,k)=(-k,v_1-\al(v_2))^T$, and introducing  
the co--states (Lagrange multipliers) $$\lam:\Om\times(0,\infty)\ra \R^{2}$$ 
and the 
(local current value) Hamiltonian
$\CH=\CH(v,\lam,k)=J_c(v,k)+\spr{\lam,D\Delta v+g_1(v,k)}$, 
by Pontryagin's Maximum Principle 
for $\tilde{\CH}=
\int_0^\infty \er^{-\rho t} \ov{\CH}(t)\dd t$ with %
$\ov{\CH}(t)=\int_\Om \CH(v(x,t),\lam(x,t),k(x,t))
\dd x$, 
an optimal solution $(v,\lam)$ has to solve the canonical system 
(first order necessary optimality conditions) 
\begin{subequations} \label{cs}
\hual{
\pa_t v&=\pa_\lam\CH=D\Delta v+g_1(v,k), \quad v|_{t=0}=v_0, \\
\pa_t \lam&=\rho\lam-\pa_v\CH=\rho\lam+g_2(v,k)-D\Delta\lam, 
}
\end{subequations}
where $\pa_{\bf n} \lam=0$ on $\pa\Om$. The control $k$ fulfills 
$k{=}\argmax_{\tilde{k}}\CH(v,\lam,\tilde{k})$, and under 
suitable concavity assumptions on $J_c$ and in the absence of 
control constraints is obtained from solving 
$\pa_k\CH(v,\lam,k){=}0$, thus here 
\huga{\label{kform}
k=k(\lam_1)=-(1+\lam_1)/\ga.
}

 Note that \reff{cs} is ill--posed 
as an initial value problem due to the backward diffusion 
in the co--states $\lam$.
Thus it seems unlikely that periodic orbits for \reff{cs} 
can be obtained via shooting methods.   
For convenience we set 
$u(t,\cdot):=(v(t,\cdot), \lam(t,\cdot)): \Om\ra\R^{4}$,  
and write \reff{cs} as 
 \hual{\label{cs2}
&\pa_t u=-G(u):=\CD\Delta u+f(u),
}
where $\ds 
\CD=$diag$(d_1,d_2,-d_1,-d_2)$, 
$\ds f(u)=\biggl(-k, v_1-\al(v_2), 
\rho\lam_1-p-\lam_2, (\rho+\al'(v_2))\lam_2+\beta\biggr)^T$. 
Besides the boundary condition $\pa_{\bf n} u=0$ on $\pa\Om$ 
and the initial 
condition $v|_{t=0}=v_0$ (only) for the states, we have 
the so called intertemporal transversality condition 
\huga{\label{tcond}
\lim_{t\to\infty}\er^{-\rho t}\int_\Om \spr{v,\lam}\dd x=0, 
}
which was already used in the derivation of \reff{cs}. 

A solution $u$ of the canonical system \reff{cs2} is called a 
\emph{canonical path}, 
and a steady state of \reff{cs2} (which automatically fulfills \reff{tcond}) 
is called a \emph{canonical steady state (CSS)}. 
A first step for OC problems of type \reff{oc1} 
is to find canonical steady states and canonical paths connecting 
to some CSS $u^*$. To find such connecting orbits to $u^*$ we may 
choose a cut--off 
time $T_1$ and require that $u(\cdot,T_1)$ is in the stable manifold 
$W_s(u^*)$ of $u^*$, which we approximate by the associated stable eigenspace $E_s(u^*)$. 
If we consider \reff{cs} after spatial discretization, then, 
since we have $n_u/2$ initial conditions, this 
requires that $\dim E_s(u^*)=n_u/2$. Defining the 
defect $d(u^*)$ of a CSS as 
\huga{\label{ddef1}
d(u^*)=\frac{n_u}2-\dim E_s(u^*), 
}
it turns out (see \cite[Appendix A]{GU15}) 
that always $d(u^*)\ge 0$, and we call a $u^*$ with $d(u^*)=0$ a 
saddle--point CSS. 
See \cite{grassetal2008, GU15} for more formal definitions, 
and further comments on the notions of optimal systems,  
the significance of the transversality condition \reff{tcond}, 
and the (mesh-independent) defect $d(u^*)$. 
For a saddle point CSS $u^*$ we can then compute canonical paths 
to $u^*$,  and this has for instance 
been carried out for a vegetation problem in \cite{U16},  with some 
surprising results, including the bifurcation of patterned 
{\em optimal} steady states and optimal paths.  

A natural next step is to search for 
time--periodic solutions $u_H$ of canonical systems, which obviously 
also fulfill 
\reff{tcond}. The natural generalization of \reff{ddef1} is 
\huga{\label{ddef2}
d(u_H)=\ind(u_H)-\frac{n_u}2. 
}
In the (low--dimensional) ODE case, there then exist methods to compute 
connecting orbits to (saddle point) 
periodic orbits $u_H$ with $d(u_H)=0$,   see \cite{BPS01, grassetal2008}, 
which require comprehensive information on the 
Floquet multipliers and the associated eigenspace of $u_H$. 
Our (longer term) aim is to extend these methods to periodic orbits 
of PDE OC systems. 

However, a detailed numerical analysis of \reff{oc1} and 
similar PDE optimal control problems with Hopf bifurcations, and economic 
interpretation of 
the results, will appear elsewhere. Here we only illustrate that 
\bci 
\item Hopf orbits can appear as candidates for optimal solutions 
in OC problems of the form \reff{oc1}, 
\item the computation of multipliers via the periodic 
Schur decomposition (\flb) can yield accurate results, even when 
the computation directly based on 
the product \reff{fl3} (\fla) completely fails. 
\eci 

\begin{figure}[ht]
\bce{\small
\begin{tabular}{p{35mm}p{38mm}p{35mm}p{35mm}}
\mbox{(a) spectrum of} $\pa_u G(u^*)$, 
$\rho=0.5$ &(b) bif.~diagram
&\mbox{(c) time series on h2/pt14 (spat.~homogen.~branch)}\\
\ig[width=0.2\textwidth, height=35mm]{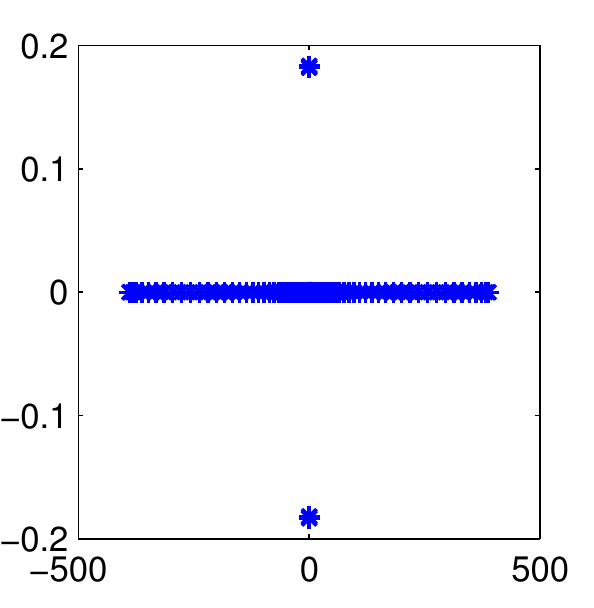}&
\ig[width=0.25\textwidth, height=35mm]{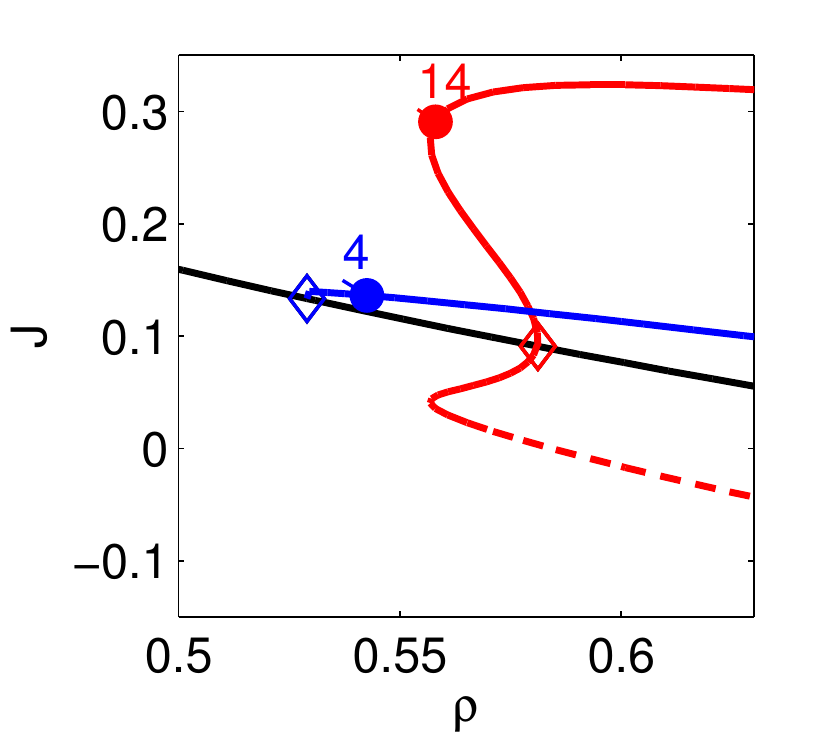}&
\ig[width=0.22\textwidth, height=35mm]{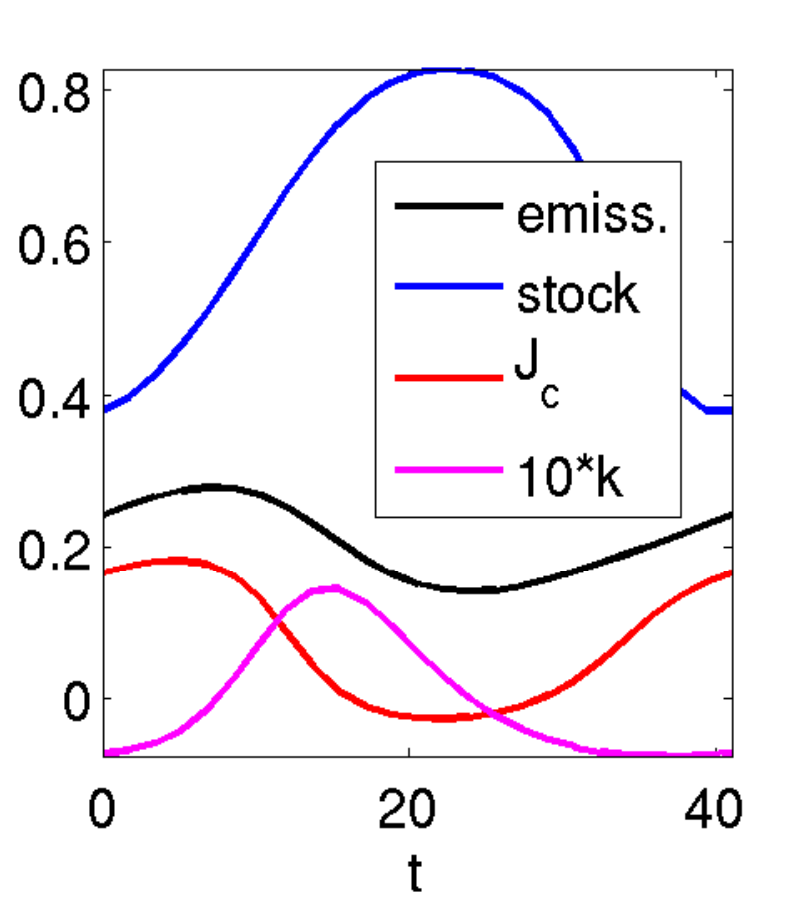}&
\ig[width=0.2\textwidth, height=35mm]{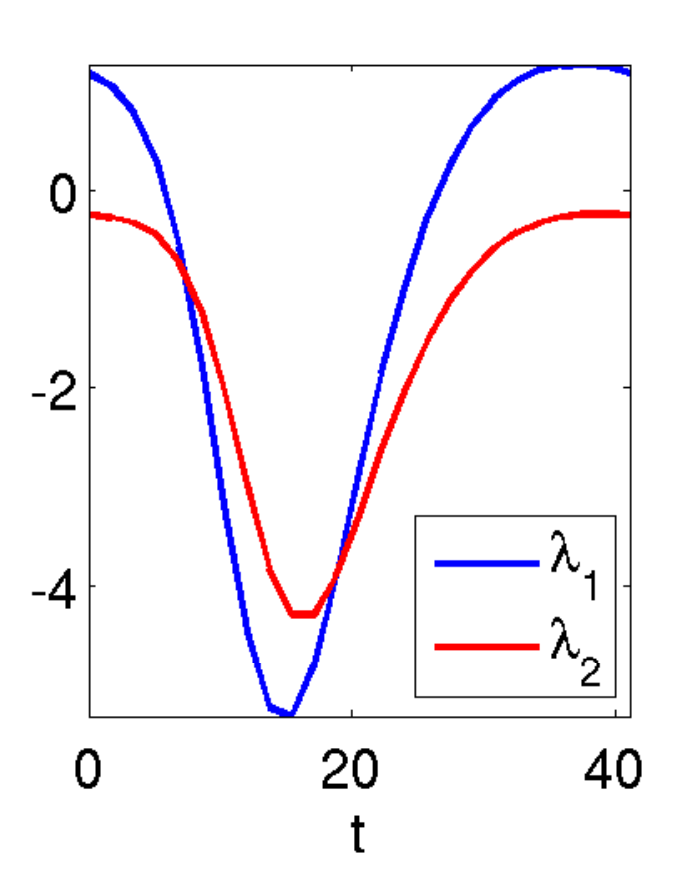}\\
\mbox{(d) example plots at h1/pt4} \\
\ig[width=0.2\textwidth, height=35mm]{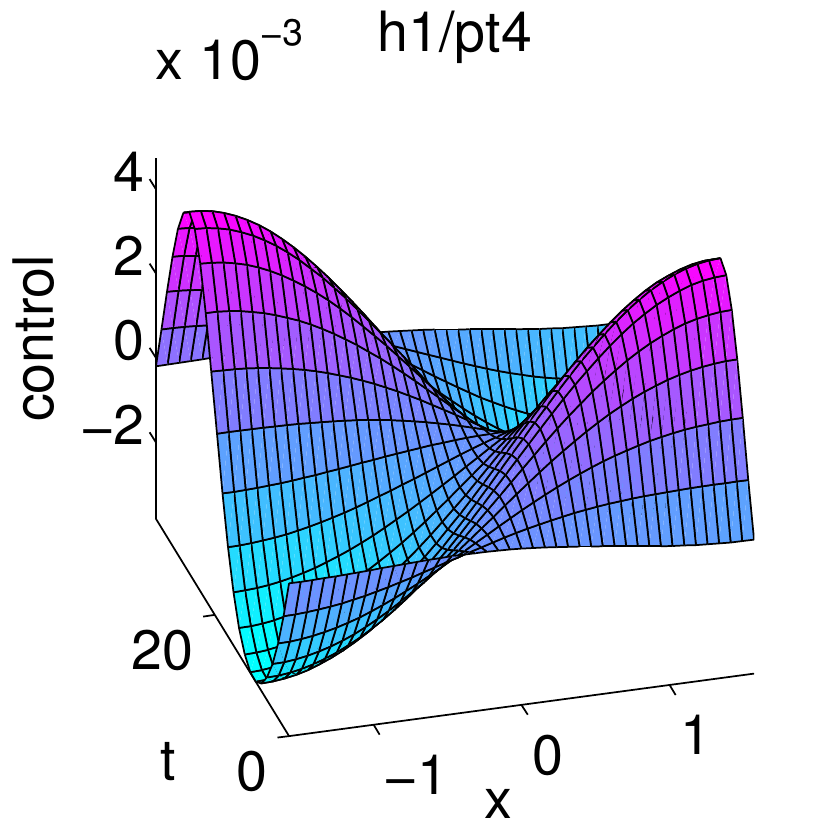}&
\ig[width=0.2\textwidth, height=35mm]{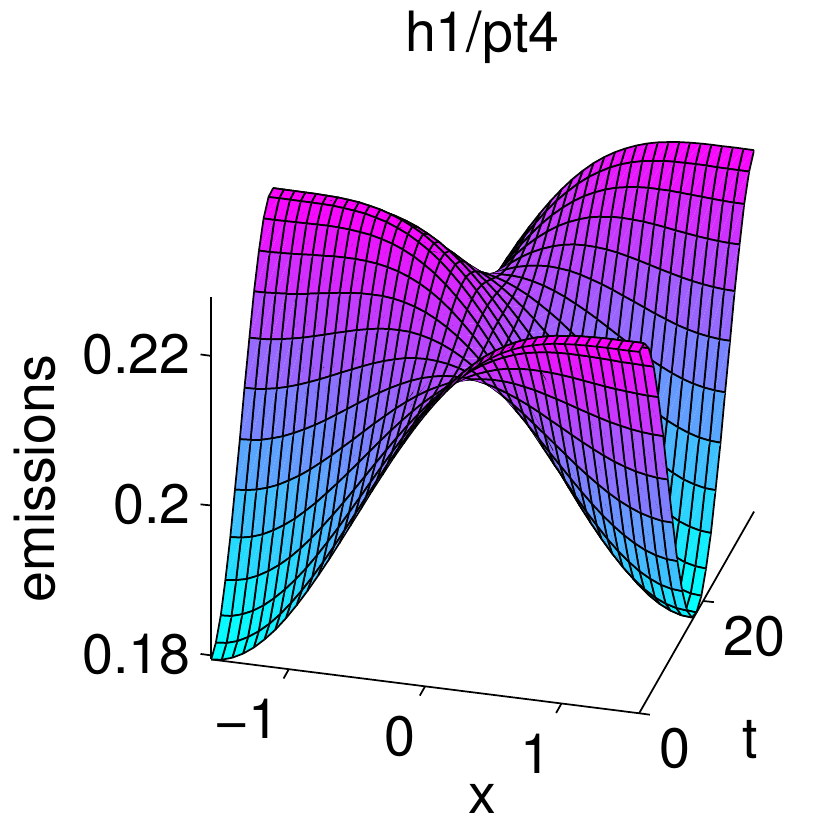}&
\ig[width=0.2\textwidth, height=35mm]{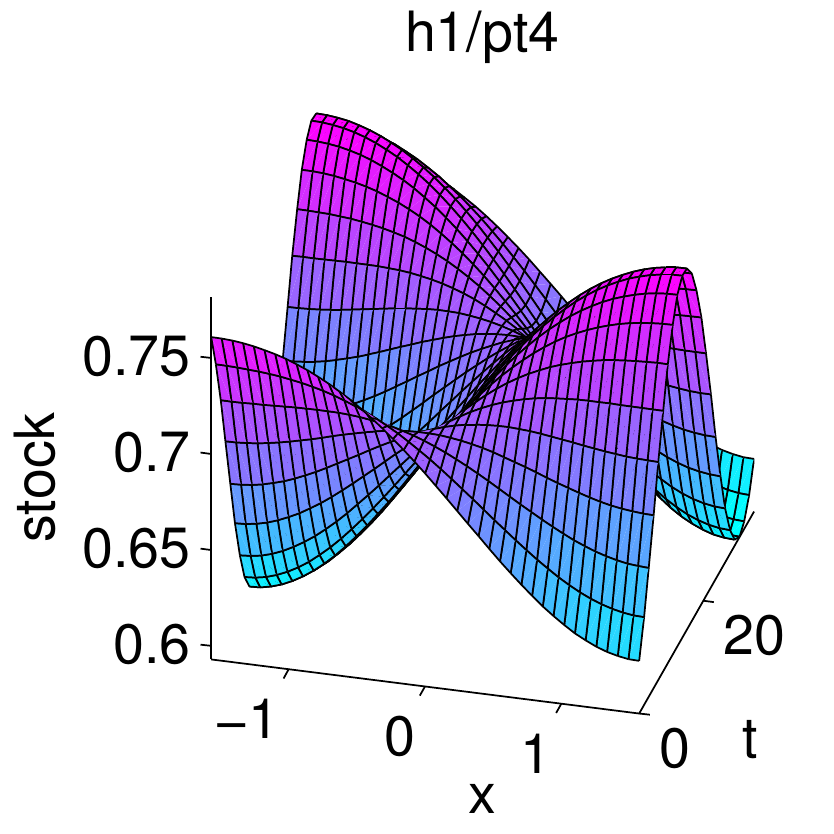}&
\ig[width=0.2\textwidth, height=35mm]{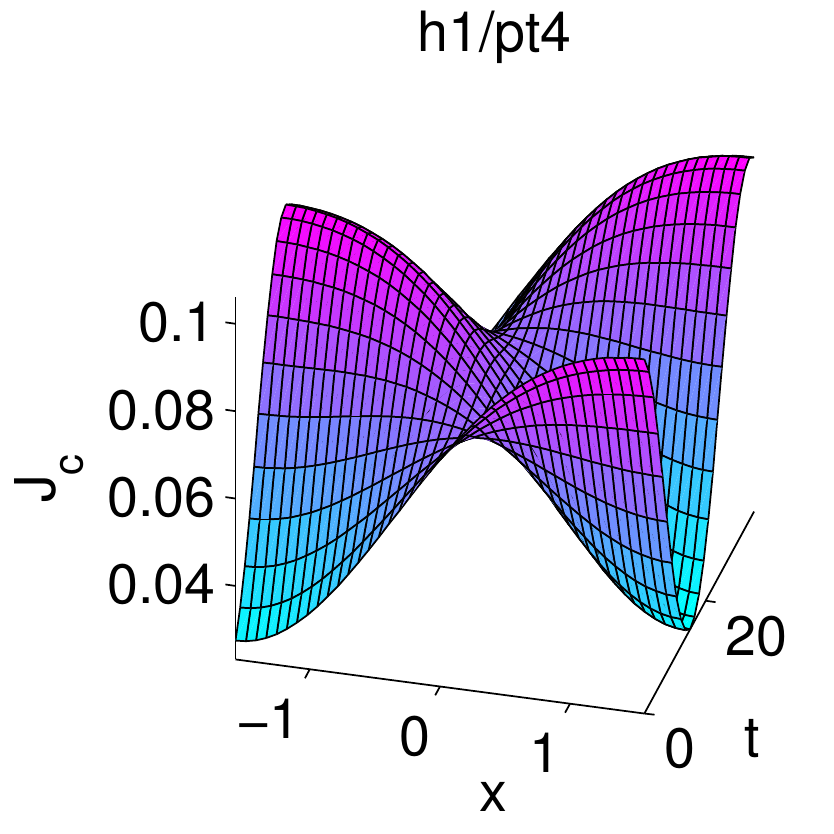}\\
\mbox{(f) the $\frac {n_u}2$ smallest $\ga_j$} at h2/pt4&
\mbox{(g) $|\ga_j|$ for the $\frac {n_u}2$} largest $\ga_j$ at h2/pt4
&\mbox{(h) the $\frac {n_u}2$ smallest $\ga_j$ at h2/pt4 and at h2/pt14}\\
\ig[height=35mm]{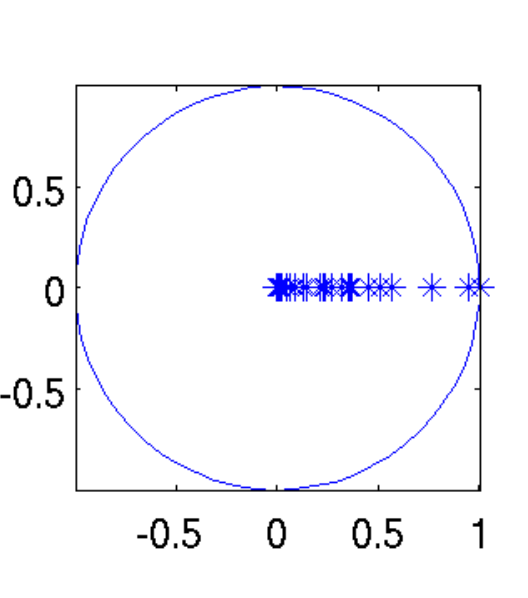}&
\ig[height=33mm]{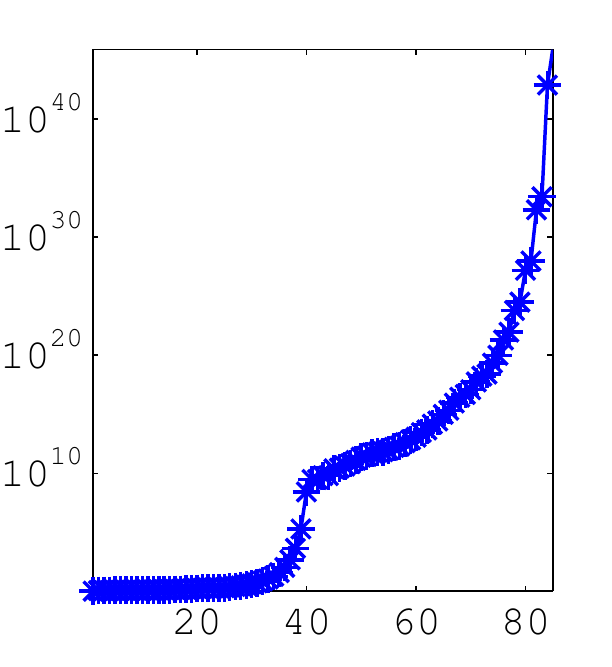}&
\mbox{\ig[height=36mm]{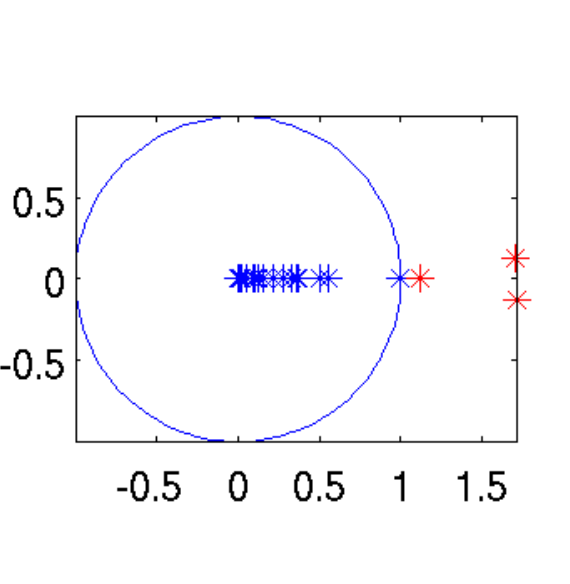}
\raisebox{3mm}{\ig[height=30mm]{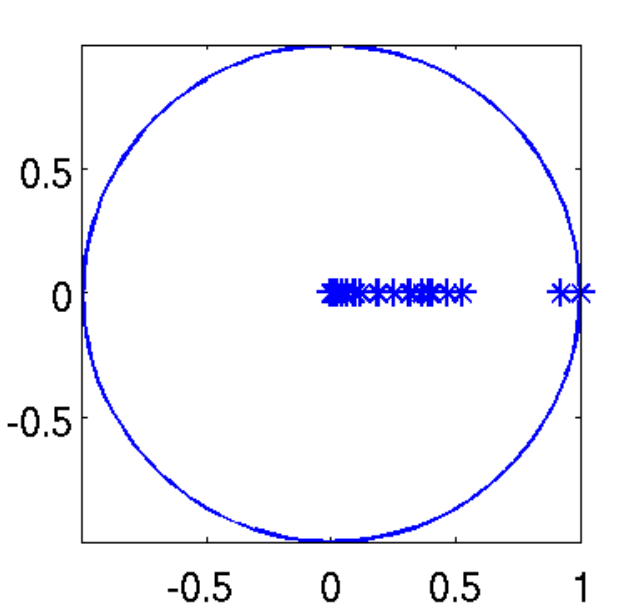}}}
\end{tabular}
}
\ece 

\vs{-5mm}
   \caption{{\small 
(a) full spectrum of the linearization of \reff{cs2} 
around $u^*$ at $\rho=0.5$ on a coarse mesh with $n_p=17$. 
(b) Bifurcation 
diagram, value $J$ over $\rho$. Black: $u^*$; blue: {\tt h1}, 
  red: {\tt h2}, $J(u^H;0)$ (full line) and $J(u^H;T/2)$ (dashed line).  
(c) Time series of a spatially homogeneous solution, including current value  
$J_c$, control $k$, and co--states $\lam_{1,2}$. 
(d,f,g) Example plots and  and multipliers of $u_H$ at {\tt h1/pt4}, which shows 
that $\ind(u_H)=0$. 
(h) multipliers at {\tt h2/pt4}, which shows that $\ind(u_H)=3$ at this point,  
while solutions on h2 become saddles after the fold. 
  \label{ocf1}}}
\end{figure}

For all parameter values, \reff{cs2} has the spatially homogeneous 
CSS 
$$u^*=(z_*(1-z_*),z_*,-1,-(p+\rho)), \quad \text{where}\quad 
z_*=\frac 1 2\left(1+\rho-\frac \beta{p+\rho}\right). 
$$
 We use similar 
parameter ranges as in \cite{wirl00}, namely 
\huga{
(p,\beta,\ga)=(1,0.2,300), \text{ and } \rho\in[0.5,0.65] \text{ as 
a continuation parameter}, 
}
consider \reff{cs2} over 
$\Om=(-\pi/2,\pi/2)$, 
and set the diffusion constants to $d_1=0.001, d_2=0.2$.%
\footnote{The motivation for this choice is to have the first 
(for increasing $\rho$) 
Hopf bifurcation to a spatially patterned branch, and the second to  
a spatially uniform Hopf branch, because the former is 
more interesting. We use that the HBPs for the model \reff{cs2} 
can be analyzed by a simple modification of 
\cite[Appendix A]{wirl00}. We find that 
for branches with spatial wave number $l\in\N$ the 
necessary condition for Hopf bifurcation, $K>0$ from \cite[(A.5)]{wirl00}, 
becomes $K=-(\al'+d_2l^2)(\rho+\al'+d_2l^2)-d_1l^2(\rho+d_1l^2)>0$. 
Since $\al'=\al'(z_*)<0$, a convenient way to first fulfill $K>0$ 
for $l=1$ is to choose $0<d_1\ll d_2<1$, such that for $l=0,1$ the factor 
$\rho+\al'+d_2l^2$ is the crucial one.}
In Figure \ref{ocf1} we give some basic results for \reff{cs2} 
with a coarse spatial discretization of $\Om$ by only $n_p=41$ 
points (and thus $n_u=164$).  
(a) shows the full spectrum of the linearization 
of \reff{cs2} around $u^*$ at $\rho=0.5$, illustrating the ill-posedness 
of \reff{cs2} as an initial value problem. 
(b) shows a basic 
bifurcation diagram. At $\rho=\rho_1\approx 0.53$ there bifurcates a 
Hopf branch {\tt h1} with spatial wave number $l=1$, and at 
$\rho=\rho_2\approx 0.58$ a 
spatially homogeneous ($l=0$) Hopf branch {\tt h2} bifurcates 
subcritically with a fold at $\rho=\rho_f\approx 0.55$. 
(c) shows the pertinent time series on h2/pt14. As should be expected, 
$J_c$ is large when the pollution stock is low and emissions are high, 
and the pollution stock follows the emissions with some delay. 

 Since ultimately we are interested in the 
values $J$ of solutions of \reff{cs2}, in (b) we plot $J$ over 
$\rho$. For the CSS $u^*$ this is simply 
$J(u^*)=\frac 1 \rho J_{c,a}(u^*)$, but for the 
periodic orbits we have to take into account the 
phase, which is free for \reff{cs2}. If $u_H$ is a time periodic solution 
of \reff{cs2}, then, for $\phi\in [0,T)$, we consider 
$$
J(u_H;\phi):=\int_0^\infty \er^{-\rho t}J_{c,a}(u_H(t+\phi))\dd t
=\frac 1 {1-\er^{-\rho T}}\int_0^T
\er^{-\rho t}J_{c,a}(u_H(t+\phi))\dd t, 
$$
which in general may depend on the phase, 
and for {\tt h2} in (c) we plot $J(u_H;\phi)$ for $\phi=0$ (full red line) 
and $\phi=T/2$ (dashed red line).  For the 
spatially periodic branch {\tt h1}, $J_{c,a}(t)$ averages out in $x$ and 
hence $J(u_H;\phi)$ only weakly depends on $\phi$.  
Thus, we first conclude that for $\rho\in(\rho_1,\rho_f)$ the 
spatially patterned periodic orbits from {\tt h1} give the 
highest $J$, while for $\rho\ge \rho_f$ this is obtained from 
{\tt h2} with the correct phase. 
The example plots (c) at {\tt h1/pt4} illustrate how the 
spatio-temporal dependence of $k$ should be chosen, and the resulting 
behaviors of $v$ and $J_c$.  

It remains to compute the defects $d(u^*)$ of the CSS and 
$d(u_H)$ of periodic orbits on the bifurcating branches. 
For $d(u^*)$ we find 
that it starts with $0$ at $\rho=0.5$, and, as expected, 
increases by 2 at each Hopf point. 
On the Hopf branches we always have $n_+\ge n_u/2$ unstable multipliers 
(computed with \flb, which yields 
$\emu<10^{-8}$ for all computations, and hence we trust it), and 
the leading multipliers are very large, i.e., on the order of $10^{40}$, 
even for the coarse space discretization. Thus, we may expect \fla\  
to fail, and indeed it does so completely. 
For instance, calling {\tt floq} to compute all multipliers typically 
returns 10 and larger for the modulus of the {\em smallest} multiplier 
(which from \flb\ is on the order of $10^{-25}$). 

On {\tt h1} we find $d(u_H)=0$ up to {\tt pt4}, 
see (e) for the $n_u/2$ smallest multipliers, and (f) 
for $|\ga_j|$ for the large 
ones, which are mostly real, and $d(u_H)\ge 1$ for larger $\rho$.  
On {\tt h2} we start with $d(u_H)=3$, see (h), but $d(u_H)=0$ after the fold 
until $\rho=\rho_1\approx 0.6$, after which $d(u_H)$ increases 
again by multipliers going though 1. 
Since on {\tt h1} we have that $J(u_H)$ is 
larger than $J(u^*)$, and since $u_H$ is a saddle point up to {\tt pt4}, 
we expect that these $u_H$ are at least locally optimal, 
and similarly we expect $u_H$ from {\tt h2} after the fold until 
$\rho_1$ to be 
locally, and probably globally, optimal. However, as already said, 
for definite answers and, e.g., to characterize the domains of 
attractions, we need to compute canonical paths connecting to these 
periodic orbits, and this will be studied elsewhere.

\section{Summary and outlook}\label{dsec}
With the {\tt hopf} library of \pdep\ we provide basic functionality for Hopf 
bifurcations and periodic orbit continuation 
for the class \reff{tform} of PDEs  over 1D, 2D and 3D domains. 
The user interfaces reuse the standard 
\pdep\ setup, and no new user functions are necessary. 
For the detection of Hopf points we check for eigenvalues crossing 
the imaginary axis near guesses $\ri \om_j$, where the $\om_j$ 
can either be set by the user (if such a priori information is 
available), or can be estimated using 
the function $g$ from \reff{reso}. 
An initial guess for a bifurcating periodic orbit is then 
obtained from the normal form \reff{honf}, and the 
continuation of the periodic orbits is based on modifications 
of routines from TOM \cite{MT04}. 

Floquet multipliers of 
periodic orbits can be computed from the monodromy 
matrix $\CM$ \reff{fl3} (\fla), 
or via a periodic Schur decomposition of the block matrices of $\CM$ (\flb). 
The former is suitable for dissipative systems, and computes 
a user chosen number of largest multipliers of %
$\CM$. This definitely fails for problems 
of the type considered in \S\ref{ocsec}, and in general we recommend 
to monitor $\emu=|\ga_1-1|$ to detect further possible inaccuracies. 
The periodic Schur decomposition is expensive, 
but has distinct advantages:  
It can be used to efficiently compute eigenspaces 
at all time--slices and hence bifurcation information in case 
of critical multipliers, and, presently most importantly for 
us, it accurately (measured by  $\emu$) computes the multipliers 
also for ill posed evolution problems.

We tested our algorithms on four example problems, 
where we believe that the second, third and fourth are close to interesting 
research problems. For instance, the numerical results on \reff{spir1},\reff{spirbc}  
seem to be the first on bifurcation of spiral waves 
out of zero over a bounded domain, 
in a reaction diffusion system without very special boundary conditions. 
Further interesting problems will be, e.g., the bifurcation {\em from} Hopf branches in this model, and in \reff{ebru}. Thus, 
as one next step we plan to  implement 
the necessary localization and branch switching routines, for which the 
cGL equation \reff{cAC0} will again provide a good test case. In \S\ref{ocsec} we give a (very basic) illustration 
of the widely unexplored field of Hopf bifurcations and time 
periodic orbits in infinite time horizon distributed optimal 
control PDE problems. For this, as a next step we plan to implement routines 
to compute canonical paths connecting to periodic orbits. 

Finally, another interesting field are Hopf bifurcations from traveling waves, or more generally in systems with continuous symmetries, see Remark \ref{qrem}. 
Examples how to treat these in \pdep, based on the setup presented here 
with minor additions, are given in \cite[\S4]{symtut} and \cite[\S5]{hotut}.

\renewcommand{\refname}{References}
\renewcommand{\arraystretch}{1.05}\renewcommand{\baselinestretch}{1}
\small
\bibliographystyle{alpha}
%
\newcommand{\etalchar}[1]{$^{#1}$}

\end{document}